\documentclass[11pt]{article}

\usepackage{amsmath}
\usepackage{amssymb}
\usepackage{latexsym}
\usepackage{amsmath, amsfonts,amssymb, amsthm, euscript,makeidx,color,mathrsfs}
\usepackage{esint}
\usepackage{cleveref}
\usepackage{enumitem}
\usepackage{tikz}

\usepackage{pgfplots,subfig,graphicx}
\pgfplotsset{compat=1.9}

\usepackage{float}

\usetikzlibrary {positioning}
\usepackage{pgf}
\usepackage{pgfplots}

\usepackage[numbers]{natbib}

\newtheorem{assumption}{Assumption}
\def\qed{ \ \vrule width.2cm height.2cm depth0cm\smallskip}

\newcommand{\hP}{\hat\dbP}

\newcommand{\ba}{\begin{array}}
\newcommand{\ea}{\end{array}}
\newcommand{\be}{\begin{equation}}
\newcommand{\ee}{\end{equation}}
\newcommand{\bea}{\begin{eqnarray}}
\newcommand{\eea}{\end{eqnarray}}
\newcommand{\beaa}{\begin{eqnarray*}}
\newcommand{\eeaa}{\end{eqnarray*}}

\def\neg{\negthinspace}

\def\o{\omega}

\def\O{\Omega}

\def\O{\Omega}
\def\cA{{\cal A}}

\def\cF{{\cal F}}
\def\cG{{\cal G}}

\def\cJ{{\cal J}}
\def\cK{{\cal K}}

\def\cN{{\cal N}}

\def\cR{{\cal R}}

\def\cZ{{\cal Z}}
\def\hB{\mathbb{B}}
\def\hC{\mathbb{C}}

\def\hE{\mathbb{E}}
\def\hF{\mathbb{F}}
\def\hG{\mathbb{G}}

\def\hL{\mathbb{L}}

\def\hN{\mathbb{N}}

\def\hP{\mathbb{P}}
\def\hQ{\mathbb{Q}}
\def\hR{\mathbb{R}}

\def\hX{\mathbb{X}}

\def\sA{\mathscr{A}}
\def\sB{\mathscr{B}}
\def\sC{\mathscr{C}}

\def\sG{\mathscr{G}}

\def\sK{\mathscr{K}}
\def\scL{\mathscr{L}}

\def\sN{\mathscr{N}}

\def\sR{\mathscr{R}}

\def\sW{\mathscr{W}}

\def\no{\noindent}

\def\ss{\smallskip}

\def\bs{\bigskip}
\def\q{\quad}

\def\cd{\cdot}

\def\lan{\langle}
\def\ran{\rangle}

\def\qed{ \hfill \vrule width.25cm height.25cm depth0cm\smallskip}

\newcommand{\basa}{\begin{assumption}}
\newcommand{\easa}{\end{assumption}}

\newcommand{\ol}{\overline}

\newcommand{\bas}{\begin{assum}}
\newcommand{\eas}{\end{assum}}

\def\lan{\mathop{\langle}}
\def\ran{\mathop{\rangle}}

\def\Leb{\mathop{\rm Leb}}

 \def\cd{\cdot}

\def\as{\hbox{\rm-a.s.{ }}}

\def\ol{\overline}

\def\cl{\mbox{\rm cl}}

\def\MS{\mbox{\rm MS}}
\def\TC{\mbox{\rm TC}}

\def\1{{\bf 1}}

\def\:{\!:\!}

\DeclareMathOperator{\aff}{aff}
\DeclareMathOperator{\dec}{dec}

\DeclareMathOperator{\co}{conv}
\DeclareMathOperator{\cone}{cone}
\DeclareMathOperator{\interior}{int}

\DeclareMathOperator{\bd}{bd}

at 9pt

\newcommand{\of}[1]{\ensuremath{\left( #1 \right)}}
\newcommand{\cb}[1]{\ensuremath{ \left\{ #1 \right\} }}
\newcommand{\sqb}[1]{\ensuremath{ \left[ #1 \right] }}
\newcommand{\norm}[1]{\ensuremath{ \left\Vert #1 \right\Vert }}
\newcommand{\ip}[1]{\ensuremath{ \left\langle #1 \right\rangle }}

\def\prehp(#1,#2){\ensuremath{  #1 \cdot #2 }}

\let\abs=\envert

\begin{document}

\newtheorem{thm}{Theorem}[section]
\newtheorem{lem}[thm]{Lemma}
\newtheorem{cor}[thm]{Corollary}
\newtheorem{prop}[thm]{Proposition}
\newtheorem{rem}[thm]{Remark}
\newtheorem{eg}[thm]{Example}
\newtheorem{defn}[thm]{Definition}
\newtheorem{assum}[thm]{Assumption}

\renewcommand {\theequation}{\arabic{section}.\arabic{equation}}
\def\thesection{\arabic{section}}

\title{\bf Path-Regularity and Martingale Properties of \\ Set-Valued Stochastic Integrals}

\author{
 \c{C}a\u{g}{\i}n Ararat\thanks{Department of Industrial Engineering, Bilkent University, Ankara,~06800, Turkey. Email: cararat@bilkent.edu.tr. This author is supported by T{\"{U}}B\.{I}TAK~2219 Program and by the Fulbright Scholar Program of the U.S. Department of State, sponsored by the Turkish Fulbright Commission. This work was partly completed while the author was visiting University of Southern California, whose hospitality is greatly appreciated.}
 ~ and ~ Jin Ma\thanks{Department of
Mathematics, University of Southern California, Los Angeles, CA, 90089, USA.
Email: jinma@usc.edu. This author is supported in part by US NSF grants DMS\#1908665 and \#2205972.}
}

%\date{}
\date{\today}
\maketitle

\begin{abstract}
	In this paper we study the path-regularity and martingale properties of the set-valued stochastic integrals defined in our previous work \cite{amw}. Such integrals have some fundamental differences from the well-known Aumann-It\^o stochastic integrals, and are much better suitable for representing set-valued martingales, whence potentially useful in the study of {\it set-valued backward stochastic differential equations}. However, similar to the Aumann-It\^o integral, the new integral is only a set-valued submartingale  in general, and there is very limited knowledge about the path regularity of the related indefinite integral, much less the sufficient conditions under which the integral is a true martingale. In this paper, we first establish  the existence of right- and left-continuous modifications of set-valued submartingales in continuous time, and apply the  results to set-valued stochastic integrals. Moreover, we show that a set-valued stochastic integral yields a martingale if and only if the set of terminal values of the stochastic integrals associated to the integrand is closed and decomposable. Finally, as a particular example, we study the set-valued martingale in the form of the conditional expectation of a set-valued random variable. We show that when the random variable is a convex random polytope, the conditional expectation of a vertex stays as a vertex of the set-valued conditional expectation if and only if the random polytope has a deterministic normal fan.
\end{abstract}

\vfill \bs

\no

{\bf Keywords.} \rm Set-valued stochastic integral, set-valued martingale, martingale representation, path-regularity, random polytope, normal fan.

\bs

\no{\it 2020 AMS Mathematics subject classification:} 60H05, 60G44, 28B20, 47H04.

% 60H05: stochastic integrals
% 60H10: stochastic ordinary differential equations
% 60G44: martingales with continuous parameter
% 28B20: set-valued set functions and measures; integration of set-valued functions; measurable selections
% 47H04: set-valued operators

\eject

\section{Introduction}\label{sec:intro}

The theory of set-valued stochastic integrals, often referred to as the \emph{Aumann-It\^{o}  integral}, as well as the associated set-valued stochastic analysis have been well-established in the literature. We refer the reader to the well-known books \cite{K, Kis2020} and the references cited therein for the history and basic knowledge of the theory; see also \cite{viability,cascosmolchanov,zachsurvey, Molchanov,econometrics} for various studies of \emph{random sets, set-valued random variables and stochastic processes}, and their applications in control theory, economics, and finance.

A simple but intriguing example of a set-valued stochastic process is the conditional expectation (in the sense of \cite{hiai-umegaki}) of a set-valued random variable with respect to a given filtration. Under sufficient integrability assumptions on the set-valued random variable, this process is well-defined and it is a \emph{set-valued martingale} (cf. \cite{hess,hess1999,hiai-umegaki}). Indeed, in a finite horizon setting, all set-valued martingales are of this form as in the case of vector-valued martingales. However, to date, there does not seem to have been any detailed study of the structure and characterization of this simplest set-valued martingale. For example, if the set-valued random variable is in the form of a random polytope with random vertices, then would this martingale be generated by the vector-valued martingales obtained by taking the conditional expectations of the random vertices?

It turns out that this is not a trivial question to answer using the currently available literature on set-valued stochastic analysis. It is well-known that in the case when the underlying filtration is ``Brownian", the celebrated martingale representation theorem provides a stochastic integral representation for every square-integrable martingale. The set-valued analog of this problem, often in terms of the \emph{Aumann-It\^{o} stochastic integral} and/or its generalized version, was studied in \cite{Kis1997,Kis2020}. However, as we shall see below, the notion of Aumann-It\^{o} stochastic integral has some serious deficiencies from the martingale representation point of view.

To make our point more clearly, let us state some facts regarding the Aumann-It\^{o} stochastic integral (cf. e.g., \cite{Kis2020}) that we have pointed out in our recent work \cite{amw}. Let us denote the \emph{indefinite} Aumann-It\^{o} integral by $I^0_t(\cZ)= \int_0^t \cZ dB=\int_0^T\cZ \1_{[0,t]} dB$, $0\le  t\le T$, where $\cZ$ is some appropriately measurable set-valued process (cf. \cite{aumann,Kis1997} or \cite{Kis2020}). We should first note that, while this integral has some similarities with the usual vector-valued It\^{o} integral, it does not have the $\hL^2$-isometry property. In fact, it is always \emph{stochastically unbounded} in the $\hL^2$-sense unless $\cZ$ is a singleton. Consequently, in order for an Aumann-It\^{o} integral to be a non-singleton stochastically bounded set-valued martingale, the so-called  \emph{generalized Aumann-It\^{o} stochastic integral} is introduced (cf. \cite{Kis2020}) so that the integrand $\cZ $ is allowed to be a set of progressively measurable processes that is not necessarily \emph{decomposable} (see Section~\ref{sec:svfunctions} for details). However, besides the immediate drawback of lacking the elementary \emph{temporal additivity} property (i.e., in general, one only has $ I^0_T( \cZ) \subseteq I^0_t( \cZ)+ I^t_T( \cZ)$ and the equality holds when $\cZ$ is decomposable), a much more serious issue is that it is strongly restricted by the so-called ``singleton test" (see, e.g, \cite{hess1999, ZhangYano}). That is, if the process $(I^0_t(\cZ))_{t\geq 0}$ is a set-valued martingale with $I_0^0(\cZ) = \{0\}$, then it must be a singleton(!). Such a result essentially nullifies the existing set-valued martingale representation theorem in \cite{MgRT}, and eliminates the possibility of building a theory of truly set-valued backward SDEs based on the existing framework of stochastic analysis.

In our recent work \cite{amw}, we proposed a further generalization of the Aumann-It\^{o} integral so that it allows for the (non-singleton) initial values and can be used to represent a truly set-valued martingale in a temporally pointwise manner. However, such a representation is still a far cry from the desired martingale representation theorem for two important reasons. First, due to the lack of the much needed ``path-regularity" of the set-valued stochastic integrals, the representing integral is unique only up to modifications. Second, since the generalized  set-valued stochastic integral only produces set-valued submartingales similar to the Aumann-It\^{o} integral, the martingale representation theorem provides only an injection from the collection of set-valued martingales to that of set-valued stochastic integrals, a much weaker result than its vector-valued counterpart.

In this paper, we are aiming at an in-depth study of set-valued (sub)martingales and stochastic integrals through the following questions that may be of independent interests in their own rights:

(i) Does a set-valued stochastic integral always admit a modification with certain path regularity in the Hausdorff sense?

(ii) When will a set-valued stochastic integral be a martingale?

(iii) How to characterize the martingales that are  in the form of conditional expectations of terminal values by the set-valued stochastic integrals, which are submartingales in general?

To elaborate more on the motivations of these questions, we first note that unlike the vector-valued stochastic integrals the path regularity of the set-valued integral is quite non-trivial, due mainly to the nature of ``decomposability". In fact, the set-valued literature lacks such results  except for some special cases. On the other hand, path-regularity results would facilitate the discussion of various types of joint measurability in stochastic analysis and the study of SDEs.

To answer question (i), we first study the path-regularity issue in a more general setting of continuous-time submartingales under an arbitrary filtration. We shall alternately use two ways of ``scalarizing" a set-valued submartingale: its \emph{``set norm" process} (i.e., the maximum norm of a point in the submartingale at each time) and its \emph{support function process} (i.e., the maximum value of a linear function over the submartingale at each time). Since both processes are real-valued submartingales, the usual Martingale Convergence Theorem can then be used to establish some fundamental path regularity results for set-valued submartingales. 
%As a fundamental tool and a subject of independent interest, we also prove a convergence result for set-valued reverse submartingales in discrete time. Then, 
Under the usual conditions on the filtration, we shall argue that a set-valued submartingale admits a right-continuous modification if and only if its expectation function is right-continuous. Surprisingly, we can also show that, under an augmented and left-continuous filtration, the same result holds for the ``left-continuous" counterpart, which is not the case for real-valued submartingales since a submartingale restricted to a compact interval may not be uniformly integrable in general. This seemingly paradoxical contrast between the theory of real-valued and set-valued submartingales underscores a crucial observation:  \emph{a real-valued submartingale does not form a (singleton) set-valued submartingale}(!) as the former is defined via the usual ordering of real numbers, whereas the latter uses the ``subset" relation between sets.

To answer question (ii), we argue that the only set-valued stochastic integrals that yield martingales are those for which the terminal values of the ``trajectory integrals" (the vector-valued stochastic integrals of the selectors of the integrand) form already a decomposable closed set. That is, the lack of martingale property for a set-valued stochastic integral is due to lack of decomposability. To the best of our knowledge, such a result, albeit conceivable, is novel.

Finally, as a continuation of (ii), question (iii) is motivated by the following simple issue: given a random polytope at the terminal time, being a martingale, its conditional expectation process can be expressed as a stochastic integral. On the other hand, we may also consider the vector-valued conditional expectations of its vertices. The question is whether the conditional expectations of the vertices remains vertices of the conditional expectation of the random polytope at all times. We shall first study the case in a greater generality, that is, the given set-valued random variable is convex but not necessarily polyhedral. We show that the conditional expectation of the given random set is the closed convex hull of the conditional expectations of its measurable selections if and only if the latter process is a set-valued martingale. Then, we narrow down to the case of a random polytope and show that these equivalent conditions are also equivalent to the random polytope having a deterministic \emph{normal fan}, i.e., while the vertices of the polytope are random, they have the same sets of supporting directions (a.k.a. \emph{normal cones}) with probability one. This observation seems to be new and it provides a bridge between set-valued stochastic analysis and convex geometry, which could be discovered further in future studies.

The rest of the paper is organized as follows. In Section~\ref{sec:prelim}, we provide some preliminaries on set-valued analysis with a special focus on set-valued martingales in discrete time. Section~\ref{sec:reverse} is on set-valued reverse submartingales in discrete time and it serves as a preparation for the right-continuity result in continuous time. In Section~\ref{sec:ctmtg}, we study the existence of the path limits and the regularity properties of set-valued submartingales in continuous-time. We focus on the set-valued stochastic integral in Section~\ref{sec:integral} and characterize the cases where it gives a set-valued martingale. In Section~\ref{sec:compare}, we compare a set-valued martingale and the closed convex hull of a sequence of vector-valued martingales, which yields a set-valued submartingale. As a special case, we consider the conditional expectation of a convex random polytope in Section~\ref{sec:polytope} and show some links between set-valued stochastic integrals and normal fans.

\section{Preliminaries on Set-Valued Analysis}\label{sec:prelim}

In this section, we recall some important notions about the convergence of sets, set-valued measurable functions, and set-valued martingales in discrete time. Throughout this paper, we consider functions whose values are subsets of the finite-dimensional Euclidean space $\hX=\hR^d$,  $d\in\hN$. However, many of the results are valid when $\hX$ is any separable Banach space.

To begin with, given nonempty sets $C, D$, a function $f\colon C\to D$, and a subset $A\subseteq C$, we write $f[A]:=\{f(x)\colon x\in A\}$ for the image of $A$ under $f$; we define its indicator function $\1_A\colon C\to \{0,1\}$ by $\1_A(x)=1$ if $x\in A$ and by $\1_A(x)=0$ if $x\in A^c:=C\setminus A$. If $C$ is a topological space, then $\cl_{C}(A)$ (resp. $\interior_{C}(A)$, $\bd_C(A)$) denotes the closure (resp. interior, boundary) of $A$ in $C$. For each $n\in\hN$, the unit simplex in $\hR^n$ is denoted by $\Delta^{n-1}$, i.e., $\Delta^{n-1}:=\cb{r\in\hR^n\colon r_1,\ldots,r_n\geq 0;\ \sum_{i=1}^n r_i=1}$. If $C$ is a vector space, then $\co(A)$ denotes the convex hull of $A$ and we have
\[
\co(A)=\cb{\sum_{i=1}^n r_ix^i\colon x^1,\ldots,x^n\in A;\ r\in\Delta^{n-1};\ n\in\hN}.
\]

\subsection{Convergence of Sets}

Let us denote $\lan\cd, \cd\ran$ and $\abs{\cdot}$ to be the inner product and Euclidean norm on $\hR^d$, respectively;  and $\hB_{\hR^d}(r):=\{x\in\hR^d\colon \abs{x}\leq r\}$, $r>0$. In particular, we write $\hB_{\hR^d}:=\hB_{\hR^d}(1)$ for simplicity. For $C\subseteq\hR^d$, we denote $\sB(C)$ to be the Borel $\sigma$-algebra on $C$; and we denote $\sC(\hR^d)$ (resp. $\sG(\hR^d)$, $\sK(\hR^d)$) to be the set of all nonempty and closed (resp. closed convex, compact convex) subsets of $\hR^d$. We shall also fix a countable dense subset $\sW$ of $\hB_{\hR^d}$ such that $0\notin\sW$.

For $C\in\sC(\hR^d)$, we define the   \emph{distance function} $d(\cdot,C)\colon\hR^d\to\ol{\hR}_+$   by
$d(x,C):=\inf_{y\in C}\abs{x-y}$, $x\in\hR^d$; and 
for $C, D\in\sC(\hR^d)$, let $\bar{h}(C,D):=\sup_{x\in C}d(x,D)$. We define the \emph{Hausdorff distance} between $C, D$ by $h(C,D):=\max\{\bar{h}(C,D),\bar{h}(D,C)\}$; we also define $\norm{C}:=h(C,\{0\})=\sup\{\abs{x}\colon x\in C\}$. It is well-known that (see \cite[\S3.2]{beer-book})
\bea\label{H-W}
h(C,D)=\sup_{x\in \hR^d}\abs{d(x,C)-d(x,D)}.
\eea 
Furthermore, we define the \emph{support function} $s(\cdot,C) \colon \hR^d \to(-\infty,+\infty]$ associated to $C$ by 
 \[
s(x^\ast,C):=\sup_{x\in C}\ip{x^\ast,x},\quad x^\ast\in\hR^d.
\]
We should note that a closed convex set is uniquely identified by its support function. Indeed, as a consequence of separating hyperplane theorem, the support function is {\it monotone} in the sense that,  for $C,D\in\sG(\hR^d)$, 
$C\subseteq D$ if and only if
% \quad \Longleftrightarrow \quad \
$ s(x^\ast,C)\leq s(x^\ast,D)$, $x^\ast\in \hR^d$. Furthermore, if $C,D$ are bounded, then we also have (see \cite[Theorem~3.2.7]{beer-book})
\bea\label{H-S}
h(C,D) = \sup_{x^\ast\in\hB_{\hR^d}}\abs{s(x^\ast,C)-s(x^\ast,D)}.
\eea 
Finally, for a given sequence of sets $(C_n)_{n\in\hN}$ in $\sC(\hR^d)$, we define the sets (with $n_0:=0$)
\beaa
&&\underset{n\rightarrow\infty}{\liminf{}} C_n:=\cb{x\in\hR^d\colon x=\underset{n\rightarrow\infty}{\lim}x_{n};\ \forall n\in\hN\colon x_n\in C_n},\\
&&\underset{n\rightarrow\infty}{\limsup{}} C_n:=\cb{x\in\hR^d\colon x=\underset{k\rightarrow\infty}{\lim}x_{n_k};\ \forall k\in\hN\colon n_k\in\hN, n_{k}>n_{k-1},\ x_{n_k}\in C_{n_k}}.
\eeaa

Using these notions, we may introduce several modes of convergence for sets.

\begin{defn}
Let $C, C_1, C_2,\ldots \in \sC(\hR^d)$. We say that

(a) $(C_n)_{n\in\hN}$ Wijsman converges to $C$ if $\lim_{n\rightarrow\infty}d(x,C_n)=d(x,C)$ for every $x\in\hR^d$,

(b) $(C_n)_{n\in\hN}$ Hausdorff converges to $C$ if $\lim_{n\rightarrow\infty}h(C_n,C)=0$,

(c) $(C_n)_{n\in\hN}$ scalarly converges to $C$ if $\lim_{n\rightarrow\infty}s(x^\ast,C_n)=s(x^\ast,C)$ for every $x^\ast\in\hB_{\hR^d}$,

(d) $(C_n)_{n\in\hN}$ Painlev\'{e}-Kuratowski converges to $C$ if $\limsup_{n\rightarrow\infty} C_n\subseteq C\subseteq  \liminf_{n\rightarrow\infty}C_n$.
\end{defn}

We now summarize some well-known connections between these modes of convergence which will be useful in our discussions below.

\begin{prop}\label{conv-finite}
	 \emph{\cite[Theorems~1,~2,~6]{wets79}, \cite[Proposition~2.3]{hess1999}} Let $C, C_1, C_2,\ldots \in \sC(\hR^d)$ and consider the following properties:
	 
	(a) $(C_n)_{n\in\hN}$ Hausdorff converges to $C$.
	
	(b) $\lim_{n\rightarrow\infty} d(x,C_n)=d(x,C)$ uniformly in $x\in\hR^d$.
	
	(c) $\lim_{n\rightarrow\infty} s(x^\ast,C_n)=s(x^\ast,C)$ uniformly in $x^\ast\in\hB_{\hR^d}$.
	
	(d) $(C_n)_{n\in\hN}$ Painlev\'{e}-Kuratowski converges to $C$.
	
	(e) $(C_n)_{n\in\hN}$ Wijsman converges to $C$.
	
	(f) $(C_n)_{n\in\hN}$ scalarly converges to $C$.
	
	(g) $\lim_{n\rightarrow\infty} s(x^\ast,C_n)=s(x^\ast,C)$ for every $x^\ast\in\sW$.
	
	Then, the following results hold.
	
	(i) In general, (a)$\Leftrightarrow$(b)$\Rightarrow$(e), (a)$\Rightarrow$(d), (c)$\Rightarrow$(f)$\Rightarrow$(g).
	
	(ii) Suppose that $C, C_1, C_2,\ldots\in \sG(\hR^d)$. Then, (a)$\Leftrightarrow$(c), (d)$\Leftrightarrow$(e) as well.
	
	(iii) Suppose that $C, C_1, C_2, \ldots \in \sK(\hR^d)$. Then, (a)$\Leftrightarrow$(f)$\Leftrightarrow$(g) as well.
	\end{prop}
	
We also have the following Bolzano-Weierstrass property for compact subsets of $\hR^d$:

\begin{prop}\label{BW}
	Let $(C_n)_{n\in\hN}$ be a sequence of nonempty compact subsets of $\hR^d$. If $(C_n)_{n\in\hN}$ is bounded, i.e., $\sup_{n\in\hN}\norm{C_n}<+\infty$, then there exists a subsequence $(C_{n_k})_{k\in\hN}$ and a nonempty compact set $C\subseteq\hR^d$ such that $(C_{n_k})_{k\in\hN}$ Hausdorff converges to $C$. If we further have that $(C_n)_{n\in\hN}$ is in $\sK(\hR^d)$, then $C\in\sK(\hR^d)$.
\end{prop}

{\it Proof}: The first claim is by \cite[Proposition~2.2]{hess1999}. The second claim follows from the fact that $\sK(\hR^d)$ can be embedded as a closed convex cone in a Banach space thanks to H\"{o}rmander's theorem; see \cite[Theorem~3.2.9]{beer-book}.
\qed

\subsection{Set-Valued Measurable Functions}\label{sec:svfunctions}

Let us fix a probability space $(\O,\cF,\hP)$. When $\cA$ is a collection of subsets of $\O$, we denote by $\sigma(\cA)$ the $\sigma$-algebra on $\O$ generated by $\cA$. Let $\cG$ be a sub-$\sigma$-algebra of $\cF$. We denote by $\Pi_\cG(\O)$ the set of all $\cG$-measurable finite partitions of $\O$. We denote by $\hL^0_\cG(\O,\hR^d)$ the set of all $\cG$-measurable random vectors $\xi\colon\O\to \hR^d$ that are distinguished up to $\hP$-almost sure (a.s.) equality. For $p\geq 1$, we denote by $\hL^p_\cG(\O,\hR^d)$ the set of all $\xi\in \hL^0_\cG(\O,\hR^d)$ that are $p$-integrable, i.e., $\hE[\abs{\xi}^p]<+\infty$; we also write $\hL^p:=\hL^p_{\cF}(\O,\hR^d)$. A set $\cK\subseteq \hL^0_\cF(\O,\hR^d)$ is called \emph{$\cG$-decomposable} if $\xi\1_A+\zeta\1_{A^c}\in \cK$ for every $\xi,\zeta\in \cK$ and $A\in\cG$. In general, the smallest $\cG$-decomposable superset of $\cK$ is given by
\bea\label{eq:dec}
\dec_{\cG}(\cK)=\cb{\sum_{i=1}^n \xi^i \1_{A_i}\colon \xi^1,\ldots,\xi^n \in \cK;\ (A_i)_{i=1}^n\in\Pi_{\cG}(\O);\ n\in\hN}
\eea
and we call it the \emph{$\cG$-decomposable hull} of $\cK$. When $\cK\subseteq \hL^p$ for some $p\geq 1$, we say that $\cK$ is \emph{bounded} in $\hL^p$ if $\sup_{\xi\in\cK}\hE[|\xi|^p]<+\infty$ and \emph{dominated} in $\hL^p$ if there exists $\xi\in \hL^p$ such that $\abs{\zeta}\leq \xi$ $\hP$-a.s. for every $\zeta\in \cK$.

The next result is concerned with the relationship between being dominated and bounded for sets of random vectors.

\begin{lem}
\label{kis-decomp-bdd}
	 Let $\cG$ be a sub-$\sigma$-algebra of $\cF$, $p\geq 1$, and $\cK\subseteq \hL^p_\cG(\O,\hR^d)$ be a nonempty set. Then, the following are equivalent:
	 
	(i) $\cK$ is dominated in $\hL^p$.
	
	(ii) $\dec_{\cG}(\cK)$ is dominated in $\hL^p$.
	
	(iii) $\cl_{\hL^p}\dec_{\cG}(\cK)$ is dominated in $\hL^p$.
	
	(iv) $\dec_{\cG}(\cK)$ is bounded in $\hL^p$.
	
	(v) $\cl_{\hL^p}\dec_{\cG}(\cK)$ is bounded in $\hL^p$.
	\end{lem}

{\it Proof}: The equivalences (i)$\Leftrightarrow$(iii)$\Leftrightarrow$(v) are given by \cite[Theorem~3.3.5]{Kis2020}. The implications (iii)$\Rightarrow$(ii), (ii)$\Rightarrow$(i), (v)$\Rightarrow$(iv) are trivial as $\cK\subseteq\dec_{\cG}(\cK)\subseteq\cl_{\hL^p}\dec_{\cG}(\cK)$. To complete the proof, we check (iv)$\Rightarrow$(v). Suppose that (iv) holds and let
$r:=\sup_{\xi\in \dec_{\cG}(\cK)}\hE[|\xi|^p]<+\infty$.
Let $\xi\in \cl_{\hL^p}\dec_{\cG}(\cK)$. Then, there exists a sequence $(\xi^n)_{n\in\hN}$ in $\dec_{\cG}(\cK)$ that converges to $\xi$ in $\hL^p$. Hence, $\hE[|\xi|^p]=\lim_{n\rightarrow\infty}\hE[|\xi^n|^p]\leq r$ and (v) follows.
\qed 

We say that a set-valued function $\Xi\colon \O\to \sC(\hR^d)$ is \emph{$\cG$-measurable} if 
$\{\o\in\O\colon \Xi(\o)\cap C\neq \emptyset \}\in\cG$, 
for every $C\in\sC(\hR^d)$; the set of all such functions is denoted by $\scL^0_\cG(\O,\sC(\hR^d))$. It can be checked that the scalar function $\norm{\Xi}$ is $\cG$-measurable whenever $\Xi$ is $\cG$-measurable. A random vector $\xi\in\hL^0_{\cG}(\O,\hR^d)$ is called a \emph{measurable selection} of $\Xi$ if $\hP\{\xi\in \Xi\}=1$; the set of 
all measurable selections of $\Xi$ is denoted by $S^0_\cG(\Xi)$. For $p\geq 1$, let $S^p_{\cG}(\Xi):= S^0_\cG(\Xi)\cap \hL^p$; we say that $\Xi$ is \emph{$p$-integrable} if $S^p_\cG(\Xi)\neq \emptyset$ and it is \emph{$p$-integrably bounded} if $\hE[\norm{\Xi}^p]<+\infty$.
We also define
\beaa 
&&\sA^p_\cG(\O,\sC(\hR^d)):=\{\Xi\in \scL^0_\cG(\O,\sC(\hR^d))\colon S^p_\cG(\Xi)\neq\emptyset\},\\
&&\scL^p_\cG(\O,\sC(\hR^d)):=\{\Xi\in \scL^0_\cG(\O,\sC(\hR^d))\colon \hE[\norm{\Xi}^p]<+\infty\}.
\eeaa 

We list some main results on set-valued measurable functions in the next theorem.

\begin{thm}\label{decomp-meas}
	Let $\cG$ be a sub-$\sigma$-algebra of $\cF$ and $p\geq 1$.
	
	(i) \emph{\cite[Theorem~1.4.1, Proposition~2.1.4]{Molchanov}} Let $\Xi\in \scL^0_\cG(\O,\sC(\hR^d))$. Then, $S^0_\cG(\Xi)\neq \emptyset$. If $\Xi$ is $p$-integrably bounded, then $S^0_\cG(\Xi)=S^p_\cG(\Xi)$ and $\Xi$ has bounded values $\hP$-a.s. In particular, $\scL^p_\cG(\O,\sC(\hR^d))\subseteq \sA^p_\cG(\O,\sC(\hR^d))$.
	
	(ii) \emph{\cite[Proposition~2.1.7]{Molchanov}} Let $\Xi\in \sA^p_\cG(\O,\sC(\hR^d))$. Then, $S^p_\cG(\Xi)$ is a nonempty (strongly) closed subset of $\hL^p$. Moreover, $\Xi$ has convex values $\hP$-a.s. if and only if $S^p_\cG(\Xi)$ is a convex set; $\Xi$ is $p$-integrably bounded if and only if $S^p_\cG(\Xi)$ is bounded in $\hL^p$.
	
	(iii) \emph{\cite[Proposition~2.1.4]{Molchanov}} Let $\Xi,\Upsilon\in \scL^0_\cG(\O,\sC(\hR^d))$ (resp. $\sA^p_\cG(\O,\sC(\hR^d))$). Then, $\Xi=\Upsilon$ $\hP$-a.s. if and only if $S^0_\cG(\Xi)=S^0_\cG(\Upsilon)$ (resp. $S^p_\cG(\Xi)=S^p_\cG(\Upsilon)$).
	
	(iv) \emph{\cite[Theorem~2.1.18]{Molchanov}} Let $\Xi\in\scL^p_\cG(\O,\sC(\hR^d))$ and suppose that $\Xi$ has convex values $\hP$-a.s. Then, $\Xi$ has compact values $\hP$-a.s. if and only if $S^p_\cG(\Xi)$ is a weakly compact set. %In particular, when $\hX=\hR^d$ for some $d\in\hN$, $S^p_\cG(\Xi)$ is a weakly compact set.
	
	(v) \emph{\cite[Theorem~2.1.10]{Molchanov}} Let $\cK\subseteq \hL^p_\cG(\O,\hR^d)$ be a nonempty closed set. Then, $\cK$ is $\cG$-decomposable if and only if there exists $\Xi\in \sA^p_\cG(\O,\sC(\hR^d))$ such that $\cK=S^p_\cG(\Xi)$.
	\end{thm}

In view of Theorem \ref{decomp-meas}(iii), we may assume that the members of $\scL^0_\cG(\O,\sC(\hR^d))$ are identified up to $\hP$-a.s. equality similar to the case of random vectors. For $p\geq 1$, we denote by $\scL^p_\cG(\O,\sK(\hR^d))$ the set of all $\Xi\in \scL^p_\cG(\O,\sC(\hR^d))$ with compact convex values $\hP$-a.s.

Finally, we review the notion of conditional expectation for set-valued random variables. Given $\Xi\in\sA^1_\cF(\O,\sC(\hR^d))$, the set $\hE[\Xi]:=\cl_{\hR^d}(\{\hE[\xi]\colon \xi\in S^1_\cF(\Xi)\})$ is called the \emph{(closed) Aumann integral} of $\Xi$. It is easy to check that the set $\{\hE[\xi|\cG]\colon \xi\in S^1_\cG(\Xi)\}$ is $\cG$-decomposable. Hence, by Theorem~\ref{decomp-meas}(v) (see also \cite[Theorem~2.1.71]{Molchanov}), there exists a   unique set-valued random variable $\hE[\Xi|\cG]\in\sA^1_\cG(\O,\sC(\hR^d))$ such that $S^1_\cG(\hE[\Xi|\cG])=\cl_{\hL^1}(\{\hE[\xi|\cG]\colon \xi\in S^1_\cG(\Xi)\})$. Moreover, if $\Xi$ is \emph{$p$-integrable} (resp. \emph{$p$-integrably bounded}), then so is $\hE[\Xi|\cG]$ and we have $S^p_\cG(\hE[\Xi|\cG])=\cl_{\hL^p}(\{\hE[\xi|\cG]\colon \xi\in S^p_\cG(\Xi)\})$ (resp. $S^p_\cG(\hE[\Xi|\cG])=\{\hE[\xi|\cG]\colon \xi\in S^p_\cG(\Xi)\}$); if $\Xi$ has convex (resp. compact convex) values $\hP$-a.s., then so does $\hE[\Xi|\cG]$.

\subsection{Set-Valued Submartingales in Discrete Time}

Let $\hG=(\cG_n)_{n\in\hN}$ be a discrete-time filtration on $(\O,\cF,\hP)$, and let $\cG_\infty:=\bigvee_{n\in\hN}\cG_n:=\sigma(\bigcup_{n\in\hN}\cG_n)$. A discrete time \emph{set-valued (stochastic) process} is a collection $(\Xi_n)_{n\in\hN}$ with $\Xi_n \in\scL^0_\cF(\O,\sC(\hR^d))$, $n\in\hN$. It is called \emph{$\hG$-adapted} if $\Xi_n\in\scL^0_{\cG_n}(\O,\sC(\hR^d))$, \emph{convex} if $\Xi_n\in\scL^0_\cF(\O,\sG(\hR^d))$, \emph{integrable} if $\Xi_n\in\sA^1_{\cF}(\O,\sC(\hR^d))$, and \emph{integrably bounded} if $\Xi_n\in\scL^1_{\cF}(\O,\sC(\hR^d))$, for each $n\in\hN$.

With these notions, we are ready to define set-valued submartingales.

\begin{defn}\label{defn:dtmtg}
A $\hG$-adapted integrable set-valued process $(\Xi_n)_{n\in\hN}$ is called a \emph{set-valued $\hG$-submartingale} (resp. \emph{set-valued $\hG$-supermartingale}) if $\hE[\Xi_{n+1}|\cG_n]\supseteq\Xi_n$ (resp. $\hE[\Xi_{n+1}|\cG_n]\subseteq\Xi_n$) $\hP$-a.s. for each $n\in\hN$, it is called a \emph{set-valued $\hG$-martingale} if it is both a set-valued $\hG$-submartingale and a set-valued $\hG$-supermartingale.
\end{defn}

\begin{rem}\label{rem:vector-vs-set}
By definition, it is clear that a vector-valued $\hG$-martingale $(\xi_n)_{n\in\hN}$ can be  viewed as a singleton  set-valued $\hG$-martingale. But there is a dramatic distinction between real-valued and set-valued submartingales. Indeed, if $(\xi_n)_{n\in\hN}$ is a real-valued $\hG$-submartingale, then it holds that $\hE[\xi_{n+1}|\cG_n]\geq \xi_n$, $n\in\hN$. But this by no means implies that $\hE[\{\xi_{n+1}\}|\cG_n]\supseteq \{\xi_n\}$ as sets, and hence $(\{\xi_n\})_{n\in\hN}$ is not a set-valued $\hG$-submartingale. It is interesting to note that as singleton set-valued process the desired ``inclusion" only holds when $(\xi_n)_{n\in\hN}$ is actually a $\hG$-martingale.
 \qed
 \end{rem}

We collect some basic observations about set-valued martingales in the next proposition.

\begin{prop}\label{mtg-ex}
	(i) \emph{\cite[Example~5.1.2(i)]{Molchanov}} Let $\Upsilon\in\scL^1_{\cF}(\O,\sG(\hR^d))$ and set $\Xi_n:=\hE[\Upsilon | \cG_n]$ for each $n\in\hN$. Then, $(\Xi_n)_{n\in\hN}$ is an integrably bounded convex set-valued $\hG$-martingale.\\
	(ii) \emph{\cite[Theorem~5.1.5(i)]{Molchanov}} Let $(\Xi_n)_{n\in\hN}$ be an integrably bounded convex $\hG$-submartingale. Then, $(\norm{\Xi_n})_{n\in\hN}$ is a real-valued $\hG$-submartingale.\\
	(iii) \emph{\cite[Theorem~5.1.5(iii)]{Molchanov}} Let $(\Xi_n)_{n\in\hN}$ be an integrably bounded convex $\hG$-martingale (resp. submartingale, supermartingale). Then, for each $x^\ast\in\hR^d$, the sequence $(s(x^\ast,\Xi_n))_{n\in\hN}$ is a real-valued $\hG$-martingale (resp. submartingale, supermartingale).
	\end{prop}
		
We would like to remark here that, unlike the real-valued setting where submartingales and supermartingales are negatives of each other, whence in a sense ``symmetric", in the set-valued case submartingales and supermartingales have to be studied separately as is clear from Definition~\ref{defn:dtmtg}. However, in this paper, we shall focus more on the set-valued submartingale for two reasons. First, the important  and frequently used ``scalarization"  procedure in this paper: $(\Xi_n)_{n\in\hN}\mapsto (\norm{\Xi_n})_{n\in\hN}$ only preserves the submartingale property, thanks to Proposition~\ref{mtg-ex}(ii). Second, the main objective of this paper is the path-regularity of the set-valued stochastic integrals, which are set-valued submartingales in general (see Lemma~\ref{lem:submtg}). Bearing these in mind, in the rest of the paper, we shall consider only set-valued submartingales and martingales, and leave some of the non-trivial supermartingale cases to future studies.

We end this section by recalling a set-valued version of \emph{Martingale Convergence Theorem}.

\begin{thm}\label{dtmtg:findim}
\emph{\cite[Theorem~4.5]{hess1999}} Let $(\Xi_n)_{n\in\hN}$ be a convex set-valued $\hG$-submartingale, such that the family $\{\norm{\Xi_n}\colon n\in\hN\}$ is bounded in $\hL^1_{\cF}(\O,\hR)$. Then, there exists a set-valued random variable $\Xi_\infty \in \scL^1_{\cG_\infty}(\O,\sK(\hR^d))$ such that $\lim_{n\rightarrow\infty}\Xi_n=\Xi_\infty$ in the Hausdorff sense, $\hP$-a.s. Furthermore, if $\{\norm{\Xi_n}\colon n\in\hN\}$ is uniformly integrable, then the Hausdorff convergence also holds in $\hL^1$.
\end{thm}
	
\section{Set-Valued Reverse Martingales in Discrete Time}\label{sec:reverse}

In this section, we establish a key tool for the study of path-regularity of continuous time set-valued martingales: the convergence theorem for \emph{reverse set-valued martingales}.

Let $\hG_{-}=(\cG_{-n})_{n\in\hN}$ be a \emph{reverse filtration}, that is, each $\cG_{-n}$ is a sub-$\sigma$-algebra of $\cF$ such that $\cG_{-1}\supseteq\cG_{-2}\supseteq\ldots$ We define $\cG_{-\infty}:=\bigcap_{n\in\hN}\cG_{-n}$. A set-valued process $(\Xi_n)_{n\in\hN}$ is called \emph{$\hG_{-}$-adapted} if $\Xi_n$ is $\cG_{-n}$-measurable for each $n\in\hN$. A $\hG_{-}$-adapted integrable set-valued process $(\Xi_n)_{n\in\hN}$ is called a \emph{set-valued $\hG_{-}$-submartingale} (resp. \emph{set-valued $\hG_{-}$-supermartingale}) if $\hE[\Xi_{n+1}|\cG_{-n}]\subseteq \Xi_n$ (resp. $\hE[\Xi_{n+1}|\cG_{-n}]\supseteq \Xi_n$) $\hP$-a.s. for each $n\in\hN$, it is called a \emph{set-valued $\hG_-$-martingale} if it is both a set-valued $\hG_-$-submartingale and a set-valued $\hG_-$-supermartingale.

It is easy to check that Remark~\ref{mtg-ex} extends naturally for reverse filtrations. The following is the reverse version of Theorem~\ref{dtmtg:findim}, as well as a set-valued generalization of the converge theorem for real-valued reverse submartingales (cf. \cite[Theorem~V.4.19]{cinlar}). To the best of our knowledge, it is new. Also, as we pointed out before, the set-valued submartingales and real-valued submartingales are fundamentally different (see Remark~\ref{rem:vector-vs-set}), so the condition in the set-valued version is actually weaker than the real-valued case. We provide a detailed proof for the interested reader.

\begin{thm}\label{thm:reverse}
	Let $(\Xi_n)_{n\in\hN}$ be a convex integrably bounded set-valued $\hG_-$-submartingale. Then, $\{\norm{\Xi_n}\colon n\in\hN\}$ is uniformly integrable and there exists a set-valued random variable $\Xi_\infty\in\scL^1_{\cG_{-\infty}}(\O,\sK(\hR^d))$ such that $\lim_{n\rightarrow\infty}\Xi_n=\Xi_\infty$, 
	in the Hausdorff sense $\hP$-a.s. and in $\hL^1$.
	\end{thm}
	
{\it Proof}: The proof follows the idea of that of  \cite[Theorem~V.4.19]{cinlar}, with some extra steps to deal with the special set-valued nature. 
To begin with, we first note that, by a reverse version of Proposition \ref{mtg-ex}(ii), $(\norm{\Xi_n})_{n\in\hN}$ is a positive reversed $\hG_-$-submartingale, thus by \cite[Theorem~V.4.19]{cinlar} it is uniformly integrable and hence converges to a random variable $m_\infty\in\hL^1_{\cG_{-\infty}}$, both $\hP$-a.s. and in $\hL^1$. We denote $\cN_0$ to be the exceptional $\hP$-null set here.

Next, for each $n\in\hN$, we have $\hE[\Xi_n]\in \sK(\hR^d)$ since $\Xi_n$ is integrably bounded. Moreover, the submartingale property implies that $\hE[\Xi_1]\supseteq\hE[\Xi_2]\supseteq\ldots$ Hence, $D:=\bigcap_{n\in\hN}\hE[\Xi_n]\neq \emptyset$, thanks to Cantor's intersection theorem.
Furthermore, for each $x^\ast\in\sW$, the reverse version of Proposition~\ref{mtg-ex}(iii) implies that $(s(x^\ast,\Xi_n))_{n\in\hN}$ is a real-valued $\hG_-$-submartingale, and by the definitions of set $D$ and the support function $s(x^*, \cd)$, it is easy to check (cf, e.g., \cite[Theorem~2.1.35]{Molchanov}) that  
\[
\inf_{n\in\hN}\hE[s(x^\ast,\Xi_n)]=\inf_{n\in\hN}s(x^\ast,\hE[\Xi_n])\geq s(x^\ast,D)>-\infty.
\]
Therefore, the real-valued reverse martingale convergence theorem implies that $(s(x^\ast,\Xi_n))_{n\in\hN}$ converges to a random variable $s^{x^\ast}_\infty\in \hL^1_{\cG_{-\infty}}$, both $\hP$-a.s. and in $\hL^1$. We denote $\cN_{x^\ast}$ to be the exceptional $\hP$-null set, for each $x^*\in\sW$.

Now let us define $\cN=\cN_0\cup\big(\bigcup_{x^*\in\sW}\cN_{x^*}\big)$, and $\ol{\O}:=\O\setminus \cN$. 
Then, $\hP(\ol{\O})=1$ and, as a convergent sequence, $(\norm{\Xi_n(\o)})_{n\in\hN}$ is bounded, for  $\o\in\ol{\O}$. Then, by Bolzano-Weierstrass Theorem (Proposition~\ref{BW}), there exists a subsequence $(\Xi_{n_k}(\o))_{k\in\hN}$ that converges to some $C(\o)\in\sK(\hR^d)$ in the Hausdorff sense. Since Hausdorff convergence implies scalar convergence, we have
\[
s^{x^\ast}_\infty(\o)=\lim_{n\rightarrow\infty}s(x^\ast,\Xi_n(\o))=\lim_{k\rightarrow\infty}s(x^\ast,\Xi_{n_k}(\o))=s(x^\ast,C(\o)), \q  x^\ast\in\sW.
\]
This in turn implies that $(\Xi_n(\o))_{n\in\hN}$ itself converges to $C(\o)$ in the Hausdorff sense, thanks to Proposition~\ref{conv-finite}. We can thus denote $\Xi_\infty(\o):=C(\o)$, $\o\in \bar{\O}$. Setting  $\Xi_\infty(\o):=\{0\}$ for $\o\in\cN$ we can extend $\Xi_\infty$ to all $\O$, so that $(\Xi_n)_{n\in\hN}$ Hausdorff converges to $\Xi_\infty$ $\hP$-a.s. Furthermore, note that $s(x^\ast,\Xi_\infty(\o))=s^{x^\ast}_{\infty}(\o)$, $\o\in\bar{\O}$, for all $x^\ast\in\sW$ and $s^{x^\ast}_{\infty}$ is $\cG_{-\infty}$-measurable. We conclude that $s(x^\ast,\Xi_\infty)$ is $\cG_{-\infty}$-measurable for all $x^\ast\in\sW$. That is, $\Xi_\infty$ is $\cG_{-\infty}$-measurable.

Finally, note that $\abs{\norm{\Xi_n}-\norm{\Xi_\infty}}\leq h(\Xi_n,\Xi_\infty)$ for each $n\in\hN$, $(\norm{\Xi_n})_{n\in\hN}$ converges to $\norm{\Xi_\infty}$ $\hP$-a.s. Then, by Fatou's lemma and the  submartingale property, we have
\[
\hE[\norm{\Xi}_\infty]=\hE\Big[\lim_{n\rightarrow\infty}\norm{\Xi_n}\Big]\leq\liminf_{n\rightarrow\infty}\hE[\norm{\Xi_n}]\leq \sup_{n\in\hN}\hE[\norm{\Xi_n}]=\hE[\norm{\Xi_1}]<+\infty,
\]
that is,  $\Xi_\infty$ is integrably bounded. Moreover,  note that $h(\Xi_n,\Xi_\infty)\leq \norm{\Xi_n}+\norm{\Xi_\infty}$ for each $n\in\hN$, and $(\norm{\Xi_n})_{n\in\hN}$ is uniformly integrable, we see that $(h(\Xi_n,\Xi_\infty))_{n\in\hN}$ is uniformly integrable as well. This, together with the fact that $(h(\Xi_n,\Xi_\infty))_{n\in\hN}$ converges to $0$,  $\hP$-a.s.,  implies that the convergence is in $\hL^1$ as well, i.e., $(\Xi_n)_{n\in\hN}$ Hausdorff converges to $\Xi_\infty$ in $\hL^1$.
\qed 

\section{Set-Valued Submartingales in Continuous Time}\label{sec:ctmtg}

In this section, we study the path-regularity for set-valued submartingales in continuous time. We begin by introducing some basic definitions. Let $D\in\sB([0,+\infty))$ be a nonempty index set. We will later take either $D=[0,+\infty)$ or $D=[0,T]$ for some $T>0$. We denote by $\Leb$ the Lebesgue measure on $(D,\sB(D))$. Let $\hF=(\cF_t)_{t\in D}$ be a filtration on $(\O,\cF,\hP)$. Unless otherwise stated, we do \emph{not} assume that $\hF$ satisfies the {\it usual conditions} (i.e., right-continuous and aumented, cf. e.g., \cite{cinlar}) in general. We denote by $\hL^0_{\hF}(D\times\O,\hR^d)$ the set of all $\hF$-progressively measurable processes $z=(z_t)_{t\in D}$ taking values in $\hR^d$ that are distinguished up to $\Leb\otimes\hP$-almost everywhere (a.e.) equality and by $\hL^2_{\hF}(D\times\O,\hR^d)$ the set of all $z\in \hL^0_{\hF}(D\times\O,\hR^d)$ such that $\hE[\int_D \abs{z_t}^2 dt]<+\infty$.

A set-valued function $\Phi\colon D\times\O\to\sC(\hR^d)$ can be identified as a family $\Phi=(\Phi_t)_{t\in D}$, where $\Phi_t(\o):=\Phi(t,\o)$ for each $(t,\o)\in D\times\O$; in this case, $\Phi$ is called a \emph{set-valued (stochastic) process} (with index set $D$) if $\Phi_t\in\scL^0_\cF(\O,\sC(\hR^d))$ for each $t\in D$. We say that a set-valued process $\Phi$ is \emph{$\hF$-progressively measurable} if it is measurable with respect to the $\hF$-progressive $\sigma$-algebra on $D\times\O$, the collection of all such processes is denoted by $\scL^0_\hF(D\times\O,\sC(\hR^d))$. The definitions of $\hF$-adaptedness, convexity, $p$-integrability and $p$-integrable boundedness ($p\geq 1$) are analogous to the discrete-time case. \emph{Set-valued $\hF$-submartingales}, \emph{supermartingales}, and \emph{martingales} with time-index $D$ can be defined naturally following Definition~\ref{defn:dtmtg}. Furthermore, it is easy to see that Proposition~\ref{mtg-ex} extends naturally to the continuous-time setting.

\subsection{Left- and Right-Limits of Set-Valued Submartingales}

In the rest of this section, we consider the case $D=[0,+\infty)$. Our first result provides the existence of all left- and right-limits of a set-valued submartingale with probability one.

\begin{lem}\label{leftlimited}
	Let $(M_t)_{t\geq 0}$ be a set-valued convex and integrably bounded set-valued $\hF$-submartingale. Then, with probability one, the following limits exist in the Hausdorff sense:
	\[
	M_{t-}:=\lim_{r\in\hQ,r\uparrow t}M_r,\quad t>0;\quad M_{t+}:=\lim_{r\in\hQ,r\downarrow t}M_r,\quad t\geq 0.
	\]
	\end{lem}
	
{\it Proof}: The scheme of the proof is similar to that of Theorem \ref{thm:reverse}. We first consider the ``scalarization" $m_t:=\norm{M_t}$, $t\geq 0$. Since $(m_t)_{t\geq 0}$ is a real-valued $\hF$-submartingale, by \cite[Proposition~V.7.14]{cinlar}, there exists a set $\O_0\in \cF$ such that $\hP(\O_0)=1$ and, for $\o\in\O_0$, the limits
$m_{t-}(\o):=\lim_{r\in\hQ,r\uparrow t}m_r(\o)$, $t>0$, and $m_{t+}(\o):=\lim_{r\in\hQ,r\downarrow t}m_r(\o)$, $t\geq 0$,  exist and are finite.

Next, we consider another ``scalarization": for $x^\ast\in \sW$, define $s^{x^\ast}_t:=s(x^\ast,M_t)$, $t\geq 0$. Again, since each $(s^{x^\ast}_t)_{t\geq 0}$ is a real-valued $\hF$-submartingale,  there exists a set $\O_{x^\ast}\in \cF$ such that $\hP(\O_{x^\ast})=1$ and the limits
$s^{x^\ast}_{t-}(\o):=\lim_{r\in\hQ,r\uparrow t}s^{x^\ast}_r(\o)$, $t>0$, and $s^{x^\ast}_{t+}(\o):=\lim_{r\in\hQ,r\downarrow t}s^{x^\ast}_r(\o)$,$t\geq 0$, exist and are finite, whenever $\o\in\O_{x^\ast}$.

Let us now focus on the left limit case, as the right limit case can be argued similarly. To this end, denote $\ol{\O}=\O_0\cap\big(\bigcap_{x^*\in\sW}\O_{x^*}\big)$, so that $\hP(\ol{\O})=1$. For $t>0$, let $(r_n(t))_{n\in\hN}\subset \hQ$ be such that $r_n(t)\uparrow t$ as $n\rightarrow\infty$. Since for $\o\in\ol{\O}$, it holds that $\lim_{n\to\infty}\norm{M_{r_n(t)}(\o)}=m_{t-}(\o)$, the sequence $(M_{r_n(t)}(\o))_{n\in\hN}$ is bounded. Applying Bolzano-Weierstrass Theorem (Proposition~\ref{BW}), we can find a Hausdorff convergent subsequence $(M_{r_{n_k}(t)}(\o))_{k\in\hN}$, and denote its limit by $A_t(\o)\in\sK(\hR^d)$. Now,  by using the relationship between Hausdorff convergence and scalar convergence (see Proposition~\ref{conv-finite}), we can show that, if $\o\in\ol{\O}$, then for every $x^\ast\in\sW$ it holds that
\beaa
s^{x^\ast}_{t-}(\o)= \lim_{r\in\hQ,r\uparrow t}s(x^\ast,M_r(\o))=\lim_{n\rightarrow\infty}s(x^\ast,M_{r_n(t)}(\o))
=\lim_{k\rightarrow\infty}s(x^\ast,M_{r_{n_k}(t)}(\o))=s(x^\ast,A_{t}(\o)). 
\eeaa

Finally, we note that if $(\tilde{r}_n(t))_{n\in\hN}\subset \hQ$ is another sequence such that $\tilde{r}_n(t)\uparrow t$ as $n\rightarrow\infty$, then the same argument shows that, for $\o\in\ol{\O}$, and every $x^\ast\in\sW$, we have
\[
\lim_{n\rightarrow\infty}s(x^\ast,M_{\tilde{r}_n(t)}(\o))=s^{x^\ast}_{t-}(\o)=\lim_{n\rightarrow\infty}s(x^\ast,M_{r_n(t)}(\o))=s(x^\ast,A_t(\o)).
\]
Since Hausdorff convergence and scalar convergence are equivalent on $\sK(\hR^d)$ (see Proposition~\ref{conv-finite}(iii)),  we conclude that $(M_{\tilde{r}_n(t)}(\o))_{n\in\hN}$ Hausdorff converges to $A_t(\o)$. In other words, the limit is free of the choice of the rational sequence converging to $t$, and thus  $M_{t-} =A_t $, $\hP$-a.s. 

An almost identical argument also shows that $M_{t+} =\tilde{A}_t\in\hL^0(\O, \sK(\hR^d))$ exists, we leave it to the interested reader.
\qed

\subsection{Existence of Right-Continuous Modifications}\label{sec:reg}

In light of the vector-valued martingale theory, we now try to turn the  temporally sectional result Lemma \ref{leftlimited} into the first path-regularity result for the set-valued submartingales. As before, although the arguments mostly follow the standard stochastic analysis, we would like to emphasize the special properties based on the nature of set-valued (multi)functions. 

Let $(M_t)_{t\geq 0}$ be a set-valued $\hF$-submartingale. Let  $(M_{t-})_{t>0}$, $(M_{t+})_{t\ge0}$ be the families of set-valued random variables in Lemma~\ref{leftlimited}, defined on the $\hP$-a.s. set $\bar{\O}$; and set $M_{t-}(\o)=M_{t+}(\o):=\{0\}$,   $\o\in\bar{\O}^c$, for definiteness.  Then, we have the following result.

\begin{thm}\label{prop:cadlag}
	Suppose that the filtration $\hF=(\cF_t)_{t\geq 0}$ satisfies the usual conditions. Let $(M_t)_{t\geq 0}$ be a set-valued convex and integrably bounded set-valued $\hF$-submartingale.
		
	(i) For each $t\geq 0$, we have $M_{t+}\in\scL^1_{\cF_{t}}(\O,\sK(\hR^d))$ and $M_t\subseteq M_{t+}$ $\hP$-a.s.; moreover, $M_t=M_{t+}$ $\hP$-a.s. if and only if the function $r\mapsto s(x^\ast,\hE[M_r])$ is right-continuous at $r=t$.
	
	(ii) The process $(M_{t+})_{t\geq 0}$ is a Hausdorff c\`{a}dl\`{a}g (right-continuous and left-limited), convex, and integrably bounded set-valued $\hF$-submartingale. If we further assume that $(M_t)_{t\geq 0}$ is a set-valued $\hF$-martingale, then so is $(M_{t+})_{t\geq 0}$.
	\end{thm}
	
{\it Proof}: Let $t\geq 0$. First note that, since $\hP(\bar{\O})=1$,  we have $\bar{\O}\in\cF_t$ by the augmentedness of $\hF$. Thus, it follows from Lemma~\ref{leftlimited} that $M_{t+}$ is $\cF_{t+}$-measurable, whence $\cF_t$-measurable  by the right-continuity of $\hF$. We now follow our ``scalarization" scheme to carry out the proof.

(i) {\it Step 1.} Let $t\geq 0$ and $(r_n(t))_{n\in\hN}\subset \hQ$ be such that $r_n(t)\downarrow t$ as $n\rightarrow\infty$. Let $\hG_-=(\cG_{-n}:=\cF_{r_n(t)})_{n\in\hN}$. Then, $(M_{r_n(t)})_{n\in\hN}$ is a convex integrably bounded $\hG_-$-submartingale. By Theorem~\ref{thm:reverse}, the sequence $(\norm{M_{r_n(t)}})_{n\in\hN}$ is uniformly integrable, $(M_{r_n(t)})_{n\in\hN}$ Hausdorff converges to $M_{t+}$ in $\hL^1$, and the limit $M_{t+}$ is convex and integrably bounded. Hence, $M_{t+}\in\scL^1_{\cF_t}(\O,\sK(\hR^d))$.

\ss
{\it Step 2.} For each $x^\ast\in\sW$, we define $s^{x^\ast}_r:=s(x^\ast,M_r)$, $r\geq 0$. Then $(s^{x^\ast}_r)_{r\geq 0}$ is a real-valued $\hF$-submartingale by Remark~\ref{mtg-ex}(iii). Hence, by \cite[Theorem V.4.19]{cinlar}, we have
\bea\label{support-order}
s^{x^\ast}_{t}\leq s^{x^\ast}_{t+}\quad \hP\as
\eea 
On the other hand, the  Hausdorff convergence of $(M_{r_n(t)})_{n\in\hN}$ is Step 1 implies that
\bea\label{support-limit}
s^{x^\ast}_{t+}=\lim_{n\rightarrow\infty}s^{x^\ast}_{r_n(t)}=s\big(x^\ast,\lim_{n\rightarrow\infty}M_{r_n(t)}\big)=s(x^\ast,M_{t+})\quad\hP\text{-a.s.},
\eea
and \eqref{support-order} reads as $s(x^\ast,M_t)\leq s(x^\ast,M_{t+})$, $\hP$-a.s. Since $\sW\subset\hB_{\hR^d}$ is dense, the continuity of the support functions yields $s(x^\ast,M_t)\leq s(x^\ast,M_{t+})$, for all $x^\ast\in \hB_{\hR^d}$, $\hP$-a.s. To wit, $M_t\subseteq M_{t+}$ $\hP$-a.s.

\ss
{\it Step 3.} Continuing from Step 2, we see that $M_t=M_{t+}$ $\hP$-a.s. if and only if $s(x^\ast,M_t)=s(x^\ast,M_{t+})$ $\hP$-a.s. for every $x^\ast\in \sW$. By \cite[Theorem V.4.19]{cinlar}, the latter condition is equivalent to the right-continuity of $r\mapsto \hE[s(x^\ast,M_r)]$ at $r=t$ for every $x^\ast\in\sW$. Since $\hE[s(x^\ast,M_r)]=s(x^\ast,\hE[M_r])$ by \cite[Theorem~2.1.35]{Molchanov}, the second part of (i) follows.

\ss
(ii) We first note that the integrable boundedness, $\hF$-adaptedness, and convexity are already argued in (i), and the path regularity is by construction. Thus, it remains to check the submartingale property of $(M_{r+})_{r\geq 0}$. To this end, let $0\leq u< t$ and $x^\ast\in \sW$. We claim that
\bea \label{scalar-submtg}
s(x^\ast,M_{u+})\leq \hE[s(x^\ast,M_{t+})|\cF_u]\quad \hP\text{-a.s.}
\eea 
By Step 1 of (i), $s(x^\ast,M_{u+})$ is $\cF_u$-measurable. Hence, to show \eqref{scalar-submtg}, we fix an arbitrary $A\in \cF_u$ and show that $\hE[1_A s(x^\ast,M_{u+})]\leq \hE[1_A \hE[s(x^\ast,M_{t+})|\cF_u]]$, or equivalently, $\hE[1_A s(x^\ast,M_{u+})]\leq \hE[1_A s(x^\ast,M_{t+})]$. Let $(r_n(u))_{n\in\hN}$, $(r_n(t))_{n\in\hN}$ be sequences in $\hQ$ such that $u<r_n(u)<t<r_n(t)$ for every $n\in\hN$ and $r_n(u)\downarrow u$, $r_n(t)\downarrow t$ as $n\rightarrow\infty$. Note that $(M_{r_n(t)})_{n\in\hN}$ Hausdorff converges to $M_{t+}$ in $\hL^1$ and, by \eqref{H-S}, for each $n\in\hN$, we have
\[
\abs{s(x^\ast,M_{r_n(t)})-s(x^\ast,M_{t+})}\leq h(M_{r_n(t)},M_{t+})\quad\hP\text{-a.s.}
\]
Hence, $(s(x^\ast,M_{r_n(t)}))_{n\in\hN}$ converges to $s(x^\ast,M_{t+})$ in $\hL^1$. Similarly, $(s(x^\ast,M_{r_n(u)}))_{n\in\hN}$ converges to $s(x^\ast,M_{u+})$ in $\hL^1$. In particular, we get
\beaa
\hE[1_A s(x^\ast,M_{u+})]&=&\lim_{n\rightarrow\infty}\hE[1_A s(x^\ast,M_{r_n(u)})]\le \lim_{n\rightarrow\infty}\hE[1_A \hE[s(x^\ast,M_{r_n(t)})| \cF_u]]\\
&=& \lim_{n\rightarrow\infty}\hE[1_A s(x^\ast,M_{r_n(t)})]  = \hE[1_A s(x^\ast,M_{t+})]\quad \hP\text{-a.s.},
\eeaa
where the inequality is by the submartingale property of $(s(x^\ast,M_r))_{r\geq 0}$. Therefore, \eqref{scalar-submtg} follows. Finally, since $s(x^\ast, \hE[M_{t+}|\cF_u])=\hE[s(x^\ast,M_{t+})|\cF_u]$ by \cite[Theorem~2.1.72]{Molchanov}, we conclude that $s(x^\ast,M_{u+})\leq s(x^\ast, \hE[M_{t+}|\cF_u])$ for every $x^\ast\in \hB_{\hR^d}$ $\hP$-a.s., that is, $M_{u+}\subseteq \hE[M_{t+}|\cF_u]$ $\hP$-a.s. Hence, $(M_{r+})_{r\geq 0}$ is an $\hF$-submartingale. When $(M_r)_{r\geq 0}$ is an $\hF$-martingale, the inequality in \eqref{support-order} becomes an equality so that $(M_{r+})_{r\geq 0}$ is an $\hF$-martingale.
\qed 

As an immediate consequence of Theorem~\ref{prop:cadlag}, we extend a well-known theorem for real-valued martingales to the set-valued setting: \emph{Every martingale has a c\`{a}dl\`{a}g modification.}

\begin{cor}\label{thm:cadlag}
	Suppose that $\hF$ is a filtration that satisfies the usual conditions. Let $(M_t)_{t\geq 0}$ be a convex and integrably bounded set-valued $\hF$-submartingale such that $t\mapsto s(x^\ast,\hE[M_t])$ is right-continuous for every $x^\ast\in\sW$. Then, it has a Hausdorff c\`{a}dl\`{a}g modification $(\tilde{M}_t)_{t\geq 0}$, i.e., $(\tilde{M}_t)_{t\geq 0}$ is a convex and integrably bounded set-valued $\hF$-submartingale with $\hP\{M_t=\tilde{M}_t\}=1$ for every $t\geq 0$ and the following properties hold with probability one:
	
	(a) $\tilde{M}_{t-}$ exists for every $t>0$,
	
	(b) $\tilde{M}_{t+}$ exists and $\tilde{M}_t=\tilde{M}_{t+}$ for every $t\geq 0$.
	
	If we further assume that $(M_t)_{t\geq 0}$ is a set-valued $\hF$-martingale, then so is $(\tilde{M}_t)_{t\geq 0}$.
	\end{cor}
	
{\it Proof}: Let $\tilde{M}_t=M_{t+}$ for each $t\geq 0$. Then, all the claimed properties for $(\tilde{M}_t)_{t\geq 0}$ follow directly from Theorem~\ref{prop:cadlag}.
\qed 

Another important fact in stochastic analysis is that all martingales have a continuous modification when the filtration $\hF$ is {\it Brownian}, thanks to the martingale representation theorem and the continuity of stochastic integrals. Although the direct imitation does not work in the set-valued case, we are nevertheless able to prove  the following result for set-valued martingales, by again taking advantage of the \emph{scalarization} procedure through support functions.

\begin{thm}\label{cor:cont}
	Let $(B_t)_{t\geq 0}$ be a standard Brownian motion with values in $\hR^m$, where $m\in\hN$. Suppose that $\hF$ is the natural filtration generated by $B$ and it satisfies the usual conditions. Let $(M_t)_{t\geq 0}$ be a convex and integrably bounded set-valued $\hF$-martingale. Then, it has a Hausdorff continuous modification $(\tilde{M}_t)_{t\geq0}$, i.e., $(\tilde{M}_t)_{t\geq 0}$ satisfies the properties in Corollary~\ref{thm:cadlag} and it has Hausdorff continuous paths with probability one.
\end{thm}

{\it Proof}: Since $(M_t)_{t\ge0}$ is an $\hF$-martingale, the mapping  $t\mapsto s(x^\ast,\hE[M_t])$ is a constant, for every $x^\ast\in\sW$. By Corollary~\ref{thm:cadlag},  $(M_t)_{t\geq0}$ has a Hausdorff c\`{a}dl\`{a}g modification $(\tilde{M}_t)_{t\geq 0}$. Next, for fixed $x^\ast\in\sW$, consider $\tilde{s}^{x^\ast}_t:=s(x^\ast,\tilde{M}_t)$, $t\geq 0$. Since $(\tilde{M}_t)_{t\geq 0}$ is Hausdorff c\`{a}dl\`{a}g, the real-valued $\hF$-martingale $(\tilde{s}^{x^\ast}_t)_{t\geq 0}$ is c\`{a}dl\`{a}g. By the classical martingale representation theorem, $(\tilde{s}^{x^\ast}_t)_{t\geq 0}$ has a continuous modification $(\bar{s}^{x^\ast}_t)_{t\geq 0}$. Since both processes are determined by their values at rational times, they are indeed indistinguishable. Hence, $(\tilde{s}^{x^\ast}_t)_{t\geq 0}$ has continuous paths $\hP$-a.s.

Finally, let $t>0$ with $t\in\hQ$ and $x^\ast\in\sW$. Since $(\tilde{s}^{x^\ast}_t)_{t\geq 0}$ has continuous paths, we have
\[
s(x^\ast,\tilde{M}_t)=\tilde{s}^{x^\ast}_t=\tilde{s}^{x^\ast}_{t-}=\lim_{r\in\hQ,r\uparrow t}s(x^\ast,\tilde{M}_r)=s(x^\ast,\tilde{M}_{t-})\quad \hP\text{-a.s.}
\]
Thus, $\tilde{M}_t=\tilde{M}_{t-}$, $t\in\hQ\cap(0,\infty)$, $\hP$-a.s. Since $(\tilde{M}_{t})_{t\geq 0}$ is Hausdorff c\`{a}dl\`{a}g, this implies that $\tilde{M}_t=\tilde{M}_{t+}=\lim_{r\in\hQ,r\downarrow t}\tilde{M}_r=\lim_{r\in\hQ,r\downarrow t}\tilde{M}_{r-}=\tilde{M}_{t-}$, for all $t>0$, $ \hP$-a.s. That is, $(\tilde{M}_{t})_{t\geq 0}$ has Hausdorff continuous paths with probability one.
\qed

\subsection{Existence of Left-Continuous Modifications}\label{sec:reg2}

Complementing the discussion of Section~\ref{sec:reg}, we now study the existence of left-continuous modifications of set-valued submartingales, which turns out to be slightly more involved. 

We begin by the following observation in the  real-valued martingale theory: although under the usual conditions, every real-valued martingale has a c\`{a}dl\`{a}g modification, the left-continuous analogue does not hold in full generality when one replaces the right-continuity of the filtration by left-continuity.

\begin{lem}\label{lem:scalarleftcont}
	Suppose that $\hF=(\cF_t)_{t\geq 0}$ is an augmented and left-continuous filtration, that is, $\cF_t=\cF_{t-}:=\sigma(\bigcup_{r<t}\cF_r)$ for each $t>0$. Let $(\eta_t)_{t\geq 0}$ be a real-valued $\hF$-submartingale.
	
	(i) For each $t>0$, we have $\eta_{t-}\leq\eta_t$ $\hP$-a.s.; moreover, $\eta_{t-}=\eta_t$ $\hP$-a.s. if and only if $r\mapsto \hE[\eta_r]$ is left-continuous at $r=t$.
	
	(ii) Suppose that $\eta_t\geq 0$ $\hP$-a.s. for each $t\geq 0$. Then, we have $\eta_{t-}\in\hL^1_{\cF_t}(\O,\hR)$ for each $t>0$. Moreover, $(\eta_{t-})_{t\geq 0}$ is a c\`{a}gl\`{a}d (left-continuous and right-limited) real-valued $\hF$-submartingale, where $\eta_{0-}:=\eta_0$.
\end{lem}

{\it Proof}: (i) Let $t>0$. By \cite[Proposition~V.7.14]{cinlar}, there exists a set $\bar{\O}\in\cF$ such that $\hP(\bar{\O})=1$ and, for every $\o\in\bar{\O}$, the limits
$\eta_{t-}:=\lim_{r\in\hQ,r\uparrow t}\eta_r(\o)$, $t>0$, and $\eta_{t+}:=\lim_{r\in\hQ,r\downarrow t}\eta_r(\o)$, $t\geq 0$, both
exist and are finite. Again, we set $\eta_{t-}(\o):=0$, $t>0$, and $\eta_{t+}(\o):=0$, $t\geq0$, for every $\o\in\bar{\O}^c$. Since $\hF$ is augmented, we have $\bar{\O}\in\cF_0\subseteq\cF_{t-}\subseteq\cF_t$. Hence, the definition of $\eta_{t-}$ implies that it is $\cF_{t-}$-measurable. Since we assume that $\cF_{t-}=\cF_t$, we conclude that $\eta_{t-}$ is $\cF_t$-measurable.

Let $(r_n(t))_{n\in\hN}\subset \hQ$ be such that $r_n(t)\uparrow t$ as $n\rightarrow\infty$. By the submartingale property, we have
\bea\label{eq:submtg-n}
\eta_{r_n(t)}\leq \hE[\eta_t|\cF_{r_n(t)}]\quad\hP\text{-a.s.}, \q n\in\hN.
\eea
Next note that $\cF_{t-}=\bigvee_{n\in\hN}\cF_{r_n(t)}$ and $\eta_t$ is  $\cF_{t-}$-measurable. By L\'{e}vy's upward theorem, we have
\[
\lim_{n\rightarrow\infty}\hE[\eta_t|\cF_{r_n(t)}]=\hE[\eta_t|\cF_{t-}]=\eta_t  \quad\hP\text{-a.s.}
\]
Since $\lim_{n\to\infty}\eta_{r_n(t)}=\eta_{t-}$ $\hP$-a.s., letting $n\rightarrow\infty$ in \eqref{eq:submtg-n} yields $\eta_{t-}\leq \eta_t$ $\hP$-a.s. In particular, $\eta_{t-}=\eta_t$ $\hP$-a.s.~if and only if $\hE[\eta_{t-}]=\hE[\eta_t]$, that is, $r\mapsto \hE[\eta_r]$ is left-continuous at $r=t$.

(ii) Let $t, (r_n(t))_{n\in\hN}$ be as in (i). Since $\eta_t$ is integrable, the family $(\hE[\eta_t|\cF_{r_n(t)}])_{n\in\hN}$ is uniformly integrable. Then, the positivity assumption and \eqref{eq:submtg-n} imply that the family $(\eta_{r_n(t)})_{n\in\hN}$ is also uniformly integrable. Hence, by the real-valued submartingale convergence theorem \cite[Theorem~V.4.5]{cinlar}, $(\eta_{r_n(t)})_{n\in\hN}$ converges to $\eta_{t-}$, both $\hP$-a.s. and in $\hL^1$. In particular, $\eta_{t-}\in\hL^1_{\cF_t}(\O,\hR)$.

The $\hF$-adaptedness and integrability of $(\eta_{r-})_{r\geq 0}$ are already established. It remains to check the submartingale property. Let $0\leq u<t$. We show that $\eta_{u-}\leq \hE[\eta_{t-}|\cF_u]$ $\hP$-a.s., which is equivalent to having $\hE[\eta_{u-}{\bf 1}_A]\leq \hE[\hE[\eta_{t-}|\cF_u]{\bf 1}_A]$, that is,
\bea\label{eq:submtg-}
\hE[\eta_{u-}{\bf 1}_A]\leq \hE[\eta_{t-}{\bf 1}_A],  \quad A\in\cF_u.
\eea

First, suppose that $u>0$. Let us fix $s<u$ and take $A\in\cF_s$. Let $(r_n(u))_{n\in\hN}$, $(r_n(t))_{n\in\hN}$ be sequences in $\hQ$ such that $s<r_n(u)<u<r_n(t)<t$ for every $n\in\hN$ and $r_n(u)\uparrow u$, $r_n(t)\uparrow t$ as $n\rightarrow\infty$. Since $(\eta_{r_n(u)})_{n\in\hN}$ converges to $\eta_{u-}$ in $\hL^1$ and $(\eta_{r_n(t)})_{n\in\hN}$ converges to $\eta_{t-}$ in $\hL^1$, we have
\[
\hE[\eta_{u-}{\bf 1}_A]=\lim_{n\rightarrow\infty}\hE[\eta_{r_n(u)}{\bf 1}_A]\leq\lim_{n\rightarrow\infty}\hE[\hE[\eta_{r_n(t)}|\cF_{r_n(u)}]{\bf 1}_A]= \lim_{n\rightarrow\infty}\hE[\eta_{r_n(t)}{\bf 1}_A]=\hE[\eta_{t-}{\bf 1}_A],
\]
where the first and last equalities follow since (strong) convergence in $\hL^1$ implies weak convergence in $\hL^1$, the inequality follows by the submartingale property of $(\eta_r)_{r\geq 0}$, and the second equality follows since $A\in\cF_s\subseteq\cF_{r_n(u)}$. Hence, \eqref{eq:submtg-} holds for every $A\in \cF_s$. It follows that \eqref{eq:submtg-} holds for every $A\in \bigcup_{s<u}\cF_s$. Let $\sA$ be the collection of all sets $A\in\cF$ for which \eqref{eq:submtg-} holds. As an immediate consequence of monotone convergence theorem, $\sA$ is closed under increasing and decreasing sequences of sets. Hence, it is a monotone class. Since $\bigcup_{s<t}\cF_s$ is an algebra of sets, by monotone class theorem, we get $\cF_u=\cF_{u-}\subseteq\sA$. Therefore, \eqref{eq:submtg-} holds for every $A\in\cF_u$ when $u>0$. When $u=0$, the same argument works in a simplified way by taking $r_n(u)=0$ for each $n\in\hN$. Hence, the proof of the submartingale property of $(\eta_{r-})_{r\geq 0}$ is complete.

Finally, suppose that $(\eta_r)_{r\geq 0}$ is an $\hF$-martingale. Then, $r\mapsto \hE[\eta_r]$ is a constant function and we have $\eta_{r-}=\eta_{r}$ $\hP$-a.s. for each $r\geq 0$ by (i). Hence, $(\eta_{r-})_{r\geq 0}$ is an $\hF$-martingale.
\qed 

Despite the lack of a general integrability result for the left-limits of a real-valued submartingale, we are able to prove the following left-continuous analog of Theorem~\ref{prop:cadlag} for set-valued submartingales. As we will see in the proof, this is possible by exploiting the uniform integrability of the positive submartingale given by the ``set norm" of the given set-valued submartingale.

\begin{thm}\label{prop:caglad}
	Suppose that $\hF=(\cF_t)_{t\geq 0}$ is an augmented and left-continuous filtration. Let $(M_t)_{t\geq 0}$ be a set-valued convex and integrably bounded set-valued $\hF$-submartingale. Let $\bar{\O}\in\cF$ be the $\hP$-a.s. set in Lemma~\ref{leftlimited}. Set $M_{0-}:=M_0$. For definiteness, set $M_{t-}(\o):=\{0\}$, $t>0$, and $M_{t+}(\o):=\{0\}$, $t\geq 0$, for every $\o\in\bar{\O}^c$. Then, we have the following results.
	
	(i) For each $t>0$, we have $M_{t-}\in\scL^1_{\cF_{t}}(\O,\sK(\hR^d))$ and $M_{t-}\subseteq M_{t}$ $\hP$-a.s.; moreover, $M_{t-}=M_{t}$ $\hP$-a.s. if and only if the function $r\mapsto s(x^\ast,\hE[M_r])$ is left-continuous at $r=t$.
	
	(ii) The process $(M_{t-})_{t\geq 0}$ is a Hausdorff c\`{a}gl\`{a}d, convex, and integrably bounded set-valued $\hF$-submartingale. Furthermore, if   $(M_t)_{t\geq 0}$ is a set-valued $\hF$-martingale, then so is $(M_{t-})_{t\geq 0}$.
\end{thm}

{\it Proof}: Let $t> 0$. Since $\bar{\O}\in\cF$ is a $\hP$-a.s. set, it belongs to $\cF_0\subseteq\cF_t$ by the augmentedness of $\hF$. This and Lemma~\ref{leftlimited} imply that $M_{t-}$ is measurable with respect to $\cF_{t}$. Since $\cF_{t-}=\cF_{t}$ by the left-continuity of $\hF$, we conclude that $M_{t-}$ is $\cF_t$-measurable.

(i) {\it Step 1.} Let $t> 0$ and $(r_n(t))_{n\in\hN}$ be a sequence in $\hQ$ such that $r_n(t)\uparrow t$ as $n\rightarrow\infty$. Let $\hG=(\cG_n:=\cF_{r_n(t)})_{n\in\hN}$. Note that $(M_{r_n(t)})_{n\in\hN}$ is a convex integrably bounded $\hG$-submartingale in discrete time. Note that, by the submartingale property of $(\norm{M_r})_{r\geq 0}$, we have $0\leq \norm{M_{r_n(t)}}\leq \hE[\norm{M_t}|\cF_{r_n(t)}]$ $\hP\text{-a.s.}$ for each $n\in\hN$. Since $\norm{M_t}\in\hL^1$, the family $\{\hE[\norm{M_t}|\cF_{r_n(t)}]\colon n\in\hN\}$ is uniformly integrable. It follows that the family $\{\norm{M_{r_n(t)}}\colon n\in\hN\}$ is also uniformly integrable. Hence, by Theorem~\ref{dtmtg:findim}, $(M_{r_n(t)})_{n\in\hN}$ converges to $M_{t-}$ in the Hausdorff sense in $\hL^1$ (in addition to $\hP$-a.s.) and we have $M_{t-}\in\scL^1_{\cG_{\infty}}(\O,\sK(\hR^d))$. Here, note that $\cG_\infty=\cF_{t-}=\cF_t$. Hence, $M_{t-}\in\scL^1_{\cF_t}(\O,\sK(\hR^d))$.

{\it Step 2.} Let $x^\ast\in\sW$ and define $s^{x^\ast}_r:=s(x^\ast,M_r)$ for each $r\geq 0$. Note that $(s^{x^\ast}_r)_{r\geq 0}$ is a real-valued $\hF$-submartingale by Remark~\ref{mtg-ex}(iii). Hence, by Lemma~\ref{lem:scalarleftcont}, we have
\bea\label{support-order2}
s^{x^\ast}_{t-}\leq s^{x^\ast}_{t}\quad \hP\as
\eea 
On the other hand, for every sequence $(r_n(t))_{n\in\hN}$ in $\hQ$ with $r_n(t)\uparrow t$ as $n\rightarrow\infty$, we have
\bea\label{support-limit2}
s^{x^\ast}_{t-}=\lim_{n\rightarrow\infty}s^{x^\ast}_{r_n(t)}=s\of{x^\ast,\lim_{n\rightarrow\infty}M_{r_n(t)}}=s(x^\ast,M_{t-})\quad\hP\text{-a.s.},
\eea 
where the second equality follows since Hausdorff convergence implies scalar convergence. Hence, \eqref{support-order2} reads as $s(x^\ast,M_{t-})\leq s(x^\ast,M_{t})$ $\hP$-a.s. Since $\sW$ is dense in $\hB_{\hR^d}$ and the support functions are continuous in $x^\ast$, we conclude that $s(x^\ast,M_{t-})\leq s(x^\ast,M_{t})$ holds for every $x^\ast\in \hB_{\hR^d}$ with probability one. Equivalently, $M_{t-}\subseteq M_{t}$ $\hP$-a.s.

{\it Step 3.} By the reasoning in Step 2, we also have that $M_{t-}=M_{t}$ $\hP$-a.s. if and only if $s(x^\ast,M_{t-})=s(x^\ast,M_{t})$ $\hP$-a.s. for every $x^\ast\in \sW$. By Lemma~\ref{lem:scalarleftcont}, the latter condition is equivalent to the left-continuity of $r\mapsto \hE[s(x^\ast,M_r)]$ at $r=t$ for every $x^\ast\in\sW$. Since $\hE[s(x^\ast,M_r)]=s(x^\ast,\hE[M_r])$ by \cite[Theorem~2.1.35]{Molchanov}, the second part of (i) follows.

(ii) Integrable boundedness, $\hF$-adaptedness, and convexity are already shown in (i). Path regularity is by construction. It remains to prove the submartingale property of $(M_{r-})_{r\geq 0}$. Let $0\leq u< t$ and $x^\ast\in \sW$. We claim that
\bea \label{scalar-submtg2}
s(x^\ast,M_{u-})\leq \hE[s(x^\ast,M_{t-})|\cF_u]\quad \hP\text{-a.s.}
\eea 
By Step 1 of (i), $s(x^\ast,M_{u-})$ is $\cF_u$-measurable. Hence, to show \eqref{scalar-submtg2}, we fix an arbitrary $A\in \cF_u$ and show that $\hE[1_A s(x^\ast,M_{u-})]\leq \hE[1_A \hE[s(x^\ast,M_{t-})|\cF_u]]$, or equivalently, $\hE[1_A s(x^\ast,M_{u-})]\leq \hE[1_A s(x^\ast,M_{t-})]$. Let $(r_n(u))_{n\in\hN}$, $(r_n(t))_{n\in\hN}$ be sequences in $\hQ$ such that $r_n(u)<u<r_n(t)<t$ for every $n\in\hN$ and $r_n(u)\uparrow u$, $r_n(t)\uparrow t$ as $n\rightarrow\infty$. Note that $(M_{r_n(t)})_{n\in\hN}$ Hausdorff converges to $M_{t-}$ in $\hL^1$ and, by \eqref{H-S}, for each $n\in\hN$, we have
\[
\abs{s(x^\ast,M_{r_n(t)})-s(x^\ast,M_{t-})}\leq h(M_{r_n(t)},M_{t-})\quad\hP\text{-a.s.}
\]
Hence, $(s(x^\ast,M_{r_n(t)}))_{n\in\hN}$ converges to $s(x^\ast,M_{t-})$ in $\hL^1$. Similarly, $(s(x^\ast,M_{r_n(u)}))_{n\in\hN}$ converges to $s(x^\ast,M_{u-})$ in $\hL^1$. In particular, we get
\beaa
\hE[1_A s(x^\ast,M_{u-})]&=&\lim_{n\rightarrow\infty}\hE[1_A s(x^\ast,M_{r_n(u)})]\le \lim_{n\rightarrow\infty}\hE[1_A \hE[s(x^\ast,M_{r_n(t)})| \cF_u]]\\
&=& \lim_{n\rightarrow\infty}\hE[1_A s(x^\ast,M_{r_n(t)})]  = \hE[1_A s(x^\ast,M_{t-})]\quad \hP\text{-a.s.},
\eeaa
where the inequality is by the submartingale property of $(s(x^\ast,M_r))_{r\geq 0}$. Therefore, \eqref{scalar-submtg2} follows. Finally, since $s(x^\ast, \hE[M_{t-}|\cF_u])=\hE[s(x^\ast,M_{t-})|\cF_u]$ by \cite[Theorem~2.1.72]{Molchanov}, we conclude that $s(x^\ast,M_{u-})\leq s(x^\ast, \hE[M_{t-}|\cF_u])$ for every $x^\ast\in \hB_{\hR^d}$ $\hP$-a.s., that is, $M_{u-}\subseteq \hE[M_{t-}|\cF_u]$ $\hP$-a.s. Hence, $(M_{r-})_{r\geq 0}$ is an $\hF$-submartingale. When $(M_r)_{r\geq 0}$ is an $\hF$-martingale, the inequality in \eqref{support-order2} becomes an equality so that $(M_{r-})_{r\geq 0}$ is an $\hF$-martingale.
\qed

The next two corollaries complement Corollary~\ref{thm:cadlag}.

\begin{cor}\label{thm:caglad}
	Suppose that $\hF$ is a filtration that is augmented and left-continuous. Let $(M_t)_{t\geq 0}$ be a convex and integrably bounded set-valued $\hF$-submartingale such that $t\mapsto s(x^\ast,\hE[M_t])$ is left-continuous for every $x^\ast\in\sW$. Then, it has a Hausdorff c\`{a}gl\`{a}d modification $(\tilde{M}_t)_{t\geq 0}$, i.e., $(\tilde{M}_t)_{t\geq 0}$ is a convex and integrably bounded set-valued $\hF$-submartingale with $\hP\{M_t=\tilde{M}_t\}=1$ for every $t\geq 0$ and the following properties hold with probability one:
	
	(a) $\tilde{M}_{t-}$ exists and $\tilde{M}_{t-}=\tilde{M}_t$ for every $t>0$.
	
	(b) $\tilde{M}_{t+}$ exists for every $t\geq 0$.
	
	If we further assume that $(M_t)_{t\geq 0}$ is a set-valued $\hF$-martingale, then so is $(\tilde{M}_t)_{t\geq 0}$.
\end{cor}

{\it Proof}: Let $\tilde{M}_t=M_{t-}$ for each $t\geq 0$. Then, all the claimed properties for $(\tilde{M}_t)_{t\geq 0}$ follow directly from Theorem~\ref{prop:caglad}.
\qed 

\begin{cor}\label{thm:cont}
	Suppose that $\hF$ is a filtration that is left-continuous and satisfies the usual conditions. Let $(M_t)_{t\geq 0}$ be a convex and integrably bounded set-valued $\hF$-submartingale such that $t\mapsto s(x^\ast,\hE[M_t])$ is continuous for every $x^\ast\in\sW$. 
	
	(i) Then, it has a Hausdorff continuous modification $(\tilde{M}_t)_{t\geq 0}$, i.e., $(\tilde{M}_t)_{t\geq 0}$ is a convex and integrably bounded set-valued $\hF$-submartingale with $\hP\{M_t=\tilde{M}_t\}=1$, $t\geq 0$, and the following property holds: $\tilde{M}_{t-}, \tilde{M}_{t+}$ exist and $\tilde{M}_{t-}=\tilde{M}_t=\tilde{M}_{t+}$ $\hP$-a.s. for every $t\geq 0$. (Here, $\tilde{M}_{0-}:=\tilde{M}_0$.)
	
	(ii) If we further assume that $(M_t)_{t\geq 0}$ is a set-valued $\hF$-martingale, then so is $(\tilde{M}_t)_{t\geq 0}$.
\end{cor}

{\it Proof}: Let $\bar{M}_t:=M_{t+}$ for each $t\geq 0$ as in the proof of Corollary~\ref{thm:cadlag}. Since $t\mapsto s(x^\ast,\hE[M_t])$ is right-continuous for every $x^\ast\in\sW$, the process $(\bar{M}_t)_{t\geq 0}$ is a Hausdorff c\`{a}dl\`{a}g modification of $(M_t)_{t\geq 0}$. Let $\tilde{M}_t:=\bar{M}_{t-}$ for each $t\geq 0$ as in the proof of Corollary~\ref{thm:caglad}. Note that $t\mapsto s(x^\ast,\hE[\bar{M}_{t}])=s(x^\ast,\hE[M_t])$ is left-continuous for every $x^\ast\in\sW$. Hence, the process $(\tilde{M}_t)_{t\geq 0}$ is a Hausdorff c\`{a}gl\`{a}d modification of $(\bar{M}_t)_{t\geq 0}$. Moreover, since $(\bar{M}_t)_{t\geq 0}$ is Hausdorff right-continuous, so is $(\tilde{M}_t)_{t\geq 0}$. It follows that $(\tilde{M}_t)_{t\geq 0}$ is a Hausdorff continuous modification of $(M_t)_{t\geq 0}$.
\qed 

\section{Set-Valued Stochastic Integrals and Martingale Properties}\label{sec:integral}
\setcounter{equation}{0}

In this section, we study the set-valued stochastic integrals recently introduced in \cite{amw}. 
The main difference between this stochastic integral and the Aumann-It\^o integral is that it allows for non-singleton initial values at $t=0$, and hence particularly useful for the study of set-valued backward stochastic differential equations. We shall first recall the definition and construction of the integral and prove some properties that will be useful for future discussions. Then, we show that, unlike their vector-valued counterparts, set-valued stochastic integrals are set-valued submartingales in general and we characterize those stochastic integrals that are true set-valued martingales.

In what follows we consider a probability space $(\O,\cF,\hP)$ on which is defined a $m$-dimensional Brownian motion $B=(B_t)_{t\in [0,T]}$, where $T>0$ is a given horizon. We assume that $\hF=(\cF_t)_{t\in[0,T]}$ is generated by the Brownian motion $B$ and it satisfies the \emph{usual conditions}.

\subsection{Set-Valued Stochastic Integration}

Let us introduce the space $\sR^{d,m}_{\hF}:=\hR^d\times \hL_{\hF}^2([0,T]\times\O,\hR^{d\times m})$. Given $(x,z)\in \sR^{d,m}_{\hF}$, we define a continuous $\hF$-martingale  $\cJ(x,z)=(\cJ_t(x,z))_{t\in[0,T]}$ by
	\bea\label{J-int}
	\cJ_t(x,z):=x+\int_0^t z_s dB_s,\quad t\in [0,T].
	\eea

\begin{defn}\label{defn:stochint}
	\emph{\cite[Section 5.3]{amw}} Let $\cR\subseteq \sR^{d,m}_{\hF}$ be a nonempty set, $t\in[0,T]$. The stochastic integral of $\cR$ over the interval $[0,t]$ is defined as the $\hP$-a.s. unique set-valued random variable $\int_{0-}^t\cR dB\in \sA^2_{\cF_t}(\O,\sC(\hR^d))$ that satisfies 
\[
S^2_{\cF_t}\of{\int_{0-}^t \cR dB}=\cl_{\hL^2}\dec_{\cF_t}(\cJ_t[\cR]).
\]
\end{defn}

\begin{rem}\label{rem:aumann-ito}
	When $\cR=\{x\}\times\cZ$ for some $x\in \hR^d$ and $\cZ\subseteq\hL^2_{\hF}([0,T]\times \O,\hR^{d\times m})$, we have
\[
\int_{0-}^t \cR dB = \{x\}+\int_0^t \cZ dB,
\]
where $\int_0^t \cZ dB \in \sA^2_{\cF_t}(\O,\sC(\hR^d))$ is the $\hP$-a.s. unique set-valued random variable such that
\[
S^2_{\cF_t}\of{\int_{0}^t \cZ dB}=\cl_{\hL^2}\dec_{\cF_t}\of{\cb{\int_0^t z_s dB_s\colon z\in\cZ}}.
\]
$\int_0^t \cZ dB$ is called the generalized Aumann-Itô stochastic integral; see \cite[\S 5.2]{Kis2020}. In \cite[Theorem~4.2]{MgRT}, it is shown that every set-valued $\hF$-martingale $M$ with $M_0=\{0\}$ can be represented as $M_t=\int_0^t \cZ dB$ $\hP$-a.s. for each $t\in[0,T]$, where $\cZ\subseteq\hL^2_{\hF}([0,T]\times \O,\hR^{d\times m})$ is a nonempty set. However, the condition that $M_0$ is a singleton is quite restrictive; indeed, it is shown in \cite[Lemma~3.1]{ZhangYano} (see also \cite[Proposition~A.1]{hess1999}) that a set-valued random variable whose expectation is a singleton if and only if it is a random singleton. Using a set of pairs $\cR\subseteq \sR^{d,m}_{\hF}$ as an integrand has the advantage of handling multiple initial values simultaneously, hence yielding a more satisfactory integral for the representation of martingales as we will see in Section~\ref{sec:representation}.
\qed
\end{rem}

We prove some basic properties of the set-valued stochastic integral next. These properties are easy to check and they extend their analogs for the generalized Aumann-Itô stochastic integral recalled in Remark~\ref{rem:aumann-ito}.

\begin{prop}\label{lem:basicprop}
	Let $\cR\subseteq \sR^{d,m}_{\hF}$ be a nonempty set and $t\in[0,T]$. The following properties hold:
	
	(i) $\int_{0-}^t \cl_{\hR^d\times\hL^2_{\hF}}(\cR)dB=\int_{0-}^t \cR dB$ $\hP$-a.s.
	
	(ii) $\int_{0-}^t \co(\cR)dB=\cl_{\hR^d}\co(\int_{0-}^t \cR dB)$ $\hP$-a.s. Further, if $\cR$ is convex, then so is $\int_{0-}^t \cR dB$.
	
	(iii) Denoting $\ip{x^\ast,\cR}:=\{(\ip{x^\ast,x},\ip{x^\ast,z})\colon (x,z)\in \cR\}$ for $x^\ast\in \hB_{\hR^d}$, we have
	\[
	s\Big(x^\ast,\int_{0-}^t \cR dB\Big)=\sup\int_{0-}^t \ip{x^\ast,\cR} dB\quad\hP\text{-a.s.}
	\]
	
	(iv) If $\cR$ is countable, say, $\cR=\{(x^i,z^i)\colon i\in\hN\}$, then
	\[
	\int_{0-}^t \cR dB = \cl_{\hR^d}(\{\cJ_t(x^i,z^i)\colon i\in\hN\})\quad \hP\text{-a.s.}
	\]
\end{prop}

{\it Proof}: (i) We first claim that $\cJ_t[\cl_{\hR^d\times\hL^2}(\cR)]\subseteq \cl_{\hL^2}(\cJ_t[\cR])$. Indeed, let $(x^n,z^n)_{n\in\hN}$ be a sequence in $\cR$ that converges to some $(x,z)\in\sR^{d,m}_{\hF}$ in $\hR^d\times\hL^2$. Then,
\beaa
&&\hE\sqb{|\cJ_t(x^n,z^n)-\cJ_t(x,z)|^2}=\hE\sqb{|\cJ_t(x^n-x,z^n-z)|^2}\\
&&\leq  2\of{|x^n-x|^2+\hE\sqb{\Big|\int_0^t(z_s^n-z_s)dB_s\Big|^2}}= 2\of{|x^n-x|^2+\hE\sqb{\int_0^t|z_s^n-z_s|^2ds}}
\eeaa 
for each $n\in\hN$. Hence, $(\cJ_t(x^n,z^n))_{n\in\hN}$ converges to $\cJ_t(x,z)$ and the claim follows. With this and the definition of the set-valued stochastic integral, we have
\[
S^2_{\cF_t}\of{\int_{0-}^t \cl_{\hL^2_{\hF}}(\cR)dB}=\cl_{\hL^2}\dec_{\cF_t}\of{\cJ_t[\cl_{\hL^2}(\cR)]}\subseteq \cl_{\hL^2}\dec_{\cF_t}(\cJ_t[\cR])=S^2_{\cF_t}\of{\int_{0-}^t \cR dB}.
\]
The inclusion $S^2_{\cF_t}(\int_{0-}^t \cl_{\hL^2_{\hF}}(\cR)dB)\supseteq S^2_{\cF_t}(\int_{0-}^t \cR dB)$ is trivial, so these sets are equal, proving (i).

(ii) By the linearity of the mapping $\cJ_t$, we have $\cJ_t[\co(\cR)]=\co(\cJ_t[\cR])$. Then,
\beaa
S^2_{\cF_t}\of{\int_{0-}^t \co(\cR)dB}&=&\cl_{\hL^2}\dec_{\cF_t}\of{\cJ_t[\co(\cR)]}= \cl_{\hL^2}\dec_{\cF_t}\co(\cJ_t[\cR])\\
&=& \cl_{\hL^2}\co \cl_{\hL^2}\dec_{\cF_t}(\cJ_t[\cR])= \cl_{\hL^2}\co\of{S^2_{\cF_t}\of{\int_{0-}^t\cR dB}}\\
&=&S^2_{\cF_t}\of{\cl_{\hL^2}\co\of{\int_{0-}^t \cR dB}},
\eeaa
where the third equality is by \cite[Corollary~7.6]{umur}, the last equality is by \cite[Theorem~2.3.5]{Kis2020}.

(iii) Let $x^\ast\in \hB_{\hR^d}$ and $A\in\cF_t$. First, note that
\beaa
\hE\sqb{s\of{x^\ast,\int_{0-}^t \cR dB}{\bf 1}_A}
&=&\hE\sqb{\sup\cb{\ip{x^\ast,x}\colon x\in \int_{0-}^t \cR dB}{\bf 1}_A}\\
&=&\sup\cb{\hE[\ip{x^\ast,\xi}{\bf 1}_A]\colon \xi\in S^2_{\cF_t}\of{\int_{0-}^t \cR dB}}\\
&=&\sup\cb{\hE[\ip{x^\ast,\xi}{\bf 1}_A]\colon \xi\in \cl_{\hL^2}\dec_{\cF_t}(\cJ_t[\cR])}\\
&=&\sup\cb{\hE[\ip{x^\ast,\xi}{\bf 1}_A]\colon \xi\in \dec_{\cF_t}(\cJ_t[\cR])},
\eeaa 
where the second equality is by \cite[Theorem~2.2]{hiai-umegaki} and the fourth is by the continuity of expectation on $\hL^2$. Then, using \eqref{eq:dec}, we obtain
\beaa 
& &\hE\sqb{s\of{x^\ast,\int_{0-}^t \cR dB}{\bf 1}_A}\\
&&=\sup\cb{\hE\sqb{\ip{x^\ast,\sum_{i=1}^n \cJ_t(x^i,z^i){\bf 1}_{A_i\cap A}}}\colon (x^i,z^i)_{i=1}^n\subseteq\cR, (A_i)_{i=1}^n \in\Pi_{\cF_t}(\O),\ n\in\hN}\\
&&=\sup\cb{\hE\sqb{\sum_{i=1}^n \cJ_t(\ip{x^\ast,x^i},\ip{x^\ast,z^i}){\bf 1}_{A_i\cap A}}\colon (x^i,z^i)_{i=1}^n\subseteq\cR, (A_i)_{i=1}^n\in \Pi_{\cF_t}(\O),\ n\in\hN}\\
&&=\sup\cb{\hE[\xi{\bf 1}_A]\colon \xi\in\dec_{\cF_t}(\cJ_t[\ip{x^\ast,\cR}])}\\
&&=\sup\cb{\hE[\xi{\bf 1}_A]\colon \xi\in S^2_{\cF_t}\of{\int_{0-}^t \ip{x^\ast,\cR}dB}}=\hE\sqb{{\bf 1}_A\sup\int_{0-}^t \ip{x^\ast,\cR}dB},
\eeaa
where the fourth equality is by the continuity of expectation on $\hL^2$ and the fifth is by \cite[Theorem~2.2]{hiai-umegaki}. Since $A\in\cF_t$ is arbitrary, we get $s(x^\ast,\int_{0-}^t \cR dB)=\sup\int_{0-}^t \ip{x^\ast,\cR}$ $\hP$-a.s.

(iv) Let $\cR=\{(x^i,z^i)\colon i\in\hN\}$. We have
\[
S^2_{\cF_t}\of{\int_{0-}^t \cR dB}=\cl_{\hL^2}\dec_{\cF_t}(\{\cJ_t(x^i,z^i)\colon i\in\hN\})=S^2_{\cF_t}\of{\cl_{\hR^d}(\{\cJ_t(x^i,z^i)\colon i\in\hN\})},
\]
where the last equality is by \cite[Theorem~2.3.2]{Kis2020}. Hence, the result follows.
\qed 

\subsection{Integral Representation of Set-Valued Martingales}\label{sec:representation}

Let $M=(M_t)_{t\in[0,T]}$ be a convex set-valued $\hF$-martingale. A $d$-dimensional $\hF$-martingale $y=(y_t)_{t\in[0,T]}$ is called a \emph{martingale selection} of $M$ if $y_t\in S^1_{\cF_t}(M_t)$ for every $t\in[0,T]$. By \cite[Proposition~2]{michta-mtgselection}, $M$ has at least one martingale selection. Moreover, as an $\hF$-martingale, it has a c\`{a}dl\`{a}g (right-continuous and left-limited) modification so that it corresponds to a unique member of $\hL^2_\hF([0,T]\times\O,\hR^d)$. Let us denote by $\MS(M)\neq \emptyset$ the collection of all such members. We also write $P_t(y):=y_t$ whenever $y$ is a vector-valued process and $t\in[0,T]$.

Before we proceed, we first give an important result regarding the connection between the selections of $M_t$ and the $t$-projection of $MS(M)$ at each fixed time $t\in[0,T]$.

\begin{lem}\emph{\cite[Lemma~3.5.3]{Kis2020}}
 \label{lem-kis-Pt}
 Let $M$ be a convex set-valued $\hF$-martingale. Then, for every $t\in[0,T]$, it holds that $S^1_{\cF_t}(M_t)=\cl_{\hL^1}\dec_{\cF_t} (P_t[\MS(M)])$.
\end{lem}

The next proposition characterizes the square-integrable boundedness of $M$.

\begin{prop}\label{mtg-bdd}
	Let $M$ be a convex set-valued $\hF$-martingale. Consider the following conditions:
	
	(i) $M$ is square-integrably bounded.
	
	(ii) $M_T$ is square-integrably bounded.
	
	(iii) $P_T[\MS(M)]$ is a dominated subset of $\hL^2$.
	
	(iv) It holds $P_T[\MS(M)]\subseteq \hL^2$.
	
	(v) For every $y\in \MS(M)$, we have $\hE[\sup_{t\in[0,T]}\abs{y_t}^2]<+\infty$.
	
	We have (i) $\Leftrightarrow$ (ii) $\Leftrightarrow$ (iii) $\Rightarrow$ (iv) $\Leftrightarrow$ (v).
\end{prop}

{\it Proof}: The implications (i) $\Rightarrow$ (ii), (iii) $\Rightarrow$ (iv), and (v) $\Rightarrow$ (iv) are trivial.

(ii) $\Rightarrow$ (i): Suppose that $M_T\in\scL^2_{\cF_T}(\O,\sK(\hR^d))$. By the continuous-time version of Remark~\ref{mtg-ex}, the real-valued process $(\norm{M_t}^2)_{t\in[0,T]}$ is an $\hF$-submartingale. Hence, $\hE[\norm{M_t}^2]\leq \hE[\norm{M_T}^2]<+\infty$ for every $t\in[0,T]$, i.e., $M$ is square-integrably bounded.

(ii) $\Rightarrow$ (iii): Suppose that $M_T\in\scL^2_{\cF_T}(\O,\sK(\hR^d))$. For each $y\in \MS(M)$, since $y_T\in S^2_{\cF_T}(M_T)$, we have $\abs{y_T}\leq \norm{M_T}$ $\hP$-a.s. Then, $M_T\in\scL^2_{\cF_T}(\O,\sC(\hR^d))$ implies that $P_T[\MS(M)]$ is dominated in $\hL^2_{\cF_T}(\O,\hR^d)$.

(iii) $\Rightarrow$ (ii): Suppose that $P_T[\MS(M)]$ is dominated in $\hL^2$. Then, by Lemma~\ref{kis-decomp-bdd}, the sets $\dec_{\cF_T} (P_T[\MS(M)])$, $\cl_{\hL^2}\dec_{\cF_T} (P_T[\MS(M)])$ are dominated and bounded in $\hL^2$. On the other hand, by Lemma~\ref{lem-kis-Pt}, we have $\cl_{\hL^1}\dec_{\cF_T} (P_T[\MS(M)])=S^1_{\cF_T}(M_T)$. We argue that
\bea\label{decL1L2}
\cl_{\hL^1}\dec_{\cF_T} (P_T[\MS(M)])=\cl_{\hL^2}\dec_{\cF_T} (P_T[\MS(M)]).
\eea
Indeed, let $(\zeta^n)_{n\in\hN}$ be a sequence in $\dec_{\cF_T}(P_T[\MS(M)])$ that converges to some $\zeta \in \hL^1_{\cF_T}(\O,\hR^d)$ in $\hL^1$, hence also in probability. Since $\dec_{\cF_T}(P_T[\MS(M)])$ is dominated in $\hL^2$, the sequence $(\abs{\zeta^n}^2)_{n\in\hN}$ is uniformly integrable. Hence, by \cite[Exercise IV.4.13]{cinlar}, $(\zeta^n)_{n\in\hN}$ is convergent in $\hL^2$. By the uniqueness of limits in probability, the $\hL^2$-limit must be $\zeta$. This completes the proof of the inclusion $\subseteq$ in \eqref{decL1L2} and the reverse inclusion is trivial. Hence,
\[
S^1_{\cF_T}(M_T)=\cl_{\hL^2}\dec_{\cF_T} (P_T[\MS(M)])\subseteq \hL^2
\]
so that $S^2_{\cF_T}(M_T)=S^1_{\cF_T}(M_T)$ and this set is bounded in $\hL^2$. Therefore, by Theorem~\ref{decomp-meas}(ii), $M_T$ is square-integrably bounded.

(iv) $\Rightarrow$ (v): Suppose that (iv) holds and let $y\in \MS(M)$. Since $y_T\in\hL^2_{\cF_T}(\O,\hR^d)$, Doob's maximal inequality implies that $\hE[\sup_{t\in[0,T]}\abs{y_t}^2]\leq 4\hE[\abs{y_T}^2]<+\infty$.
\qed 

\begin{rem}\label{unifbdd-rem}
	In the literature (cf. \cite[Section~1]{MgRT}, \cite[Section~2.3]{Kis2020}), a set-valued $\hF$-martingale $M$ is called \emph{uniformly square-integrably bounded} if $\hE[\sup_{t\in[0,T]}\norm{M_t}^2]<+\infty$. However, without any path-regularity conditions on $M$, the supremum of $(\norm{M_t})_{t\in[0,T]}$ may fail to be measurable even if $M$ is $\hF$-progressively measurable. In the cited references, this property is only needed to ensure that condition (iv) in Proposition~\ref{mtg-bdd} holds. As shown in Proposition~\ref{mtg-bdd}, it is sufficient to assume that $\hE[\norm{M_T}^2]<+\infty$ for this purpose.
	\qed
\end{rem}

Let $M=(M_t)_{t\in[0,T]}$ be a convex square-integrably bounded set-valued $\hF$-martingale. By the standard martingale representation theorem, each $y\in \MS(M)$ can be written as
$y=\cJ(x,z)$ for a unique pair $(x,z)\in\sR^{d,m}_{\hF}$. Let us define
\bea\label{RMdefn}
\cR^M:=\cb{(x,z)\in \sR^{d,m}_{\hF}\colon \cJ(x,z)\in \MS(M)}.
\eea
The next theorem provides a representation of set-valued martingales with general (non-singleton) initial values in terms of the set-valued stochastic integral.

\begin{thm}\label{MRT-initial}
	\emph{\cite[Theorem 5.6]{amw}} Let $M=(M_t)_{t\in[0,T]}$ be a convex square-integrably bounded set-valued $\hF$-martingale. Then, for each $t\in[0,T]$, it holds
	\[
	M_t = \int_{0-}^t \cR^M  dB\quad \hP\text{-a.s.}
	\]
\end{thm}

\subsection{The Martingale Property of the Stochastic Integral}

We start by establishing the submartingale property of set-valued stochastic integrals.

\begin{prop}\label{lem:submtg}
Let $\cR\subseteq \sR^{d,m}_{\hF}$ be a nonempty set. We have the following properties:

(i) $(\int_{0-}^t \cR dB)_{t\in [0,T]}$ is a set-valued square-integrable set-valued $\hF$-submartingale.

(ii) $(\int_{0-}^t \cR dB)_{t\in [0,T]}$ is square-integrably bounded if and only if $\cJ_T[\cR]$ is dominated in~$\hL^2$.

(iii) If $\cR$ is finite, then $(\int_{0-}^t \co(\cR)dB)_{t\in [0,T]}$ is a convex and square-integrably bounded set-valued $\hF$-submartingale.
\end{prop}

{\it Proof}: (i) By definition, the indefinite integral is $\hF$-adapted and square-integrable. To check the submartingale property, let $0\leq s<t\leq T$. Note that
\beaa
S^2_{\cF_s}\of{\hE\sqb{\int_{0-}^t\cR dB\mid \cF_s}}&=&\cl_{\hL^2}\cb{\hE[\xi|\cF_s]\colon \xi\in S^2_{\cF_s}\of{\int_{0-}^t \cR dB}}\\
&=&\cl_{\hL^2}\cb{\hE[\xi|\cF_s]\colon \xi\in \cl_{\hL^2}\dec_{\cF_t}(\cJ_t[\cR])}\\
&\supseteq &\cl_{\hL^2}\dec_{\cF_s}(\cb{\hE[\xi|\cF_s]\colon \xi\in \cJ_t[\cR]})\\
&=&\cl_{\hL^2}\dec_{\cF_s}(\cJ_s[\cR])=S^2_{\cF_s}\of{\int_{0-}^s \cR dB},
\eeaa
where the $\supseteq$ relation holds since $\{\hE[\xi|\cF_s]\colon \xi\in \cl_{\hL^2}\dec_{\cF_t}(\cJ_t[\cR])\}$ is $\cF_s$-decomposable and it is a superset of $\{\hE[\xi|\cF_s]\colon \xi\in \cJ_t[\cR]\}$. It follows that $\hE[\int_{0-}^t \cR dB|\cF_s]\supseteq \int_{0-}^s \cR dB$ $\hP$-a.s.

(ii) For each $t\in[0,T]$, by Theorem~\ref{decomp-meas}(ii), $\int_{0-}^t \cR dB$ is square-integrably bounded if and only if $\dec_{\cF_t}(\cJ_t[\cR])$ is a bounded set in $\hL^2$. By Lemma~\ref{kis-decomp-bdd}, this holds if and only if $\cJ_t[\cR]$ is dominated in $\hL^2$. Hence, the ``if" part of the result follows. For the ``only if" part, it is sufficient to note that $\cJ_t[\cR]$ is dominated in $\hL^2$ whenever $\cJ_T[\cR]$ is dominated in $\hL^2$ by conditional Jensen's inequality.

(iii) Suppose that $\cR=\{(x^1,z^1),\ldots,(x^n,z^n)\}$ with $n\in\hN$. Let $\xi:=\sum_{i=1}^n \abs{\cJ_T(x^i,z^i)}\in\hL^2$. Let $(x,z)\in \co(\cR)$. Hence, $(x,z)=\sum_{i=1}^n \lambda_i (x^i,z^i)$ for some $\lambda\in\Delta^{n-1}$. Then,
\[
\abs{\cJ_T(x,z)}=\abs{\sum_{i=1}^n\lambda_i\cJ_T(x^i,z^i)}\leq \sum_{i=1}^n \lambda_i\abs{\cJ_T(x^i,z^i)}\leq\sum_{i=1}^n\lambda_i \xi =\xi
\]
so that $\cJ_T[\co(\cR)]$ is dominated in $\hL^2$. By (iii), it follows that $(\int_{0-}^t \co(\cR)dB)_{t\in[0,T]}$ is square-integrably bounded. The result follows by (i) and (ii).
\qed

Our aim is to characterize all integrands for which the indefinite integral is a convex square-integrably bounded set-valued $\hF$-martingale. We begin with a collection of useful properties regarding the selections of a set-valued martingale.

\begin{prop}\label{lem:mtgsel}
	Let $M=(M_t)_{t\in[0,T]}$ be a convex square-integrably bounded set-valued $\hF$-martingale. Then, the following properties hold:
	
(i)	Let $\xi\in S^2_{\cF_T}(M_T)$ and define an $\hF$-martingale $y=(y_t)_{t\in [0,T]}$ by
\[
y_t = \hE[\xi|\cF_t],\quad t\in [0,T].
\]
Then, $y_t\in S^2_{\cF_t}(M_t)$ for each $t\in[0,T]$. In particular, $y\in \MS(M)$.

(ii) It holds $\cR^M = \{(x,z)\in \sR^{d,m}_{\hF} \colon \cJ_T(x,z)\in S^2_{\cF_T}(M_T)\}$.

(iii) For each $t\in[0,T]$, it holds $S^2_{\cF_t}(M_t)=\cJ_t[\cR^M]$; in particular, $\cJ_t[\cR^M]$ is closed and $\cF_t$-decomposable.
\end{prop}

{\it Proof}: (i) By martingale property and the definition of conditional expectation, we have
	\[
	S^2_{\cF_t}(M_t)=S^2_{\cF_t}(\hE[M_T|\cF_t])=\{\hE[\xi|\cF_t]\colon \xi \in S^2_{\cF_T}(M_T)\}, \q t\in [0,T]. 	\]
	Hence, the result follows.

(ii) As an immediate consequence of (i), we have
\beaa
\cR^M &=& \{(x,z)\in \sR^{d,m}_{\hF} \colon \cJ(x,z)\in \MS(M)\}\\ 
&=& \{(x,z)\in \sR^{d,m}_{\hF} \colon \cJ_t(x,z)\in S^2_{\cF_t}(M_t) \text{ for every }t\in[0,T]\}\\
&=& \{(x,z)\in \sR^{d,m}_{\hF} \colon \cJ_T(x,z)\in S^2_{\cF_T}(M_T)\}.
\eeaa

(iii) First, let $t=T$. The $\supseteq$ part of the equality follows directly from (ii). To prove the $\subseteq$ part, let $\xi\in S^2_{\cF_T}(M_T)$. Let $y$ be the martingale defined by (i). We have $y\in \MS(M)$. Moreover, $y=\cJ(x,z)$ for some $(x,z)\in\sR^{d,m}_{\hF}$ by the standard martingale representation theorem. Hence, by (ii), we have $(x,z)\in \cR^M$ and $\xi=y_T\in\cJ_T[\cR^M]$.

Next, let $t\in[0,T]$. As in the proof of (i), we have
\beaa
S^2_{\cF_t}(M_t) &=& \{\hE[\xi|\cF_t]\colon \xi \in S^2_{\cF_T}(M_T)\}= \{\hE[\xi|\cF_t]\colon \xi \in \cJ_T[\cR^M]\} \\
&=& \{\hE[\cJ_T(x,z)|\cF_t]\colon (x,z)\in\cR^M\}=\cJ_t[\cR^M],
\eeaa
where the second equality is by the case $t=T$ and the third is by the definition of $\cR^M$.
\qed

The next proposition complements Proposition~\ref{lem:mtgsel}(iii).

\begin{prop}\label{prop:decTt}
	Let $\cR\subseteq\sR^{d,m}_{\hF}$ be a nonempty set, let $t\in [0,T]$. If $\cJ_T[\cR]$ is $\cF_T$-decomposable, then $\cJ_t[\cR]$ is $\cF_t$-decomposable.
	\end{prop}

{\it Proof}: Suppose that $\cJ_T[\cR]$ is $\cF_T$-decomposable. Let $(x^1,z^1), (x^2,z^2)\in\cR$ and $A\in \cF_t$. Since $A\in \cF_T$ as well, we may write
$\cJ_T(x^1,z^1){\bf 1}_A + \cJ_T(x^2,z^2){\bf 1}_{A^c}=\cJ_T(x^3,z^3)$, 
for some $(x^3,z^3)\in\cR$. Then, taking conditional expectations with respect to $\cF_t$ yields
$\cJ_t(x^1,z^1){\bf 1}_A + \cJ_t(x^2,z^2){\bf 1}_{A^c} = \cJ_t(x^3,z^3)\in \cJ_t[\cR]$.
That is, $\cJ_t[\cR]$ is $\cF_t$-decomposable.
\qed

The next theorem is the main result of this subsection and it characterizes all integrands for which the set-valued stochastic integral is a set-valued martingale.

\begin{thm}\label{thm:mtgint}
	Let $\cR\subseteq \sR^{d,m}_{\hF}$ be a nonempty set such that $\cJ_T[\cR]$ is convex and dominated in $\hL^2$. Then, the following are equivalent:
	
	(i) $\cJ_T[\cR]$ is closed and $\cF^B_T$-decomposable.
	
	(ii) $\cJ_t[\cR]$ is closed and $\cF_t$-decomposable for every $t\in[0,T]$.
	
	(iii) $(\int_{0-}^t \cR dB)_{t\in [0,T]}$ is a set-valued $\hF$-martingale.
	
	In this case, $(\int_{0-}^t \cR dB)_{t\in [0,T]}$ is a convex square-integrably bounded set-valued $\hF$-martingale that admits a Hausdorff continuous modification.
	\end{thm}

{\it Proof}: The implications (iii)$\Rightarrow$(i) and (iii)$\Rightarrow$(ii) follow by Proposition~\ref{lem:mtgsel}(iii). The implication (ii)$\Rightarrow$(i) is obvious. We check (i)$\Rightarrow$(iii) next. Let $M_t:=\int_{0-}^t \cR dB$, $t\in[0,T]$. Let $t\in[0,T]$. Under (i), we have
\beaa 
S^2_{\cF_t}(\hE[M_T|\cF_t]) &=& \cl_{\hL^2} \{\hE[\xi|\cF_t]\colon \xi \in S^2_{\cF_T}(M_T)\}= \cl_{\hL^2} \{\hE[\xi|\cF_t]\colon \xi \in \cl_{\hL^2}\dec_{\cF_t}(\cJ_T[\cR])\} \\
&=&\cl_{\hL^2}\{\hE[\xi|\cF_t]\colon \xi \in \cJ_T[\cR]\}=\cl_{\hL^2}(\cJ_t[\cR])=\cl_{\hL^2}\dec_{\cF_t}(\cJ_t[\cR])= S^2_{\cF_t}(M_t),
\eeaa
where the third and fifth equalities are by Proposition~\ref{prop:decTt}. Hence, $M_t=\hE[M_T|\cF_t]$ and (iii) holds. The last sentence of the theorem is by Proposition~\ref{lem:submtg}(ii) and Theorem~\ref{cor:cont}.
\qed

\subsection{Regularity of the Stochastic Integral}\label{sec:int-reg}

According to Definition~\ref{defn:stochint}, the set-valued stochastic integral over $[0,t]$ is defined as an $\cF_t$-measurable set-valued random variable that is $\hP$-a.s. unique for each fixed $t\in [0,T]$. In particular, the martingale representation in Theorem~\ref{MRT-initial} works $\hP$-a.s. for each $t\in [0,T]$. Without further regularity, it is not clear how to make sense of the collection $(\int_{0-}^t \cR  dB)_{t\in [0,T]}$ as an indefinite integral, i.e., a measurable set-valued stochastic process.

In this subsection, using the general regularity results of Sections~\ref{sec:reg},~\ref{sec:reg2}, we will prove that set-valued stochastic integrals have regular, hence measurable, modifications in situations where the integrand can be expressed as the convex hull of a countable family.

We first prove an auxiliary result that is helpful in calculating the expected value of the maximum of finitely many random variables.

\begin{lem}\label{lem:randomization}
	Let $\cG$ be a sub-$\sigma$-algebra of $\cF$ and let $\zeta_1,\ldots,\zeta_n\in \hL^2_{\cG}(\O,\hR)$ with $n\in\hN$. Then,
	\[
	\sup_{(A_i)_{i=1}^n\in\Pi_{\cG}(\O)}\hE\sqb{\sum_{i=1}^n {\bf 1}_{A_i}\zeta_i}=\hE\sqb{\max_{i\in\{1,\ldots,n\}}\zeta_i}=	\sup_{\eta\in S^2_{\cG}(\Delta^{n-1})}\hE\sqb{\sum_{i=1}^n \eta_i\zeta_i}
	\]
	and both suprema are attained.
\end{lem}

{\it Proof}: The $\leq$ part of the first equality is trivial. To see the $\geq$ part, let us define
\[
B_i:=\cb{\zeta_i=\max_{j\in\{1,\ldots,n\}}\zeta_j}\in\cG,\quad i\in\{1,\ldots,n\}.
\]
Then, we may write $\max_{i\in\{1,\ldots,n\}}\zeta_i=\sum_{i=1}^n {\bf 1}_{A_i}\zeta_i$, where $(A_i)_{i=1}^n\in\Pi_{\cG}(\O)$ is defined by
\bea\label{eq:part}
A_1:=B_1\in\cG,\quad A_2:=B_2\setminus A_1 \in\cG,\quad \ldots, \quad A_n:=B_n\setminus(A_1\cup\ldots\cup A_{n-1})\in\cG.
\eea
Hence, the $\geq$ part of the first equality of the lemma follows. To see the second equality, we apply \cite[Theorem~2.2]{hiai-umegaki} to get
\[
\sup_{\eta\in S^2_{\cG}(\Delta^{n-1})}\hE\sqb{\sum_{i=1}^n \eta_i\zeta_i}=\hE\sqb{\sup_{\alpha\in\Delta^{n-1}}\sum_{i=1}^n \alpha_i\zeta_i}=\hE\sqb{\max_{i\in\{1,\ldots,n\}}\zeta_i},
\]
where the second equality holds since the maximum value of the convex combination of finitely many real numbers is the maximum of these numbers. This completes the proof of the equalities. As a byproduct, we have shown that the first supremum is attained at $(A_i)_{i=1}^n$ defined in \eqref{eq:part} and the second supremum is attained at the corresponding vector $\eta=({\bf 1}_{A_1},\ldots,{\bf 1}_{A_n})\in S^2_{\cG}(\Delta^{n-1})$.
\qed

Next, we prove the continuity of the set-valued stochastic integral when the integrand is generated by a finite set.

\begin{lem}\label{lem:reg-int}
Let $\cR\subseteq \sR^{d,m}_{\hF}$ be a nonempty finite set. Then, $(\int_{0-}^t \co(\cR)dB)_{t\in[0,T]}$ has a Hausdorff continuous modification.
\end{lem}

{\it Proof}: Let $\cR:=\{(x^i,z^i)\colon i\in\{1,\ldots,n\}\}$ for some $n\in\hN$. For each $x^\ast\in\hB_{\hR^d}$ and $t\in [0,T]$, let us define
\[
g_{x^\ast}(t):=s\of{x^\ast,\hE\sqb{\int_{0-}^t \co(\cR) dB}}.
\]
By the submartingale property of stochastic integrals (Proposition~\ref{lem:submtg}(i)), $g_{x^\ast}$ is increasing.

Let $x^\ast\in\hB_{\hR^d}$ and $t\in[0,T]$. We have
\bea
g_{x^\ast}(t)&=&\hE\sqb{s\of{x^\ast,\int_{0-}^t \co(\cR)dB}}
\notag = \hE\sqb{s\of{x^\ast,\cl_{\hR^d}\co\of{\int_{0-}^t \cR dB}}}\notag\\
& =& \hE\sqb{s\of{x^\ast,\int_{0-}^t \cR dB}}\notag =\hE\sqb{s\of{x^\ast,\{\cJ_t(x^i,z^i)\colon i\in\{1,\ldots,n\}\}}}\notag \\
&=&\hE\sqb{\max_{i\in\{1,\ldots,n\}}\of{\ip{x^\ast,x^i}+\int_0^t \ip{x^\ast,z^i_s}dB_s}},\label{eq:fincase}
\eea
where the first equality is by \cite[Theorem~2.1.35]{Molchanov}, the second is by Proposition~\ref{lem:basicprop}(ii), the third is by the fact that the support function of a set and that of its closed convex hull are identical, and the fourth is by Proposition~\ref{lem:basicprop}(iv).

By Proposition~\ref{lem:submtg}, $(\int_{0-}^t \co(\cR) dB)_{t\in[0,T]}$ is a convex square-integrably bounded set-valued $\hF$-submartingale. Let $x^\ast\in\sW$. We claim that $g_{x^\ast}$ is continuous so that $(\int_{0-}^t \co(\cR) dB)_{t\in[0,T]}$ has a Hausdorff continuous modification by Corollary~\ref{thm:cont}.

Let $t\geq 0$, $\delta>0$ be such that $t+\delta\leq T$. Since the function $g_{x^\ast}$ is increasing, we have
\beaa
&&\left| g_{x^\ast}(t+\delta)-g_{x^\ast}(t)\right|
=  g_{x^\ast}(t+\delta)-g_{x^\ast}(t)\\
&&=\hE\sqb{\max_{i\in\{1,\ldots,n\}}\of{\ip{x^\ast,x^i}+\int_0^{t+\delta} \ip{x^\ast,z^i_s}dB_s}} -\hE\sqb{\max_{i\in\{1,\ldots,n\}}\of{\ip{x^\ast,x^i}+\int_0^t \ip{x^\ast,z^i_s}dB_s}}\\
&&\leq \hE\sqb{\max_{i\in\{1,\ldots,n\}}\abs{\int_t^{t+\delta}\ip{x^i,z_s^i}dB_s}}.
\eeaa 
Let $(A_i)_{i=1}^n\in\Pi_{\cF_{t+\delta}}(\O)$ be a partition such that 
\[
\hE\sqb{\max_{i\in\{1,\ldots,n\}}\abs{\int_t^{t+\delta}\ip{x^i,z_s^i}dB_s}}=\hE\sqb{\sum_{i=1}^n {\bf 1}_{A_i}\abs{\int_t^{t+\delta}\ip{x^\ast,z_s^i}dB_s}},
\]
whose existence is guaranteed by Lemma~\ref{lem:randomization}. Then, by Cauchy-Schwarz-Bunyakovski inequality and It\^{o} isometry, we have
\[
\hE\sqb{\sum_{i=1}^{n}{\bf 1}_{A_i}\abs{\int_t^{t+\delta}\ip{x^\ast,z_s^i}dB_s}}\leq \sum_{i=1}^{n}\of{\hP(A_i)\hE\sqb{\int_t^{t+\delta}|\ip{x^\ast,z_s^i}|^2ds}}^{\frac12}.
\]
Letting $\delta\rightarrow 0$, the sum on the right of the above inequality converges to zero. This establishes the right-continuity of $g^n_{x^\ast}$. By replacing $t+\delta$ with $t-\delta$ and following similar arguments, its left-continuity follows as well.
\qed

The next theorem is the main result on the regularity of set-valued stochastic integrals.

\begin{thm}\label{thm:reg-int}
	Let $\cR\subseteq \sR^{d,m}_{\hF}$ be a nonempty countable set.
	
	(i) Suppose that $\cJ_T[\cR]$ is dominated in $\hL^2$. Then, $(\int_{0-}^t \co(\cR)dB)_{t\in[0,T]}$ has a Hausdorff c\`{a}gl\`{a}d modification.
	
	(ii) Suppose that the following summability condition holds:
	\be\label{eq:summable}
	\hE\sqb{\sum_{i=1}^{\infty}\abs{\cJ_T(x^i,z^i)}^2}=\sum_{i=1}^\infty \of{\abs{x^i}^2+\hE\sqb{\int_0^T\abs{z_s^i}^2ds}}<+\infty.
	\ee
	Then, $(\int_{0-}^t \co(\cR)dB)_{t\in[0,T]}$ has a Hausdorff continuous modification.
\end{thm}

{\it Proof}: Let us write $\cR=\{(x^i,z^i)\colon i\in\hN\}$. For convenience, we define $\cR^n:=\{(x^i,z^i)\colon i\in\{1,\ldots,n\}\}$. For each $x^\ast\in\hB_{\hR^d}$ and $t\in[0,T]$, we also define
\[
g^n_{x^\ast}(t):=s\of{x^\ast,\hE\sqb{\int_{0-}^t \co(\cR^n) dB}},\quad n\in\hN;\quad g_{x^\ast}(t):=s\of{x^\ast,\hE\sqb{\int_{0-}^t \co(\cR) dB}}.
\]

(i) By the submartingale property of stochastic integrals (Proposition~\ref{lem:submtg}(i)), $g_{x^\ast}^n$, $n\in\hN$, and $g_{x^\ast}$ are increasing functions. By the proof of Lemma~\ref{lem:reg-int}, $g^n_{x^\ast}$ is also continuous for each $n\in\hN$. Let $t\in[0,T]$. By \eqref{eq:fincase} in the same proof, we have
\[
g^n_{x^\ast}(t)=\hE\sqb{\max_{i\in\{1,\ldots,n\}}\of{\ip{x^\ast,x^i}+\int_0^t \ip{x^\ast,z^i_s}dB_s}}.
\]
Similarly, we also have
\be\label{eq:infcase}
g_{x^\ast}(t)=\hE\sqb{\sup_{i\in\hN}\of{\ip{x^\ast,x^i}+\int_0^t \ip{x^\ast,z^i_s}dB_s}}.
\ee
Moreover,
\beaa
\lim_{n\rightarrow\infty}g^n_{x^\ast}(t)&=&\lim_{n\rightarrow\infty}\hE\sqb{\max_{i\in\{1,\ldots,n\}}\of{\ip{x^\ast,x^i}+\int_0^t \ip{x^\ast,z^i_s}dB_s}}\\
&=&\hE\sqb{\sup_{i\in\hN}\of{\ip{x^\ast,x^i}+\int_0^t \ip{x^\ast,z^i_s}dB_s}}= g_{x^\ast}(t),
\eeaa 
where the first equality is by \eqref{eq:fincase}, the second is by monotone convergence theorem, and the last is by \eqref{eq:infcase}; the limit here is indeed a supremum since $(g^n_{x^\ast}(t))_{n\in\hN}$ is an increasing sequence. Then, for each $r\in\hR$, we have
\[
\cb{t\in[0,T]\colon g_{x^\ast}(t)\leq r} =\bigcap_{n\in\hN}\cb{t\in[0,T]\colon g^n_{x^\ast}(t)\leq r}.
\] 
Since $g^n_{x^\ast}$ is a continuous function for each $n\in\hN$ as shown in the proof of (i), the above intersection is a closed set. It follows that $g_{x^\ast}$ is an increasing lower semicontinuous function; hence, it is a left-continuous function.

By Proposition~\ref{lem:submtg}, $(\int_{0-}^t \co(\cR) dB)_{t\in[0,T]}$ is a convex square-integrably bounded set-valued $\hF$-submartingale. Therefore, by Corollary~\ref{thm:caglad}, $(\int_{0-}^t \co(\cR) dB)_{t\in[0,T]}$ has a Hausdorff c\`{a}gl\`{a}d modification.

(ii) Let $\hC([0,T])$ be the set of all continuous real-valued functions on $[0,T]$ equipped with the supremum norm. Note that $g^n_{x^\ast}\in\hC([0,T])$ for each $n\in\hN$ and $(g^n_{x^\ast})_{n\in\hN}$ converges to $g$ pointwise. We show that $(g^n_{x^\ast})_{n\in\hN}$ is a Cauchy sequence in $\hC([0,T])$ so that the convergence to $g$ is also uniform. Let $m,n\in\hN$ with $m>n$. Then, we have
\beaa
\sup_{t\in[0,T]}\abs{g^m_{x^\ast}(t)-g^n_{x^\ast}(t)}
&=&\sup_{t\in[0,T]}\of{g^m_{x^\ast}(t)-g^n_{x^\ast}(t)}\\
&=&\sup_{t\in[0,T]}\of{\hE\sqb{\max_{i\in\{1,\ldots,m\}}\ip{x^\ast,\cJ_t(x^i,z^i)}} - \hE\sqb{\max_{i\in\{1,\ldots,n\}}\ip{x^\ast,\cJ_t(x^i,z^i)}}}.
\eeaa 
By Lemma~\ref{lem:randomization}, there exists $(A^t_i)_{i=1}^m\in\Pi_{\cF_t}(\O)$ such that
\[
\hE\sqb{\max_{i\in\{1,\ldots,m\}}\ip{x^\ast,\cJ_t(x^i,z^i)}} = \hE\sqb{\sum_{i=1}^m {\bf 1}_{A_i^t}\ip{x^\ast,\cJ_t(x^i,z^i)}};
\]
and we also have
\[
\hE\sqb{\max_{i\in\{1,\ldots,n\}}\ip{x^\ast,\cJ_t(x^i,z^i)}}=\sup_{\eta\in S^2_{\cF_t}(\Delta^{n-1})}\hE\sqb{\sum_{i=1}^n \eta_i\ip{x^\ast,\cJ_t(x^i,z^i)}}.
\]
Let us define $\bar{\eta}^t\in S^2_{\cF_t}(\Delta^{n-1})$ by
$\bar{\eta}^t_i:=\frac{{\bf 1}_{A^t_i}}{{\bf 1}_{A^t_1}+\ldots+{\bf 1}_{A^t_n}}$, $i\in\{1,\ldots,n\}$.
Then, we have
\beaa
& &\sup_{t\in[0,T]}\abs{g^m_{x^\ast}(t)-g^n_{x^\ast}(t)}\\
&&=\sup_{t\in[0,T]}\of{\hE\sqb{\sum_{i=1}^m {\bf 1}_{A_i^t}\ip{x^\ast,\cJ_t(x^i,z^i)}}-\sup_{\eta\in S^2_{\cF_t}(\Delta^{n-1})}\hE\sqb{\sum_{i=1}^n \eta_i\ip{x^\ast,\cJ_t(x^i,z^i)}}}\\
&&\leq \sup_{t\in[0,T]}\hE\sqb{\sum_{i=1}^m {\bf 1}_{A_i^t}\ip{x^\ast,\cJ_t(x^i,z^i)}-\sum_{i=1}^n \bar{\eta}^t_i\ip{x^\ast,\cJ_t(x^i,z^i)}}\\
&&=\sup_{t\in[0,T]}\hE\sqb{\sum_{i=1}^n ({\bf 1}_{A_i^t}-\bar{\eta}^t_i)\ip{x^\ast,\cJ_t(x^i,z^i)}+\sum_{i=n+1}^m {\bf 1}_{A_i^t}\ip{x^\ast,\cJ_t(x^i,z^i)}}.
\eeaa 
Note that, for each $i\in\{1,\ldots,n\}$, we have
\[
{\bf 1}_{A_i^t}-\bar{\eta}^t_i={\bf 1}_{A^t_i}\of{1-\frac{{\bf 1}_{A^t_1}+\ldots+{\bf 1}_{A^t_m}}{{\bf 1}_{A^t_1}+\ldots+{\bf 1}_{A^t_n}}}=-{\bf 1}_{A^t_i}\frac{{\bf 1}_{A^t_{n+1}}+\ldots+{\bf 1}_{A^t_m}}{{\bf 1}_{A^t_1}+\ldots+{\bf 1}_{A^t_n}}=0
\]
since $A^t_i\cap A^t_j=0$ for every $j\in\{n+1,\ldots,m\}$. Therefore,
\beaa 
\sup_{t\in[0,T]}\abs{g^m_{x^\ast}(t)-g^n_{x^\ast}(t)}&\leq& \sup_{t\in[0,T]}\hE\sqb{\sum_{i=n+1}^m {\bf 1}_{A_i^t}\ip{x^\ast,\cJ_t(x^i,z^i)}}\\
&\leq& \sup_{t\in[0,T]}\hE\sqb{\sum_{i=n+1}^m \abs{\cJ_t(x^i,z^i)}}\leq \sup_{t\in[0,T]}\of{\hE\sqb{\sum_{i=n+1}^m \abs{\cJ_t(x^i,z^i)}^2}}^{\frac12}.
\eeaa 
Note that $\hE[|\cJ_t(x^i,z^i)|^2]=|x^i|^2+\hE[\int_0^t |z_s^i|^2ds]$ for each $i\in\{n+1,\ldots,m\}$. Then, we get
\[
\sup_{t\in[0,T]}\abs{g^m_{x^\ast}(t)-g^n_{x^\ast}(t)}\leq \of{\sum_{i=n+1}^m |x^i|^2+\hE\sqb{\int_0^T |z_s^i|^2ds}}^{\frac12}=\of{\hE\sqb{\sum_{i=n+1}^{m}\abs{\cJ_T(x^i,z^i)}^2}}^{\frac12}.
\]
The last upper bound converges to $0$ as $n,m\rightarrow\infty$ thanks to the summability condition \eqref{eq:summable}. Hence, $(g^n_{x^\ast})_{n\in\hN}$ is a Cauchy sequence in $\hC([0,T])$ and its limit pointwise limit $g_{x^\ast}$ is also its uniform limit, thus a continuous function. This verifies the hypothesis of Corollary~\ref{thm:cont} and the existence of a Hausdorff continuous modification follows.
\qed 

\section{Set-Valued Stochastic Integral vs.~Conditional Expectation}\label{sec:compare}

Let $\Xi\in\scL^2_{\cF_T}(\O,\sK(\hR^d))$ be a given set-valued random variable. It induces a convex and square-integrably bounded set-valued $\hF$-martingale $M=(M_t)_{t\in[0,T]}$ by
$M_t:=\hE[\Xi|\cF_t]$, $t\in[0,T]$.
Thanks to the set-valued martingale representation theorem (Theorem~\ref{MRT-initial}) and the path continuity result (Corollary~\ref{cor:cont}), there exists an integrand $\cR^M\subseteq\sR^{d,m}_{\hF}$ such that
$M_t=\int_{0-}^t \cR^M dB$, $t\in[0,T]$,
with probability one. However, the set $\cR^M$ is typically uncountable and it is desirable to rewrite $\int_{0-}^t \cR^M dB$ in terms of the stochastic integral of a countable family in $\sR^{d,m}_{\hF}$. It is easy to see that
$S^2_{\cF_T}(\Xi)=P_T[\MS(M)]$ and $\MS(M)=\{(\hE[\xi|\cF_t])_{t\in[0,T]}\colon \xi\in S^2_{\cF_T}(\Xi)\}$.

On the other hand, thanks to the convexity of $\Xi$, one can find a sequence $(\xi^i)_{i\in\hN}$ in $S^2_{\cF_T}(\Xi)$ that gives a convex Castaing representation of $\Xi$, i.e.,
$\Xi(\o)=\cl_{\hR^d}\co(\{\xi^i(\o)\colon i\in\hN\})$ for
$\hP$-a.e. $\o\in\O$ (cf. \cite[Theorem~2.3.1]{Kis2020}). Then, by \cite[Corollary~2.3.2, Theorem~2.3.5]{Kis2020},
\bea\label{eq:selMT}
S^2_{\cF_T}(\Xi)=\cl_{\hL^2}\co\dec_{\cF_T}(\{\xi^i\colon i\in\hN\}).
\eea

Moreover, for each $i\in\hN$, there exists $(x^i,z^i)\in\cR^M$ such that $\cJ_t(x^i,z^i)=\hE[\xi^i|\cF_t]$, $t\in [0,T]$, with probability one. Let
$\cR^C:=\{(x^i,z^i)\colon i\in\hN\}$. 
In particular, $\{\xi^i\colon i\in\hN\}=\cJ_T[\cR^C]$. By Proposition~\ref{lem:basicprop}(ii), we have 
\[
\int_{0-}^t \co(\cR^C)dB=\cl_{\hR^d}\co\of{\int_{0-}^t \cR^C dB}\quad \hP\text{-a.s.}, ~t\in[0,T].
\]

\begin{rem}\label{rem:reg}
	Since $(\int_{0-}^t \cR^M dB)_{t\in[0,T]}$ is a convex and square-integrably bounded set-valued $\hF$-martingale, by Theorem~\ref{cor:cont}, it has a Hausdorff continuous modification. Moreover, by Theorem~\ref{thm:reg-int}, $(\int_{0-}^t \co(\cR^{C})dB)_{t\in[0,T]}$ has a Hausdorff c\`{a}gl\`{a}d modification. \qed
	\end{rem}

We will first show that the set-valued stochastic integrals of $\cR^M$ and $\cR^C$ coincide at time $T$. As a preparation, the next lemma shows that the closed convex hull of the decomposable hull of a set in $\hL^2$ consists of random convex combinations of its elements. 

\begin{lem}\label{lem:chcd}
	Let $\cG$ be a sub-$\sigma$-algebra of $\cF$ and let $\cK\subseteq\hL^2_{\cG}(\hR^d)$ be a nonempty set. Then,
	\beaa
	\cl_{\hL^2}\co\dec_{\cG}(\cK)&=&\cl_{\hL^2}\dec_{\cG}\co(\cK)\\
	&=&\cl_{\hL^2}\cb{\sum_{i=1}^n \lambda_i \zeta^i\colon \lambda \in S^2_{\cG}(\Delta^{n-1});\ \zeta^1,\ldots,\zeta^n\in\cK;\ n\in\hN}.
	\eeaa
	\end{lem}
	
{\it Proof}: The first identity is given by \cite[Corollary 7.6]{umur}. Let $\cA$ denote the set of all random convex combinations of $\cK$, i.e., the last set in the lemma. By \cite[Corollary 7.5]{umur}, the closed convex $\cG$-decomposable hull of a subset of $\hL^2_{\cG}(\hR^d)$ contains all random convex combinations of its elements. This implies that $\cA\subseteq \cl_{\hL^2}\dec_{\cG}\co(\cK)$. To see the converse, let us observe that $\cA$ is a convex, $\cG$-decomposable set and we have $\cK\subseteq\cA$. Hence, $\dec_{\cG}\co(\cK)\subseteq \cA$ so that $\cl_{\hL^2}\dec_{\cG}\co(\cK)\subseteq\cl_{\hL^2}(\cA)$. With this, the second identity in the lemma follows.
\qed 

\begin{prop}\label{prop:twointegralsT}
	We have $\int_{0-}^T \cR^M dB= \int_{0-}^T \co(\cR^C)dB = \Xi$ $\hP$-a.s. Moreover,
	\bea\label{eq:twointegralsT}
	S^2_{\cF_T}(\Xi)=\cl_{\hL^2}\cb{\sum_{i=1}^n \lambda_i \cJ_T(x^i,z^i)\colon \lambda\in S^2_{\cF_T}(\Delta^{n-1}),\ n\in\hN}.
	\eea
	\end{prop}

{\it Proof}: By \eqref{eq:selMT} and Lemma~\ref{lem:chcd}, \eqref{eq:twointegralsT} follows immediately. On the other hand,
\beaa
S^2_{\cF_T}\of{\int_{0-}^T\co(\cR^C) dB}&=&\cl_{\hL^2}\co\of{S^2_{\cF_T}\of{\int_{0-}^T \cR^C dB}}\\
&=& \cl_{\hL^2}\co \dec_{\cF_T} \of{\{\cJ_T(x^i,z^i)\colon i\in\hN\}}\\
&=&\cl_{\hL^2}\cb{\sum_{i=1}^n \lambda_i \cJ_T(x^i,z^i)\colon \lambda\in S^2_{\cF_T}(\Delta^{n-1}),\ n\in\hN},
\eeaa
where the last equality follows by Lemma~\ref{lem:chcd}. Since $\cJ_T(x^i,z^i)=\xi^i$ $\hP$-a.s. for each $i\in\hN$, the integrals $\int_{0-}^T \cR^M dB$, $\int_{0-}^T \co(\cR^C)dB $ have the same set of square-integrable selections; hence, they are equal $\hP$-a.s.
\qed 

The next lemma is helpful in some calculations to follow. It extends the following simple fact to the set-valued case: for integrable real-valued random variables $\alpha,\beta$ satisfying $\alpha\leq \beta$ $
\hP$-a.s., we have $\alpha=\beta$ $\hP$-a.s. if and only if $\hE[\alpha]=\hE[\beta]$.

\begin{lem}\label{lem:nonnegative}
	Let $\Xi^1,\Xi^2\in\scL^1_{\cF_T}(\O,\sK(\hR^d))$ with $\Xi^1\subseteq \Xi^2$ $\hP$-a.s. Then, $\Xi^1=\Xi^2$ $\hP$-a.s. if and only if $\hE[\Xi^1]=\hE[\Xi^2]$.
	\end{lem}
	
{\it Proof}: Clearly, $\Xi^1=\Xi^2$ $\hP$-a.s. implies $\hE[\Xi^1]=\hE[\Xi^2]$. Conversely, suppose that $\hE[\Xi^1]=\hE[\Xi^2]$. Let $x^\ast\in\sW$. Since $\Xi^1\subseteq \Xi^2$ $\hP$-a.s., we have $s(x^\ast,\Xi^1)\leq s(x^\ast,\Xi^2)$ $\hP$-a.s. Since $\hE[\Xi^1]=\hE[\Xi^2]$, by \cite[Theorem~2.1.35]{Molchanov}, we have
$\hE[s(x^\ast,\Xi^1)]=s(x^\ast,\hE[\Xi^1])=s(x^\ast,\hE[\Xi^2])=\hE[s(x^\ast,\Xi^2)]$.
Therefore, $s(x^\ast,\Xi^1)= s(x^\ast,\Xi^2)$ $\hP$-a.s. Since $x^\ast\in\sW$ is arbitrary, we get $\Xi^1=\Xi^2$ $\hP$-a.s.
\qed 

We collect some easy consequences of Proposition~\ref{prop:twointegralsT} in the next corollary.

\begin{cor}\label{cor:easy}
	We have the following results:
	
	(i) For every $t\in[0,T]$,
	\[
	\hE\sqb{\int_{0-}^t \cR^M dB}=M_0=\hE[\Xi]=\cl_{\hR^d}\cb{\sum_{i=1}^n \hE[\lambda_i \xi^i]\colon \lambda\in S^2_{\cF_T}(\Delta^{n-1}),\ n\in\hN}.
	\]
	
	(ii) It holds $\int_{0-}^0 \co(\cR^C) dB=\cl_{\hR^d}\co(\{x^i\colon i\in\hN\})$ $\hP$-a.s.
\end{cor}

{\it Proof}: (i) This is an immediate consequence of \eqref{eq:twointegralsT} in Proposition~\ref{prop:twointegralsT}.

(ii) By the definition of the set-valued stochastic integral, we have
\[
S^2_{\cF_0}\of{\int_{0-}^0 \co(\cR^C)dB} = \cl_{\hL^2}\cb{\sum_{i=1}^n \eta_i x^i\colon \eta\in S^2_{\cF_0}(\Delta^{n-1}),n\in\hN}=\cl_{\hR^d}\co(\{x^i\colon i\in\hN\})
\]
since $\cF_0$-measurable random variables are deterministic $\hP$-a.s.
\qed

Due to Proposition~\ref{prop:twointegralsT}, it is very tempting to expect that the integrals $\int_{0-}^t \cR^M dB$ and $\int_{0-}^t \co(\cR^C)dB$ are equal at each $t<T$ as well, thus yielding a representation of the set-valued martingale $M$ in terms of a countable family of integrands. However, we will show that the former integral is a superset of the latter in general (Theorem~\ref{prop:twointegrals}) and we will show that the equality works only in some exceptional cases (Theorem~\ref{prop:twointegrals2}, Theorem~\ref{cor:twointegrals2}).

\begin{thm}\label{prop:twointegrals}
	Let $t\in[0,T]$. Then, with probability one, we have
	\beaa
	 \int_{0-}^t \cR^C dB &=& \cl_{\hR^d}\of{\cb{\hE[\xi^i|\cF_t]\colon i\in\hN}},\\
	 \int_{0-}^t \cR^M dB &=& \bigcap_{x^\ast\in\hB_{\hR^d}}\cb{x\in\hR^d\colon \ip{x^\ast,x}
		\leq \hE\sqb{\sup_{i\in\hN}\ip{x^\ast,\xi^i}\mid \cF_t}};
	\eeaa
	in particular, we have
	\bea\label{eq:inpart}
	\int_{0-}^t \cR^M dB = \hE[\Xi| \cF_t]\supseteq \cl_{\hR^d}\co\of{\cb{\hE[\xi^i| \cF_t]\colon i\in\hN}}=\int_{0-}^t \co(\cR^C)dB.
	\eea
\end{thm}

{\it Proof}: The first equality is a direct consequence of Proposition~\ref{lem:basicprop}(iv) since $\cJ_t(x^i,z^i)=\hE[\xi^i|\cF_t]$ $\hP$-a.s. for each $i\in\hN$ and $t\in[0,T]$. %, and $(\int_{0-}^t \co(\cR^C)dB)_{t\in[0,T]}$ has Hausdorff c\`{a}dl\`{a}g paths.

Let $t\in [0,T]$. Using \eqref{eq:twointegralsT} in Proposition~\ref{prop:twointegralsT}, we obtain
\beaa
S^2_{\cF_t}\of{\int_{0-}^t \cR^M dB}&=& S^2_{\cF_t}(\hE[\Xi|\cF_t])= \cl_{\hL^2}\{\hE[\xi|\cF_t]\colon \xi\in S^2_{\cF_T}(\Xi)\}\\
&=& \cl_{\hL^2}\cb{\hE\sqb{\sum_{i=1}^n \lambda_i \xi^i |\cF_t}\colon \lambda\in S^2_{\cF_T}(\Delta^{n-1}),\ n\in\hN}.
\eeaa
Let $x^\ast\in\sW$. Then, we have
\beaa
\hE\sqb{s\of{x^\ast,\int_{0-}^t \cR^M dB}{\bf 1}_A}\negthinspace&=&\negthinspace\hE\sqb{\sup\cb{\ip{x^\ast,x}\colon x\in \int_{0-}^t \cR^M dB}{\bf 1}_A}\\
&=&\negthinspace \sup\cb{\hE[\ip{x^\ast,\eta}{\bf 1}_A]\colon \eta\in S^2_{\cF_t}\of{\int_{0-}^t \cR^M dB}}\\
&=&\negthinspace \sup\cb{\hE[\ip{x^\ast,\eta}{\bf 1}_A]\colon \eta\in S^2_{\cF_t}\of{\hE[\Xi|\cF_t]}}\\
&=&\negthinspace \sup\cb{\hE[\ip{x^\ast,\hE[\xi|\cF_t]}{\bf 1}_A]\colon \xi\in S^2_{\cF_T}(\Xi) }\\
&=&\negthinspace \sup\cb{\hE[\ip{x^\ast,\hE[\xi|\cF_t]}{\bf 1}_A]\colon \xi\in \cl_{\hL^2}\co\dec_{\cF_T}(\{\xi^i\colon i\in\hN\}) },
\eeaa
where the second equality is by \cite[Theorem~2.2]{hiai-umegaki}, the fourth is by the properties of conditional expectation, and the fifth is by \eqref{eq:selMT}. Then, by Lemma~\ref{lem:chcd} and the continuity of expectation,
\beaa
\hE\sqb{s\of{x^\ast,\int_{0-}^t \cR^M dB}{\bf 1}_A} \negthinspace\negthinspace\negthinspace\negthinspace
&=&\negthinspace\negthinspace\negthinspace\negthinspace \sup\cb{\hE\sqb{\ip{ x^\ast,\hE\sqb{\sum_{i=1}^n \lambda_i \xi^i | \cF_t}}{\bf 1}_A}\colon \lambda\in S^2_{\cF_T}(\Delta^{n-1}), n\in\hN}\\
&=&\negthinspace\negthinspace \sup\cb{\hE\sqb{\ip{x^\ast,\sum_{i=1}^n \lambda_i \xi^i }{\bf 1}_A}\colon \lambda\in S^2_{\cF_T}(\Delta^{n-1}),\ n\in\hN}\\
&=&\negthinspace\negthinspace\negthinspace\negthinspace\sup_{n\in\hN}\sup\cb{\hE\sqb{\ip{x^\ast,\sum_{i=1}^n \lambda_i \xi^i }{\bf 1}_A}\colon \lambda\in S^2_{\cF_T}(\Delta^{n-1})}.
\eeaa
To evaluate the inner supremum in the last line, we use \cite[Theorem~2.2]{hiai-umegaki} to interchange the expectation and the supremum. Then, we obtain
\beaa 
\hE\sqb{s\of{x^\ast,\int_{0-}^t \cR^M dB}{\bf 1}_A}
&=&\sup_{n\in\hN}\hE\sqb{\sup_{\alpha\in \Delta^{n-1}}\ip{x^\ast,\sum_{i=1}^n \alpha_i \xi^i}{\bf 1}_A}\\
&=&\sup_{n\in\hN}\hE\sqb{\sup_{\alpha\in \Delta^{n-1}}\sum_{i=1}^n \alpha_i \ip{x^\ast, \xi^i}{\bf 1}_A}\\
&=&\sup_{n\in\hN}\hE\sqb{\max\{\ip{x^\ast, \xi^1},\ldots,\ip{x^\ast,\xi^n}\}{\bf 1}_A}\\
&=&\hE\sqb{\sup_{i\in\hN}\ip{x^\ast,\xi^i}{\bf 1}_A}=\hE\sqb{\hE\sqb{\sup_{i\in\hN}\ip{x^\ast,\xi^i}|\cF_t}{\bf 1}_A},
\eeaa
where the third equality is by the simple observation that the maximum value of the convex combination of finitely many real numbers is their maximum, the fourth is by monotone convergence theorem. Since $A\in\cF_t$ is arbitrary, we conclude that $s(x^\ast,\int_{0-}^t \cR^M dB)=\hE[\sup_{i\in\hN}\ip{x^\ast,\xi^i}|\cF_t]$ $\hP$-a.s. Since $\sW$ is countable, we have this equality for every $x^\ast\in\sW$ with probability one; let $\O_0\in\cF_t$ be the $\hP$-a.s. set here. To extend this property to the whole set $\hB_{\hR^d}$, let $x^\ast\in \hB_{\hR^d}$. Then, there exists a sequence $(x^{\ast,n})_{n\in\hN}$ in $\sW$ that converges to $x^\ast$. Note that $\sup_{i\in\hN}\ip{x^{\ast,n},\xi^i}\leq \norm{\Xi}$ $\hP$-a.s. for each $n\in\hN$ and we have $\norm{\Xi}\in\hL^2_{\cF_T}(\O,\hR)$. Hence, by the conditional version of dominated convergence theorem and the continuity of $s(\cdot,\int_{0-}^t \cR^M dB)$, we obtain
\beaa
s\of{x^\ast,\int_{0-}^t \cR^M dB}&=&\lim_{n\rightarrow\infty}s\of{x^{\ast,n},\int_{0-}^t \cR^M dB}\\
&=&\lim_{n\rightarrow\infty} \hE\sqb{\sup_{i\in\hN}\ip{x^{\ast,n},\xi^i}|\cF_t}=\hE\sqb{\sup_{i\in\hN}\ip{x^\ast,\xi^i}|\cF_t}
\eeaa
holds on $\O_0$. Therefore, $s(x^\ast,\int_{0-}^t \cR^M dB)=\hE[\sup_{i\in\hN}\ip{x^\ast,\xi^i}|\cF_t]$ for every $x^\ast\in \hB_{\hR^d}$ on $\O_0$. For every closed convex set $C\subseteq\hR^d$, we may write
$C=\bigcap_{x^\ast\in \hB_{\hR^d}}\{x\in\hR^d\colon \ip{x^\ast,x}\leq s(x^\ast,C)\}$.
Applying this for $C=\int_{0-}^t \cR^M dB$ yields the second equality in the theorem.

Finally, we show \eqref{eq:inpart}. The first equality is by the set-valued martingale representation theorem, the last equality is by Proposition~\ref{lem:basicprop}(ii). Since $\Xi=\cl_{\hR^d}\co(\{\xi^i\colon i\in\hN\})$ $\hP$-a.s., the superset relation in the theorem follows immediately.
\qed 

We now characterize the cases where the stochastic integrals of $\cR^M$ and $\co(\cR^C)$ coincide.

\begin{thm}\label{prop:twointegrals2}
	The following are equivalent:
	
	(i) It holds $\int_{0-}^t \cR^M dB = \int_{0-}^t \co(\cR^C)dB$ for every $t\in[0,T]$ $\hP$-a.s.
	
	(ii) $(\int_{0-}^t \co(\cR^C)dB)_{t\in[0,T]}$ is an $\hF$-martingale.
	
	(iii) It holds $\hE[\int_{0-}^t \co(\cR^C)dB]=\cl_{\hR^d}\co(\{x^i\colon i\in\hN\})$ for every $t\in[0,T]$.
	
	(iv) It holds $\hE[\int_{0-}^T \co(\cR^C)dB]=\cl_{\hR^d}\co(\{x^i\colon i\in\hN\})$.
	
	(v) For every $x^\ast\in\hR^d$ and $t\in[0,T]$, it holds $\hE[\sup_{i\in\hN}\ip{x^\ast,\hE[\xi^i|\cF_t]}]=\sup_{i\in\hN}\ip{x^\ast,x^i}$.
	
	(vi) For every $x^\ast\in\hR^d$, it holds $\hE[\sup_{i\in\hN}\ip{x^\ast,\xi^i}]=\sup_{i\in\hN}\ip{x^\ast,x^i}$.
\end{thm}

{\it Proof}:  The implications (i)$\Rightarrow$(ii), (iii)$\Rightarrow$(iv), (v)$\Rightarrow$(vi) are obvious, (ii)$\Rightarrow$(iii) is an immediate consequence of Corollary~\ref{cor:easy}(ii).

(ii)$\Rightarrow$(i): Suppose that $(\int_{0-}^t \co(\cR^C)dB)_{t\in[0,T]}$ is an $\hF$-martingale. Then, by Proposition~\ref{prop:twointegrals}(ii),
\[
\int_{0-}^t \cR^M dB=\hE\sqb{\int_{0-}^T \cR^M dB\mid \cF_t}=\hE\sqb{\int_{0-}^T \co(\cR^C) dB\mid \cF_t}=\int_{0-}^t \co(\cR^C)dB\quad \hP\text{-a.s.}
\]
By the path-regularity of these set-valued processes, (i) follows.

(iv)$\Rightarrow$(ii): Suppose that $\hE[\int_{0-}^T \co(\cR^C)dB]=\cl_{\hR^d}\co(\{x^i\colon i\in\hN\})$, i.e., 
\[
\hE\sqb{\int_{0-}^T \co(\cR^C)dB}=\int_{0-}^0 \co(\cR^C)dB.
\]
Since $(\int_{0-}^t\co(\cR^C)dB)_{t\in[0,T]}$ is an $\hF$-submartingale, we have
\[
\int_{0-}^0 \co(\cR^C)dB\subseteq \hE\sqb{\int_{0-}^t \co(\cR^C)dB}\subseteq \hE\sqb{\int_{0-}^T \co(\cR^C)dB}
\]
so that, by the tower property of set-valued conditional expectation (cf. \cite[Theorem~2.1.74]{Molchanov}),
\[
\hE\sqb{\int_{0-}^t \co(\cR^C)dB}=\hE\sqb{\int_{0-}^T \co(\cR^C)dB}=\hE\sqb{\hE\sqb{\int_{0-}^T \co(\cR^C)dB\mid\cF_t}}
\]
for each $t\in[0,T]$. On the other hand, the submartingale property gives
\[
\int_{0-}^t \co(\cR^C)dB\subseteq \hE\sqb{\int_{0-}^T \co(\cR^C)dB\mid\cF_t}\quad\hP\text{-a.s.}
\]
By Lemma~\ref{lem:nonnegative}, we get $\int_{0-}^t \co(\cR^C)dB= \hE[\int_{0-}^T \co(\cR^C)dB|\cF_t]$ $\hP$-a.s. Hence, (ii) follows.

(iii)$\Leftrightarrow$(v), (iv)$\Leftrightarrow$(vi): Let $x^\ast\in\sW$ and $t\in[0,T]$. Note that the support function of a set and that of its convex hull are identical. Hence, by Proposition~\ref{lem:basicprop}(ii,iii), we have
\beaa
s\of{x^\ast, \int_{0-}^t \co(\cR^C)dB}&=& s\of{x^\ast,\int_{0-}^t \cR^C dB}= \sup\int_{0-}^t \ip{x^\ast,\cR^C}dB\\
&=& \sup_{i\in\hN}\cJ^0_t(\ip{x^\ast,x^i},\ip{x^\ast,z^i})=\sup_{i\in\hN}\ip{x^\ast,
	\hE[\xi^i|\cF_t]}
\eeaa 
with probability one. Using \cite[Theorem~2.1.35]{Molchanov}, we calculate the support function of the expectation as
\[
s\of{x^\ast,\hE\sqb{\int_{0-}^t \co(\cR^C) dB}}=\hE\sqb{s\of{x^\ast,\int_{0-}^t \co(\cR^C)dB}}=\hE\sqb{\sup_{i\in\hN}\ip{x^\ast,\hE[\xi^i|\cF_t}}.
\]
On the other hand, we have
\[
s\of{x^\ast,\cl_{\hR^d}\co(\{x^i\colon i\in\hN\})}=\sup_{i\in\hN}\ip{x^\ast,x^i}=\sup_{i\in\hN}\hE[\ip{x^\ast,\hE[\xi^i|\cF_t]}].
\]
We have $\hE[\int_{0-}^t \co(\cR^C)dB]=\cl_{\hR^d}\co(\{x^i\colon i\in\hN\})$ if and only if $s(x^\ast,\hE[\int_{0-}^t \co(\cR^C)dB])=s(x^\ast,\cl_{\hR^d}\co(\{x^i\colon i\in\hN\}))$ for every $x^\ast\in\sW$. Hence, (iii)$\Leftrightarrow$(v) and (iv)$\Leftrightarrow$(vi) follow.
\qed 

\section{Conditional Expectations of Random Polytopes}\label{sec:polytope}

In this section, we focus on the conditional expectation of a random polytope $\Xi$, i.e.,
\bea\label{eq:randompoly}
\Xi= \co(\{\xi^1,\ldots,\xi^n\})\quad \hP\text{-a.s.}
\eea
for some $\xi^1,\ldots,\xi^n\in \hL^2_{\cF_T}(\hR^d)$. For $\hP$-a.e. $\o\in\O$, the set of all vertices of $\Xi(\o)$ is a subset of $\{\xi^1(\o),\ldots,\xi^n(\o)\}$. We are interested in the vertices of the set-valued conditional expectation $\hE[\Xi|\cF_t](\o)$ for $t\in [0,T)$ and we raise the following question: does the set $\{\hE[\xi^1|\cF_t](\o),\ldots,$ $\hE[\xi^n|\cF_t](\o)\}$ cover all vertices of $\hE[\Xi|\cF_t](\o)$ as $t$ decreases from $T$ to $0$? In view of Theorem~\ref{prop:twointegrals}, the answer is negative in general. As a follow-up on Theorem~\ref{prop:twointegrals2}, we will provide a geometric characterization of the case where the answer is positive. Finally, we will provide a concrete example for the latter case.

\subsection{Normal Fans and Type Cones}\label{sec:polytopeprelim}

The geometric characterization announced above is related to the notion of a \emph{normal fan} for a polytope. For this reason, we first pause our discussion on set-valued stochastic integrals and review some preliminary notions from polyhedral theory. We refer the reader to \cite{normalcone,typecone1,typecone2} for more details about this subject.

Let $D\subseteq\hR^d$ be a set. It is called a \emph{cone} if $r x\in D$ for every $r>0$ and $x\in D$. The \emph{convex conic hull} of $D$ is the smallest superset of $D$ that is a convex cone and it is given by
\[
\cone(D)=\cb{\sum_{i=1}^n r_i x^i\colon r_1,\ldots,r_n\geq 0;\ x^1,\ldots,x^n\in D;\ n\in\hN}.
\]
The set $D$ is called \emph{affine} if $r x^1+(1-r)x^2\in D$ for every $r\in\hR$ and $x^1,x^2\in D$. The \emph{affine hull} of $D$ is the smallest affine superset of $D$ and it is given by
\[
\aff(D)=\cb{\sum_{i=1}^n r_i x^i \colon r_1,\ldots, r_n\in\hR;\, \sum_{i=1}^n r_i=1;\ x^1,\ldots, x^n\in D;\ n\in\hN}.
\]
The dimension of an affine set is defined as that of the subspace $\aff(D-\{x\})$, where $x\in D$. The set $D$ is called a \emph{halfspace} (resp. \emph{hyperplane}) if it is of the form $D=\{x\in\hR^d\colon \ip{x^\ast,x}\leq r\}$ (resp. $D=\{x\in\hR^d\colon\ip{x^\ast,x}=r\}$) for some $x^\ast\in\hR^d$ and $r\in\hR$. For a set $C\subseteq\hR^d$, we say that $D$ is a \emph{supporting hyperplane} of $C$ if $D=\{x\in\hR^d\colon \ip{x^\ast,x}=s(x^\ast,C)\}$ for some $x^\ast\in\hR^d$ and $C\cap D\neq \emptyset$.

The intersection of a finite number of halfspaces in $\hR^d$ is called a \emph{convex polyhedron (polyhedral set)} in $\hR^d$. Let $P$ be a convex polyhedron in $\hR^d$. The dimension of $P$ is defined as that of $\aff(S)$ and it is denoted by $\dim(P)$. If $D$ is a supporting hyperplane of $P$, then $P\cap D$ is called a \emph{face} of $P$. Clearly, every face of $P$ is also a convex polyhedron. A zero-dimensional face (or its unique element) is called a \emph{vertex}, a one-dimensional face is called an \emph{edge}, and a $(d-1)$-dimensional face is called a \emph{facet} of $P$. A face $F$ is called \emph{proper} if $F\neq P$. If $P$ is also a cone, then it can be written as $P=\cone(\{x^1,\ldots,x^k\})$ for some $x^1,\ldots,x^k\in P\setminus\{0\}$ with $k\in\hN$; in this case, $P$ is said to be \emph{generated by} these vectors. A convex polyhedral cone is called \emph{simplicial} if it is generated by linearly independent vectors.

Before proceeding further, we provide the definition of a central object for our discussions.

\begin{defn}\label{defn:fan}
\cite[Definition~1.1]{normalcone} Let $\sN$ be a collection of convex polyhedral cones in $\hR^d$. We say that $\sN$ is a \emph{fan} if it satisfies the following properties:

(i) If $C\in\sN$ and $F$ is a face of $C$, then $F\in\sN$.

(ii) If $C_1, C_2 \in \sN$, then $C_1\cap C_2$ is both a face of $C_1$ and a face of $C_2$.

In this case, $\sN$ is called \emph{complete} if $\bigcup_{C\in\sN}C=\hR^d$; it is called \emph{simplicial} if every $C\in\sN$ is simplicial; it is called \emph{essential} if $\{0\}\in\sN$. A cone in $\sN$ is called \emph{maximal} if it is not a proper subset of an element of $\sN$. Two maximal cones in $\sN$ are called \emph{adjacent} if their intersection is a facet of each.
\end{defn}

Let $F\subseteq P$ be a set. Then, the set
$\sN(F,P):=\{x^\ast\in\hR^d\colon \ip{x^\ast,x}=s(x^\ast,P)\text{ for each }x\in F\}$
is called the \emph{normal cone} of $P$ at $F$. If $F=\{v\}$ for some $v\in P$, then we simply write $\sN(v,P):=\sN(\{v\},P)$. If $F$ is a face of $P$, then $\sN(F,P)$ is a convex polyhedral cone and
\[
\dim(\sN(F,P))+\dim(F)=d.
\]
Moreover, if $F_1, F_2$ are faces of $P$ such that $F_1\subseteq F_2$, then $\sN(F_1,P)\supseteq \sN(F_2,P)$. The collection
$\sN^P:=\{\sN(F,P)\colon F\text{ is a proper face of }P\}$
is called the \emph{normal fan} of $P$. It can be checked that $\sN^P$ is a fan. Note that the maximal elements of $\sN^P$ are precisely the normal cones of $P$ at its vertices. Moreover, $\dim(P)=d$ if and only if $\sN^P$ is essential.

The next lemma will be useful later and should be well-known. Due to lack of a clear reference, we provide a proof.

\begin{lem}\label{lem:normalconeint}
	Let $P\subset\hR^d$ be a convex polyhedron and $v\in P$ be a vertex of it. For every $x^\ast\in \interior_{\hR^d}(\sN(v,P))$, the vertex $v$ is the unique maximizer of $\ip{x^\ast,\cdot}$ over $P$, i.e., for every $x\in P\setminus\{v\}$, we have $\ip{x^\ast,v}>\ip{x^\ast,x}$.
\end{lem}

{\it Proof}: To get a contradiction, suppose that there exist $x^\ast\in\interior_{\hR^d}(\sN(v,P))$ and $x\in P\setminus\{v\}$ such that $\ip{x^\ast,v}=\ip{x^\ast,x}=s(x^\ast,P)$. Note that $P$ itself is a face that contains both $v$ and $x$. Moreover, the intersection of faces of $P$ is again a face of $P$. Let $F$ be the intersection of all faces that contain both $v$ and $x$, which is the smallest such face. Then, $x^\ast\in \sN(F,P)$. Since $\dim(F)\geq 1$, we have $\dim(\sN(F,P))\leq d-1$. By Definition~\ref{defn:fan}, $\sN(F,P)\cap \sN(v,P)$ is a face of $\sN(v,P)$. Since $\sN(F,P)\cap \sN(v,P)$, $\sN(v,P)$ have different dimensions, they cannot be equal. Hence, $x^\ast\in \sN(F,P)\cap \sN(v,P)\subseteq \bd_{\hR^d}(\sN(v,P))$, a contradiction to $x^\ast\in\interior_{\hR^d}(\sN(v,P))$.
\qed

If $P$ is bounded, then it is called a \emph{polytope}. A set is a polytope if and only if it is the convex hull of a finite subset of $\hR^d$. In this case, every face of $P$ is also a polytope. Moreover, a convex polyhedron is a polytope if and only if its normal fan is complete.

We may consider two polytopes as equivalent if they have the same normal fan. This establishes an equivalence relation on the set of all polytopes, whose equivalence classes are given in the next definition.

\begin{defn}
	Let $\sN$ be a complete fan. The set of all polytopes whose normal fan is $\sN$ is called the \emph{type cone} of $\sN$.
\end{defn}

Let $\sN$ be a finite complete fan and let $\{x^{\ast,1},\ldots,x^{\ast,\ell}\}\subset\hR^d\setminus\{0\}$ be the set of all generating vectors of the members of $\sN$, where $\ell\in\hN$. For each $C\in\sN$, we denote by $J_C\subseteq\{1,\ldots,\ell\}$ the set of indices of the generating vectors of $C$. Let $X^\ast\in\hR^{\ell\times d}$ be the matrix whose rows are $x^{\ast,1},\ldots,x^{\ast,\ell}$. Then, every member of the type cone of $\sN$ can be written as
\bea\label{eq:ptype}
P_h:=\{x\in\hR^d\colon \ip{x^{\ast,j},x}\leq h_j\text{ for every }j\in\{1,\ldots,\ell\}\}=\{x\in\hR^d\colon X^\ast x\leq h\}
\eea
for some $h=(h_1,\ldots,h_\ell)\in\hR^\ell$. In this case, we indeed have
$h_j=s(x^{\ast,j},P_h)$, for every $j\in\{1,\ldots,\ell\}$. However, not every choice of $h\in\hR^\ell$ yields a member $P_h$ of the type cone. We say that $h\in\hR^d$ is \emph{admissible} if $P_h$ is in the type cone of $\sN$. Hence, we may index the members of the type cone by admissible vectors by defining
\[
\TC(\sN) := \{h\in\hR^\ell\mid \sN^{P_h}=\sN\}.
\]

Suppose that $\sN$ is complete, simplicial, and essential. Then, each maximal cone in $\sN$ is $d$-dimensional and it is generated by $d$ linearly independent vectors. Hence, if $C_1,C_2\in \sN$ are adjacent maximal cones, then there exist $j_1,j_2\in \{1,\ldots,\ell\}$ such that
$J_{C_1}\setminus\{j_1\}=J_{C_2}\setminus\{j_2\}$, 
and the set of generating vectors corresponding to $J_{C_1}\cup J_{C_2}=(J_{C_1}\cap J_{C_2})\cup\{j_1,j_2\}$ is linearly dependent; let $\alpha_{C_1,C_2}(j)\in\hR$, $j\in J_{C_1}\cup J_{C_2}$, be the unique coefficients such that
\bea\label{eq:lindep}
\sum_{j\in J_{C_1}\cup J_{C_2}}\alpha_{C_1,C_2}(j)x^{\ast,j}=0
\eea
and $\alpha_{C_1,C_2}(j_1)+\alpha_{C_1,C_2}(j_2)=2$. (The latter condition is there to ensure uniqueness.)

The next result provides an algebraic characterization of $\TC(\sN)$.

\begin{thm}\label{thm:admissible}
	\cite[Proposition~1.1]{typecone1} Let $\sN$ be a complete, simplicial, and essential fan. Let $h\in\hR^d$. Then, $h\in \TC(\sN)$ if and only if, for every choice of adjacent maximal cones $C_1,C_2\in\sN$, it holds that
$	\sum_{j\in J_{C_1}\cup J_{C_2}}\alpha_{C_1,C_2}(j)h_j>0$.
\end{thm}

\subsection{Random Polytopes with Deterministic Normal Fans}\label{sec:polytopemain}

We go back to our discussion on set-valued stochastic integrals and consider a random polytope $\Xi$ as given in \eqref{eq:randompoly}. With the notation of Section~\ref{sec:compare}, we have 
\[
\int_{0-}^t \cR^M dB = \hE[\Xi|\cF_t],\quad \int_{0-}^t \cR^C dB = \co(\{\hE[\xi^1|\cF_t],\ldots,\hE[\xi^n|\cF_t]\}),\quad t\in[0,T],
\]
with probability one, where $\cR^M$ is provided by the set-valued martingale representation theorem and $\cR^C=\{(x^1,z^1),\ldots,(x^n,z^n)\}$ is such that $\cJ^0_T(x^i,z^i)=\xi^i$ for each $i\in\{1,\ldots,n\}$.

Our aim is to show that the above set-valued stochastic integrals coincide if and only if $\Xi$ has a deterministic normal fan, up to a $\hP$-null set. Before stating the main result, we prove the following lemma of independent interest.

\begin{lem}\label{lem:volume}
	Let $C, D\in\sG(\hR^d)$ be with $\interior_{\hR^d}(C)\neq \emptyset$, $C\subseteq D$, and $\Leb(D\neg\setminus\neg C)=0$. Then, $C=D$.
	\end{lem}

{\it Proof}:
	We first claim that $\interior_{\hR^d}(C)=\interior_{\hR^d}(D)$. Since $C\subseteq D$, we have $\interior_{\hR^d}(C)\subseteq \interior_{\hR^d}(D)$. To show the reverse inclusion, it is enough to check that $\interior_{\hR^d}(D)\subseteq C$. Suppose that there exists $x\in\interior_{\hR^d}(D)$ with $x\notin C$, that is, $x\in \interior_{\hR^d}(D)\cap C^c$. Then, we can find $r>0$ such that $\{x\}+\hB_{\hR^d}(r)\subseteq \interior_{\hR^d}(D)\cap C^c\subseteq D\setminus C$. Since $\Leb(\{x\}+\hB_{\hR^d}(r))>0$, we get a contradiction. Therefore, we have $\interior_{\hR^d}(C)=\interior_{\hR^d}(D)$. Since $C, D$ are closed convex sets with nonempty interior, we conclude that $C=\cl_{\hR^d}(\interior_{\hR^d}(C))=\cl_{\hR^d}(\interior_{\hR^d}(D))=D$ by \cite[Lemma~5.13]{aliprantis}.
\qed 

\begin{thm}\label{cor:twointegrals2}
	The following are equivalent:
	
	(i) It holds $\int_{0-}^t \cR^M dB = \int_{0-}^t \co(\cR^C)dB$ for every $t\in[0,T]$ $\hP$-a.s.
	
	(ii) $(\int_{0-}^t \co(\cR^C)dB)_{t\in[0,T]}$ is an $\hF$-martingale.
	
	(iii) It holds $\hE[\int_{0-}^t \co(\cR^C)dB]=\co(\{x^1,\ldots,x^n\})$ for every $t\in[0,T]$.
	
	(iv) It holds $\hE[\int_{0-}^T \co(\cR^C)dB]=\co(\{x^1,\ldots,x^n\})$.
	
	(v) For every $x^\ast\in\hR^d$, $t\in[0,T]$, it holds $\hE[\max_{i\in\{1,\ldots,n\}}\ip{x^\ast,\hE[\xi^i|\cF_t]}]=\max_{i\in\{1,\ldots,n\}}\ip{x^\ast,x^i}$.
	
	(vi) For every $x^\ast\in\hR^d$, it holds $\hE[\max_{i\in\{1,\ldots,n\}}\ip{x^\ast,\xi^i}]=\max_{i\in\{1,\ldots,n\}}\ip{x^\ast,x^i}$.
	
	(vii) For every $x^\ast\in \hR^d$, there exists $i(x^\ast)\in\{1,\ldots,n\}$ such that it holds $\ip{x^\ast,\hE[\xi^{i(x^\ast)}|\cF_t]}\geq \ip{x^\ast,\hE[\xi^i|\cF_t]}$ $\hP$-a.s. for every $i\in\{1,\ldots,n\}$ and $t\in[0,T]$.
	
	(viii) For every $x^\ast\in \hR^d$, there exists $i(x^\ast)\in\{1,\ldots,n\}$ such that it holds $\ip{x^\ast,\xi^{i(x^\ast)}}\geq \ip{x^\ast,\xi^i}$ $\hP$-a.s. for every $i\in\{1,\ldots,n\}$.
	
	(ix) For each $i\in\{1,\ldots,n\}$ and $t\in[0,T]$, the normal cone of $\hE[\Xi|\cF_t]$ at $\hE[\xi^i|\cF_t]$ is $\hP$-a.s. deterministic and free of $t$, i.e., there exists a convex polyhedral cone $\sN_i\subseteq\hR^d$ such that, every $t\in[0,T]$, we have $\sN_i=\sN(\hE[\xi^i|\cF_t](\o),\hE[\Xi|\cF_t](\o))$ for $\hP$-a.e. $\o\in\O$.
	
	(x) For each $i\in\{1,\ldots,n\}$, the normal cone of $\Xi$ at $\xi^i$ is $\hP$-a.s. deterministic, i.e., there exists a convex polyhedral cone $\sN_i\subseteq\hR^d$ such that $\sN_i=\sN(\xi^i(\o),\Xi(\o))$ for $\hP$-a.e. $\o\in\O$.
	
	(xi) For each $t\in[0,T]$, $\hE[\Xi|\cF_t]$ has a $\hP$-a.s. deterministic normal fan that is free of $t$, i.e., for each nonempty $I\subseteq\{1,\ldots,n\}$, there exists a convex polyhedral cone $\sN_I\subseteq\hR^d$ such that, for every $t\in[0,T]$, we have $\sN_I=\sN(F_{I,t}(\o),\hE[\Xi|\cF_t](\o))$ for $\hP$-a.e. $\o\in\O$, where $F_{I,t}:=\co(\{\hE[\xi^i|\cF_t]\colon i\in I\})$.
	
	(xii) $\Xi$ has a $\hP$-a.s. deterministic normal fan, that is, for each nonempty $I\subseteq\{1,\ldots,n\}$, there exists a convex polyhedral cone $\sN_I\subseteq\hR^d$ such that $\sN_I=\sN(F_I(\o),\Xi(\o))$ for $\hP$-a.e. $\o\in\O$, where $F_I:=\co(\{\xi^i\colon i\in I\})$.  Moreover,  for $I\subseteq\{1,\ldots,n\}$ and	$i\in\{1,\ldots,n\}$, one can take
	\bea
	\label{eq:W(I)}
	\sN_I&=&\sN(\co(\{x^i\colon i\in I\}),\co(\{x^1,\ldots,x^n\})),\\
        \label{eq:W(i)}
	\sN_i&=&\sN_{\{i\}}=\sN(x^i,\co(\{x^1,\ldots,x^n\})). 
	\eea 
	\end{thm}
{\it Proof.} The equivalence of (i)--(vi) follows from Theorem~\ref{prop:twointegrals2}. The equivalence (vii)$\Leftrightarrow$(viii) and the implications (ix)$\Rightarrow$(x), (xi)$\Rightarrow$(xii), (xii)$\Rightarrow$(x) are trivial.

(vi)$\Leftrightarrow$(viii): Let $x^\ast\in\sW$. Let $i(x^\ast)\in\{1,\ldots,n\}$ be such that
$\langle x^\ast,  x^{i(x^\ast)}\rangle =\max_{i\in\{1,\ldots,n\}}\langle x^\ast, x^i \rangle$,
that is, $\hE[\langle x^\ast, \xi^{i(x^\ast)} \rangle]=\max_{i\in\{1,\ldots,n\}}\hE[\langle x^\ast, \xi^i \rangle]$.
Since we have $\langle x^\ast, \xi^{i(x^\ast)} \rangle\leq \max_{i\in\{1,\ldots,n\}}\langle x^\ast, \xi^i \rangle$ $\hP$-a.s.,
these random variables have the same expectation if and only if they are equal $\hP$-a.s. Therefore, (vi) is equivalent to (viii).

(vii)$\Rightarrow$(ix): For each $i\in\{1,\ldots,n\}$, let $\sN_i$ be defined by \eqref{eq:W(i)}, let $\sW_i$ denote a countable dense subset of $\sN_i\cap \hB_{\hR^d}$. Note that $\co(\{x^1,\ldots,x^n\})$ is a convex polytope, the set $\{x^1,\ldots,x^n\}$ contains all vertices of it, and the normal cones of these vertices cover $\hR^d$. Hence, we have $\bigcup_{i=1}^n \sN_i=\hR^d$. Moreover, $\bigcup_{i=1}^n \sW_i$ is a countable dense subset of $\hB_{\hR^d}$.

Let us fix $i\in\{1,\ldots,n\}$. Let $x^\ast\in\sW_i$. By (vii), there exists $i(x^\ast)\in\{1,\ldots,n\}$ such that $\ip{x^\ast,\hE[\xi^{i(x^\ast)}|\cF_t]}=\max_{j\in\{1,\ldots,n\}} \ip{x^\ast,\hE[\xi^j|\cF_t]}$ $\hP$-a.s. for every $t\in[0,T]$. In particular, we have $\ip{x^\ast,\hE[\xi^i|\cF_t]}\leq \ip{x^\ast,\hE[\xi^{i(x^\ast)}|\cF_t]}$ $\hP$-a.s. and
\beaa 
\hE[\ip{x^\ast,\hE[\xi^i|\cF_t]}]&=&\ip{x^\ast,x^i}=\max_{i^\prime\in\{1,\ldots,n\}}\lan x^\ast,x^{i^\prime}\ran\\
&=& \hE\sqb{\max_{i^\prime\in\{1,\ldots,n\}} \ip{x^\ast,\hE[\xi^{i^\prime}|\cF_t]}}=\hE\sqb{\ip{x^\ast,\hE[\xi^{i(x^\ast)}|\cF_t]}}
\eeaa
for every $t\in[0,T]$. Here, the third equality is by (vi), which is already shown to be equivalent to (viii). Hence, we conclude that $\ip{x^\ast,\hE[\xi^i|\cF_t]}=\ip{x^\ast,\hE[\xi^{i(x^\ast)}|\cF_t]}$ $\hP$-a.s., that is,
\[
\ip{x^\ast,\hE[\xi^i|\cF_t]}=\max_{i^\prime\in\{1,\ldots,n\}} \ip{x^\ast,\hE[\xi^{i^\prime}|\cF_t]}=s(x^\ast,\hE[\Xi|\cF_t])\quad \hP\text{-a.s.}
\]
for each $t\in[0,T]$. Here, the second equality is by (i), which is already shown to be equivalent to (vii). Let $t\in[0,T]$. Since $x^\ast\in\sW_i$ is arbitrary, $\sW_i$ is dense in $\sN_i\cap\hB_{\hR^d}$, and $s(\cdot,\hE[\Xi|\cF_t])$ is continuous and positively homogeneous, we conclude that $\ip{x^\ast,\hE[\xi^i|\cF_t]}=s(x^\ast,\hE[\Xi|\cF_t])$ for every $x^\ast\in \sN_i$ on a set $\O_{i,t}\in\cF_t$ with probability one. In particular, for every $\o\in\O_{i,t}$, we have
\bea\label{eq:normalsubset}
\sN_i\subseteq \sN(\hE[\xi^i|\cF_t](\o),\hE[\Xi|\cF_t](\o)).
\eea 

Without loss of generality, suppose that the vertices of $\co(\{x^1,\ldots,x^n\})$ are $x^1,\ldots,x^k$, where $k\in\{1,\ldots,n\}$. Then, we have $\dim(\{x^i\})=0$ so that $\dim(\sN_i)=d$ and $\Leb(\sN_i)>0$ for each $i\in\{1,\ldots,k\}$. On the other hand, for each $i\in\{k+1,\ldots,n\}$, $x^i$ belongs to a nonzero-dimensional face of $\co(\{x^1,\ldots,x^n\})$ whose normal cone is $\sN_i$ so that we have $\dim(\sN_i)\leq d-1$ and $\Leb(\sN_i)=0$. Let $\o\in\bigcap_{i=1}^k\O_{i,t}$. Since $\bigcup_{i=1}^k \sN_i=\hR^d$, by \eqref{eq:normalsubset}, we get
\bea\label{eq:cover}
\bigcup_{i=1}^k (\sN_i\cap \hB_{\hR^d}) = \hB_{\hR^d}=\bigcup_{i=1}^k \sN(\hE[\xi^i|\cF_t](\o),\hE[\Xi|\cF_t](\o)).
\eea
Moreover, for each $i_1,i_2\in\{1,\ldots,k\}$ with $i_1\neq i_2$, we have
\[
\Leb(\sN(\hE[\xi^{i_1}|\cF_t](\o),\hE[\Xi|\cF_t](\o))\cap \sN(\hE[\xi^{i_1}|\cF_t](\o),\hE[\Xi|\cF_t](\o)))=0
\]
since the set is a face of $\sN(\hE[\xi^{i_1}|\cF_t](\o),\hE[\Xi|\cF_t](\o))$. By \eqref{eq:normalsubset}, we also have $\Leb(\sN_{i_1}\cap\sN_{i_2})=0$. Hence, we get
\beaa
& &\sum_{i=1}^k\Leb\of{\of{\sN(\hE[\xi^i|\cF_t](\o),\hE[\Xi|\cF_t](\o))\cap \hB_{\hR^d}}\setminus\of{\sN_i\cap \hB_{\hR^d}}}\\
&&=\sum_{i=1}^k\of{\Leb\of{\sN(\hE[\xi^i|\cF_t](\o),\hE[\Xi|\cF_t](\o))\cap \hB_{\hR^d}}-\Leb{\of{\sN_i\cap \hB_{\hR^d}}}}\\
&&=\Leb\of{\bigcup_{i=1}^k\of{\sN(\hE[\xi^i|\cF_t](\o),\hE[\Xi|\cF_t](\o))\cap \hB_{\hR^d}}}-\Leb\of{\bigcup_{i=1}^k \of{\sN_i\cap \hB_{\hR^d}}}=0,
\eeaa 
where the last equality is by \eqref{eq:cover}. Hence,
\[
\Leb\of{\of{\sN(\hE[\xi^i|\cF_t](\o),\hE[\Xi|\cF_t](\o))\cap \hB_{\hR^d}}\setminus\of{\sN_i\cap \hB_{\hR^d}}}=0
\]
for each $i\in\{1,\ldots,k\}$. Since $\interior_{\hR^d}(\sN_i\cap \hB_{\hR^d})\neq\emptyset$ as well, by Lemma~\ref{lem:volume}, we obtain that
\[
\sN_i\cap \hB_{\hR^d}=\sN(\hE[\xi^i|\cF_t](\o),\hE[\Xi|\cF_t](\o))\cap \hB_{\hR^d}
\]
and therefore
\bea\label{eq:equalityvert}
\sN_i=\{x^\ast\in\hR^d\colon \ip{x^\ast,\hE[\xi^i|\cF_t](\o)}=s(x^\ast,\hE[\Xi|\cF_t](\o))\}
\eea
for each $i\in\{1,\ldots,k\}$ since these sets are cones.

Since $\interior_{\hR^d}(\sN_i)=\interior_{\hR^d}(\sN(\hE[\xi^i|\cF_t](\o),\hE[\Xi|\cF_t](\o)))\neq\emptyset$, the point $\hE[\xi^i|\cF_t](\o)$ is a vertex of $\hE[\Xi|\cF_t](\o)$ for each $i\in\{1,\ldots,k\}$. We claim that these are all the vertices of $\hE[\Xi|\cF_t](\o)$. To get a contradiction, suppose that $\hE[\xi^i|\cF_t](\o)$ is a vertex of $\hE[\Xi|\cF_t](\o)$ for some $i\in\{k+1,\ldots,n\}$. Without loss of generality, we assume that $i=k+1$. Then,
\[
\interior_{\hR^d}(\sN(\hE[\xi^{k+1}|\cF_t](\o),\hE[\Xi|\cF_t](\o)))\neq\emptyset.
\]
Let $x^\ast\in\hB_{\hR^d}$ be an element of this set. Then, by Lemma~\ref{lem:normalconeint}, $\hE[\xi^{k+1}|\cF_t](\o)$ is the unique maximizer of the function $\ip{x^\ast,\cdot}$. In particular, $\ip{x^\ast,\hE[\xi^i|\cF_t](\o)}<s(x^\ast,\hE[\Xi|\cF_t](\o))$ for each $i\in\{1,\ldots,k\}$. On the other hand, by \eqref{eq:cover}, there exists $i^{\prime}\in\{1,\ldots,k\}$ such that $\langle x^\ast,\hE[\xi^{i^\prime}|\cF_t](\o)\rangle =s(x^\ast,\hE[\Xi|\cF_t](\o))$, which is a contradiction. Therefore, the vertices of $\hE[\Xi|\cF_t](\o)$ are precisely $\hE[\xi^1|\cF_t](\o), \ldots, \hE[\xi^n|\cF_t](\o)$.

Let $i\in\{k+1,\ldots,n\}$. Then, $x^i$ belongs to a nonzero-dimensional face of $\co(\{x^1,\ldots,x^n\})$. Since the normal cone of a face of a polytope is determined by the normal cones of the vertices, by \eqref{eq:equalityvert}, we conclude that $\sN_i=\sN(\hE[\xi^i|\cF_t](\o),\hE[\Xi|\cF_t](\o))$ for each $i\in\{k+1,\ldots,n\}$ as well.

(ix)$\Rightarrow$(xi), (x)$\Rightarrow$(xii): These are immediate consequences of the fact that the normal cone of a face of a polytope is determined by the normal cones of its vertices, similar to the argument in the previous paragraph.

(x)$\Rightarrow$(viii): Let $x^\ast\in\hR^d$ and let $\O_0$ be the $\hP$-a.s. set in (x). Let $\o\in\O_0$. Since the vertices of $\Xi(\o)$ form a subset of $\{\xi^1(\o),\ldots,\xi^n(\o)\}$, we have $\bigcup_{i=1}^n \sN_i=\hR^d$. Hence, $x^\ast\in\sN_{i(x^\ast)}$ for some $i(x^\ast)\in\{1,\ldots,n\}$. Then, 
$
\lan x^\ast,\xi^{i(x^\ast)}\ran=s(x^\ast,\Xi(\o))=\max_{i\in\{1,\ldots,n\}}\ip{x^\ast,\xi^i(\o)}$, proving (viii).
\qed

Theorem~\ref{cor:twointegrals2} implies that the indefinite stochastic integral ($\int_{0-}^t \co(\cR^C)dB)_{t\in[0,T]}$, which is the convex hull of the (vector-valued) conditional expectations of $\xi^1,\ldots,\xi^n$, is an $\hF$-martingale if and only if the set-valued random variable $\Xi$ has a $\hP$-a.s. deterministic normal fan, which is the normal fan of $\int_{0-}^0 \co(\cR^C)dB= \co(\{x^1,\ldots,x^n\})$. Let $\sN$ denote this normal fan. We denote its elements by $\{x^{\ast,1},\ldots,x^{\ast,\ell}\}\subset\hR^d\setminus\{0\}$. The quantities $X^\ast$, $P_h$, $J_C$, $\alpha_{C_1,C_2}(j)$ are defined as in Section~\ref{sec:polytopeprelim}. Using Theorem~\ref{thm:admissible}, we may obtain an algebraic characterization of the deterministic normal fan condition as given by the next corollary.

\begin{cor}\label{cor:algebraic}
	Let $\sN$ be the normal fan of $\co(\{x^1,\ldots,x^n\})$. Suppose that $\sN$ is simplicial and essential. Then, the conditions in Theorem~\ref{cor:twointegrals2} are also equivalent to the following:
	
	(xiii) For every choice of adjacent maximal cones $C_1,C_2\in\sN$, it holds
	\[
	\sum_{j\in J_{C_1}\cup J_{C_2}}\alpha_{C_1,C_2}(j)\max_{i\in\{1,\ldots,n\}}\ip{x^{\ast,j},\xi^i}>0\quad \hP\text{-a.s.}
	\]
	\end{cor}
	
{\it Proof}: For $\hP$-a.e. $\o\in\O$, we may write $\Xi(\o)=P_{\eta(\o)}=\{x\in\hR^d\colon X^\ast x \leq \eta(\o)\}$, where $\eta_j(\o)=s(x^{\ast,j},\Xi(\o))=\max_{i\in\{1,\ldots,n\}}\ip{x^{\ast,j},\xi^i(\o)}$ for each $j\in\{1,\ldots,\ell\}$. Then, by Theorem~\ref{thm:admissible}, condition (xii) in Theorem~\ref{cor:twointegrals2} is equivalent to condition (xiii) above.
\qed

\begin{rem}\label{cor:algebraic2}
	Without reference to the random vectors $\xi^1,\ldots,\xi^n$, Corollary~\ref{cor:algebraic} can be restated as follows: ``Let $\sN$ be a complete, simplicial, and essential fan and $\Xi$ be a random polytope. The normal fan of $\Xi$ equals $\sN$ $\hP$-a.s. if and only if there exists $\eta\in\hL^2_{\cF_T}(\hR^\ell)$ such that $\Xi=P_{\eta}$ $\hP$-a.s. and for every choice of adjacent maximal cones $C_1,C_2\in\sN$, it holds
	\[
	\sum_{j\in J_{C_1}\cup J_{C_2}}\alpha_{C_1,C_2}(j)\eta_j>0\quad \hP\text{-a.s.}"
	\]
	These inequalities define a (deterministic) open convex cone $\TC(\sN)\subseteq \hR^\ell$ and the above condition says that $\eta$ is a square-integrable selection of $\TC(\sN)$, i.e., $\hP\{\eta\in \TC(\sN)\}=1$. Since $\TC(\sN)$ has a nonempty closed subset, by the standard measurable selection theorem, it has at least one measurable selection $\eta\in \hL^2_{\cF_T}(\hR^\ell)$.
	\qed
	\end{rem}
	
\begin{eg}\label{ex:triangle}
	We construct a simple example of a random polytope with deterministic normal fan. Suppose that $d=2$, $n=\ell=3$, $x^{\ast,1}=(-1,0)$, $x^{\ast,2}=(0,-1)$, $x^{\ast,3}=(1,1)$. The maximal cones of the corresponding normal fan $\sN$ have index sets $J_{C_1}=\{1,2\}$, $J_{C_2}=\{2,3\}$, $J_{C_3}=\{1,3\}$; see Figure~\ref{fig}(a). Then, the linear dependence system given in \eqref{eq:lindep} reads as
	\beaa
	& -\alpha_{C_1,C_2}(1)+\alpha_{C_1,C_2}(3)=0,\quad 
	 -\alpha_{C_1,C_2}(2)+\alpha_{C_1,C_2}(3)=0,
	& \alpha_{C_1,C_2}(1)+\alpha_{C_1,C_2}(3)=2,\\
	& -\alpha_{C_2,C_3}(1)+\alpha_{C_2,C_3}(3)=0,\quad
	-\alpha_{C_2,C_3}(2)+\alpha_{C_2,C_3}(3)=0,
	 & \alpha_{C_2,C_3}(2)+\alpha_{C_2,C_3}(1)=2,\\
	& -\alpha_{C_1,C_3}(1)+\alpha_{C_1,C_3}(3)=0,\quad
	-\alpha_{C_1,C_3}(2)+\alpha_{C_1,C_3}(3)=0,
	& \alpha_{C_1,C_3}(2)+\alpha_{C_1,C_3}(3)=2,
	\eeaa
	whose solution is given by $\alpha_{C_1,C_2}(j)=\alpha_{C_2,C_3}(j)=\alpha_{C_1,C_3}(j)=1$ for each $j\in\{1,2,3\}$. Hence, the type cone is given by
$	\TC(\sN)=\{h\in\hR^3\colon h_1+h_2+h_3>0\}$.
	In view of Remark~\ref{cor:algebraic2}, we conclude that a random polytope $\Xi$ has $\sN$ as its normal fan $\hP$-a.s. if and only if it is a random triangle of the form
$	\Xi = \{x\in\hR^2\colon -x_1\leq \eta_1,\ -x_2\leq \eta_2,\ x_1+x_2\leq \eta_3\}$,
	or equivalently,
	\[
	\Xi = \co(\{\xi^1=(-\eta_1,-\eta_2),\xi^2=(\eta_2+\eta_3,-\eta_2),\xi^3=(-\eta_1,\eta_1+\eta_3)\}),
	\] 
	where $\eta_1,\eta_2,\eta_3\in\hL^2_{\cF_T}(\O,\hR)$ are random variables such that $\eta_1+\eta_2+\eta_3>0$ $\hP$-a.s.
	
	As a special case, let us further assume that $m=3$, $T=1$, and take the conditional expectation processes as exponential martingales of Brownian motion, say, $\hE[\eta_i |\cF_t]=e^{\alpha B^i_t -\frac12 \alpha^2 t}$ for $i\in\{1,2,3\}$, $t\in[0,1]$, where $\alpha=\frac12$. We simulate the corresponding random triangle process by a random walk approximation of the Brownian motion with $10^4$ steps, see Figure~\ref{fig}(b). We also provide two two-dimensional plots that together carry the same information. In Figure~\ref{fig}(c), we show the sample path of the vertex corresponding to the right angle of the triangle as a curve parametrized by $t\in[0,1]$. Figure~\ref{fig}(d) shows the sample path of the length of the hypotenuse.
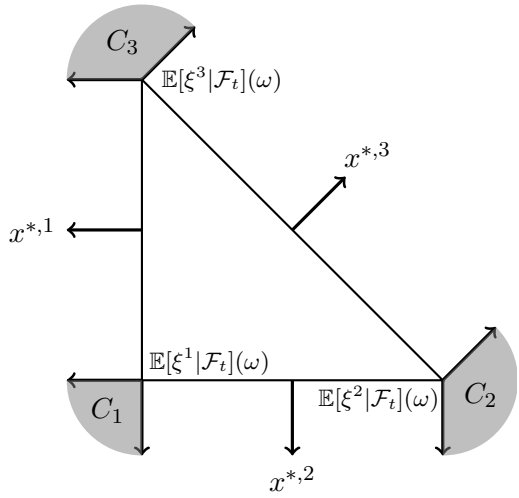
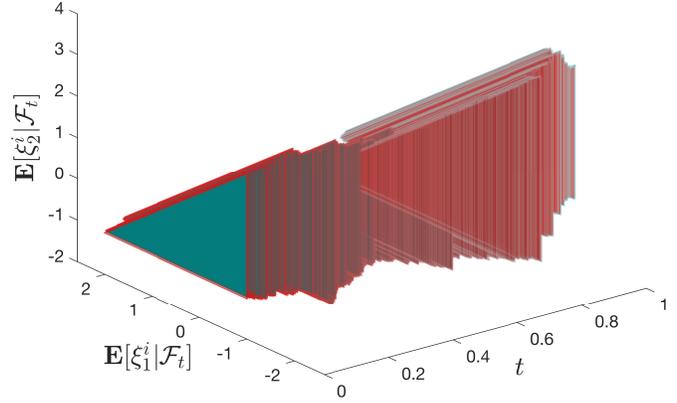
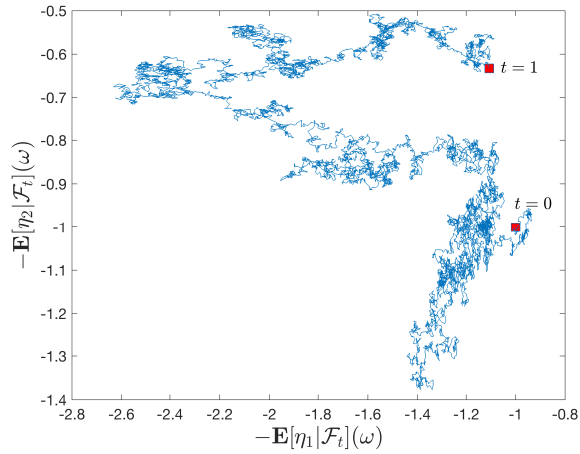
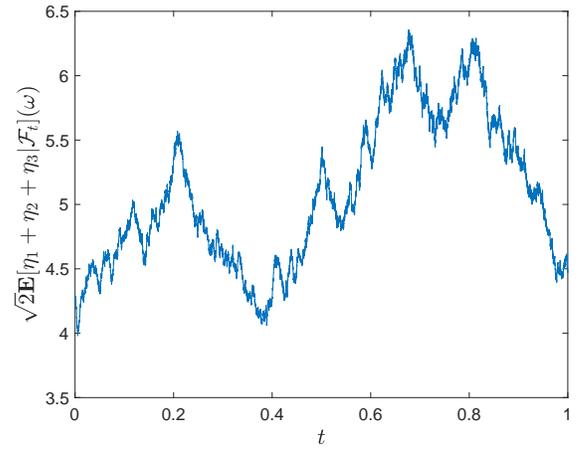
\begin{figure}[!htbp]
	\captionsetup[subfigure]{justification=centering}
	\centering
		\subfloat[Realization of the random triangle at an arbitrary time $t$ and its deterministic normal fan]{\begin{tikzpicture}
			\draw[thick] (0,0) -- (4,0) -- (0,4) -- cycle;
			\draw[very thick][->] (0,0) -- (0,-1);
			\draw[very thick][->] (0,0) -- (-1,0);
			\draw[very thick][->] (4,0) -- (4,-1);
			\draw[very thick][->] (4,0) -- (4.707,0.707);
			\draw[very thick][->] (0,4) -- (0.707,4.707);
			\draw[very thick][->] (0,4) -- (-1,4);
			\draw[very thick][->] (0,2) -- (-1,2) node[anchor=east] {$x^{\ast,1}$};
			\draw[very thick][->] (2,0) -- (2,-1) node[anchor=north] {$x^{\ast,2}$};
			\draw[very thick][->] (2,2) -- (2.707,2.707);
			\draw (3,3) node[anchor=center] {$x^{\ast,3}$};
			\fill[opacity=0.5,fill=gray] (0,0) -- (-1,0) arc (180:270:1);
			\fill[opacity=0.5,fill=gray] (0,4) -- (0.707,4.707) arc (45:180:1);
			\fill[opacity=0.5,fill=gray] (4,0) -- (4,-1) arc (-90:45:1);
			\draw (-0.1,-0.4) node[anchor=east] {$C_1$};
			\draw (4.15,-0.2) node[anchor=west] {$C_2$};
			\draw (0.1,4.5) node[anchor=east] {$C_3$};
			\draw (0.9,0.25) node[anchor=center,font=\footnotesize] {$\hE[\xi^1|\cF_t](\omega)$};
			\draw (2,4) node[anchor=east,font=\footnotesize] {$\hE[\xi^3|\cF_t](\omega)$};
			\draw (3.15,-0.25) node[anchor=center,font=\footnotesize] {$\hE[\xi^2|\cF_t](\omega)$};
			\end{tikzpicture}}
	\subfloat[Sample path of the random triangle process]{\includegraphics[width=0.6\textwidth]{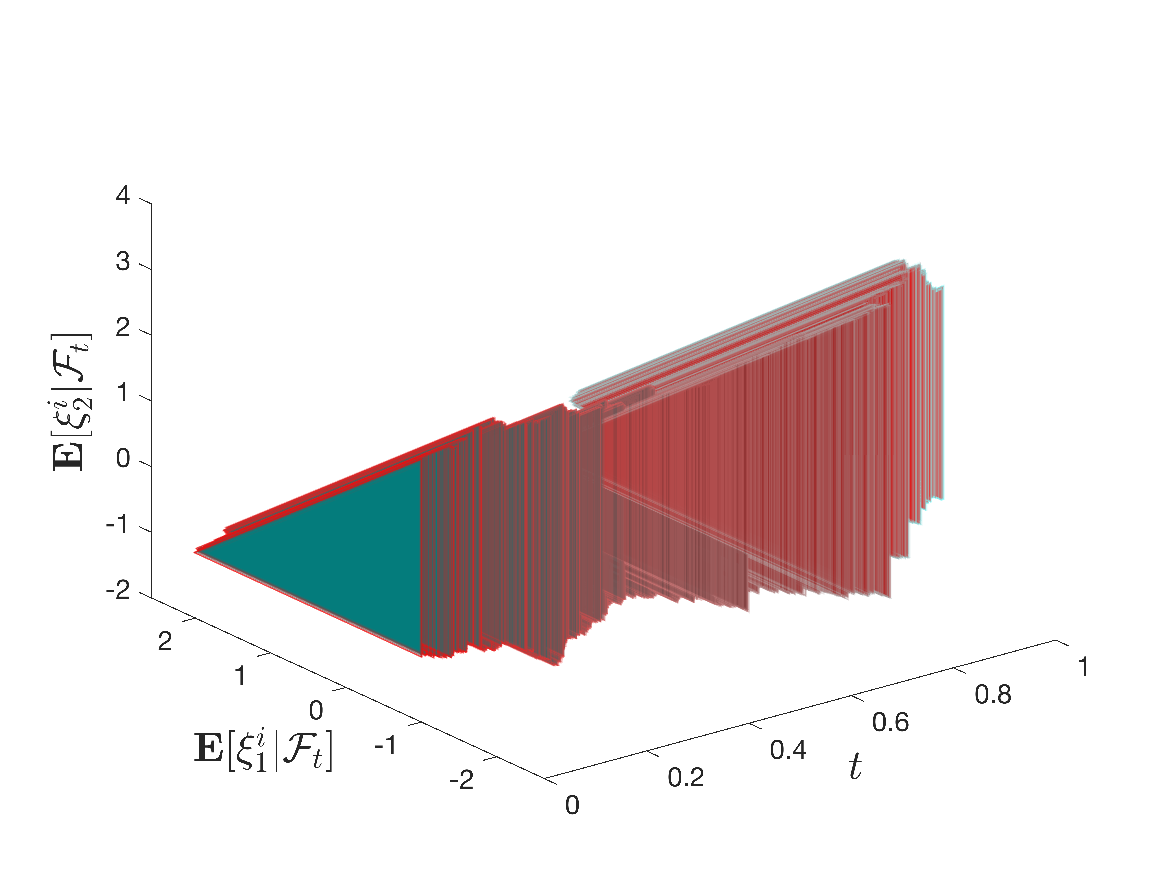}}\\
	\subfloat[Sample path of $(\hE\sqb{\xi^1|\cF_t})_{t\in\sqb{0,1}}$ as a parametric curve]{\includegraphics[width=0.52\textwidth]{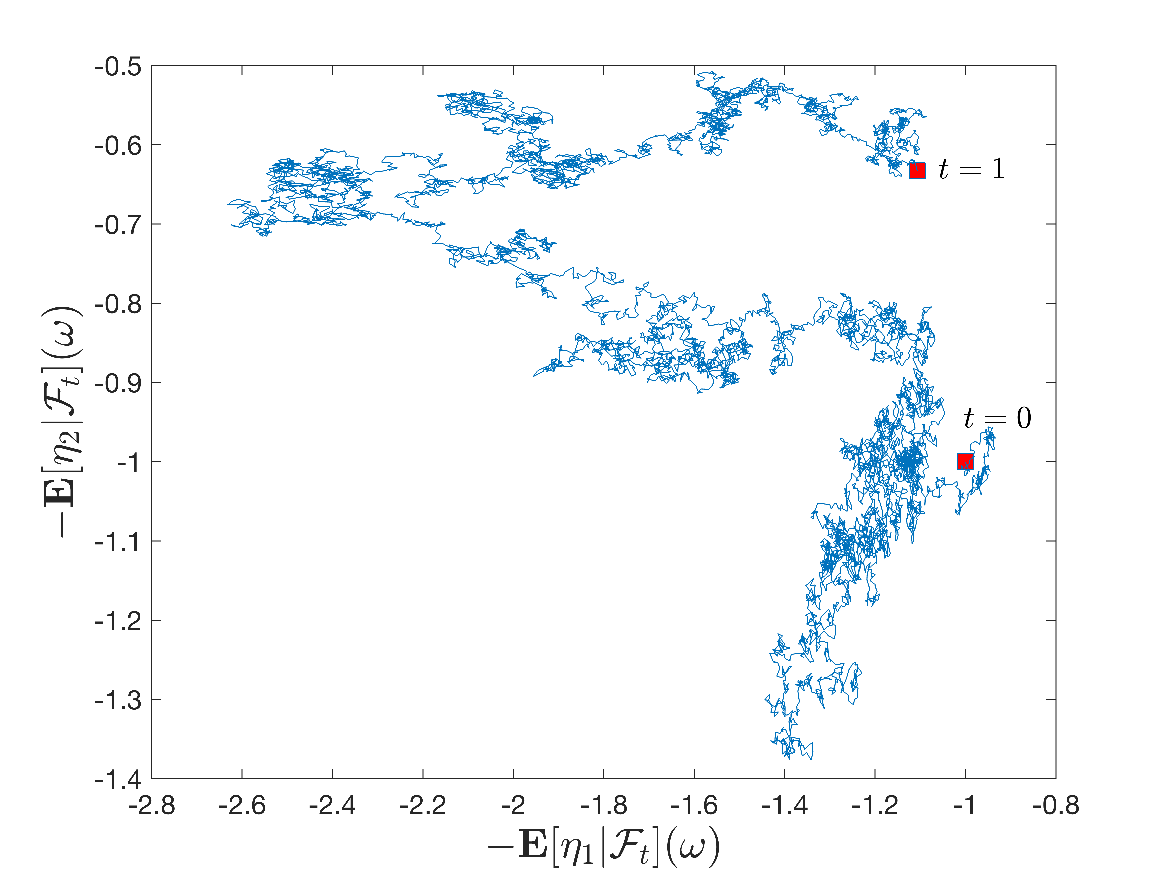}}
	\subfloat[Sample path of the length of the hypotenuse]{\includegraphics[width=0.52\textwidth]{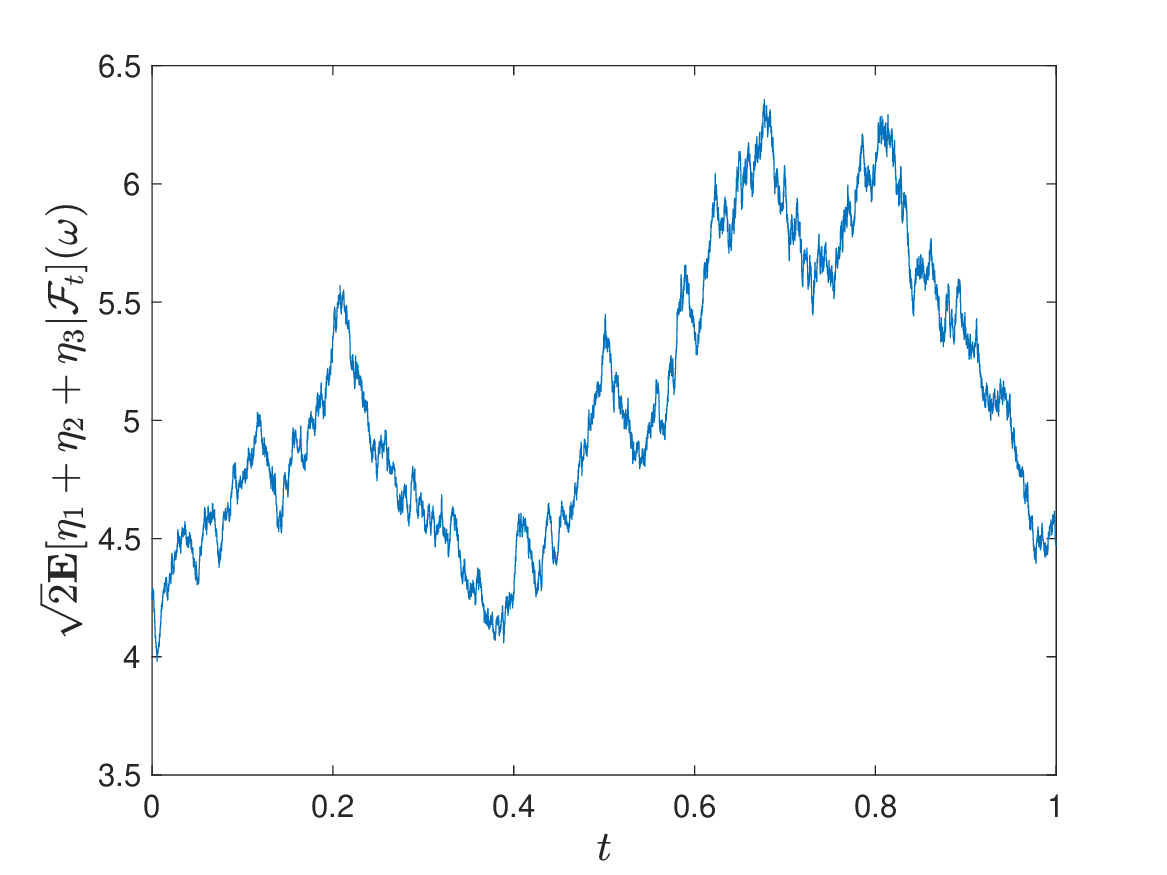}}
	\caption{Figures for Example~\ref{ex:triangle}.}
	\label{fig}
\end{figure}
	\end{eg}

%\no{\Large \bf Acknowledgments}

%\ms

\bibliographystyle{plainnat.bst}
\bibliography{AMbibliography.bib}
 \end{document}